\documentclass[onecolumn,draftcls,10pt]{IEEEtran}
\usepackage{cite}

\usepackage{color, colortbl}
\usepackage[x11names]{xcolor} 

\usepackage[pdftex]{graphicx}
\graphicspath{{fig/}{jpeg/}}
\usepackage[cmex10]{amsmath}
\usepackage{amssymb}
\usepackage{epstopdf}
\usepackage{multirow}
\usepackage{pifont}
\newcommand{\cmark}{\ding{51}}%
\newcommand{\xmark}{\ding{55}}%
\usepackage{algorithm}
\usepackage{tikz}
\usetikzlibrary{shapes,arrows,matrix}
\usetikzlibrary{decorations.pathmorphing} 
\usetikzlibrary{fit}					
\usetikzlibrary{backgrounds}	
\usepackage{verbatim}
\usepackage{algpseudocode}
\usepackage{amsmath,amssymb,lipsum}
\usepackage{hyperref}
\usepackage{comment}
\usepackage{graphicx}
\usepackage[font=small]{caption}
\usepackage{subcaption}
\usepackage{mathtools,enumerate}
\input{mysymbol.sty}
\usepackage{needspace}

\usepackage{theorem}
\newtheorem{thm}{Theorem}
\newtheorem{lemma}{Lemma}

\newtheorem{corollary}{Corollary}

\theoremstyle{definition}
\newtheorem{assumption}{Assumption}

\def \QED {$\blacksquare$}
\def \b {{\mathbf{b}}}
\def \W {{\mathbf{W}}}

\def \w {{\mathbf{w}}}
\def \x {{\mathbf{x}}}

\def \u {{\mathbf{u}}}
\def \y {{\mathbf{y}}}

 \providecommand{\Ex}[1]{\mathbb{E}\left[#1\right]}
 \providecommand{\Ek}[1]{\mathbb{E}_k\left[#1\right]} 
 \providecommand{\En}[1]{\mathbb{E}_k\|#1\|^2} 	
 \providecommand{\abs}[1]{\left|#1\right|}
 \providecommand{\norm}[1]{\left\|#1\right\|}
 \providecommand{\ip}[1]{\boldsymbol{\langle}#1\boldsymbol{\rangle}}
 \providecommand{\vect}[1]{\text{vec}\left(#1\right)}

    \def \eps {\epsilon}
    \def \mbE {{\mathbb{E}}}
    \def \mbR {{\mathbb{R}}}
    \def \bbx {{\mathbf{x}}}
    
    \def \bbz {{\mathbf{z}}}
    \def \ccalX {{\mathcal{X}}}
    \def \ccalC {{\mathcal{C}}}

    \def \bbs {{\mathbf{s}}}
        \def \bbv {{\mathbf{v}}}
    \def \bbG {{\mathbf{G}}}

    \def \bbtheta {{\boldsymbol{\theta}}}
    \def \bblambda {{\boldsymbol{\lambda}}}
        
    \def \Rn {{\mathbb{R}}}
    
    \def \lam {{\bblambda}}
    
        \def \v {{\bbv}}
    \def \G {{\bbG}}
    
    \def \z {{\bbz}}
    \def \r {{\bbr}}
    \def \cX {{\ccalX}}
    \def \cc {{\ccalC}}
    \def \cD {{\mathcal{D}}}
    \def \O {{\mathcal{O}}}

    \def \M {{\mathbf{M}}}
    \def \X {{\mathbf{X}}}

    \def \sX {{\mathsf{X}}}
    \def \sY {{\mathsf{Y}}} 
    \def \pd {{\mathcal{P}}}
    \def \md {{\text{MOST-FW }}}
    \def \ms {{\text{MOST-FW}$^+$ }}
    \def \ps {{\mathcal{P}^{+}}}
    \def \pdm {{\mathcal{P}_{\mu}}}
    
    \def \T {{\mathsf{T}}}
    \def \ssX {{\mathsf{\tilde{X}}}}
    \def \ssY {{\mathsf{\tilde{Y}}}} 
        \def \one {{\mathbf{1}}}
        
    \def \s {{\bbs}}

    \def \nf {{\nabla F_{\mu_k}(\x_k)}}

    \def \nff {{\nabla \hat{F}_{\mu_k}(\x_k)}}

    \def \n {\nabla}
    \def \nt {{\tilde{\nabla}}}

    \def \ffb {{\nabla \hat{F}_{\mu_k}(\x_k)}}

    \def \ffgbz {{\nabla \tilde{F}_{\mu_k}(\x_k)}}

    \def \Ec {{\mathcal{E}}}
    
  \tikzstyle{agent}=[circle,
  thick,
  minimum size=.5cm,
  draw=blue!80,
  fill=blue!20]

  \tikzstyle{neighbor}=[circle,
  thick,
  minimum size=.5cm,
  draw=cyan!85!black,
  fill=cyan!40,
  decorate,
  decoration={random steps,
  	segment length=2pt,
  	amplitude=2pt}]
  
  \tikzstyle{local_nat}=[rectangle,
  thick,
  minimum size=.5cm,
  draw=red!80,
  fill=red!20]
  
  \tikzstyle{glob_nat}=[rectangle,
  thick,
  minimum size=.5cm,
  draw=red!100,
  fill=red!40]      
  
  \tikzstyle{background}=[rectangle,
  fill=gray!10,
  inner sep=0.2cm,
  rounded corners=5mm]
  
  \tikzstyle{background2}=[rectangle,
  fill=green!20,
  inner sep=0.2cm,
  rounded corners=5mm]
{\tiny }
\title{\vspace{-0cm}{Zeroth and First Order Stochastic Frank-Wolfe Algorithms for Constrained Optimization}}

\author{Zeeshan~Akhtar,
	and~Ketan~Rajawat
	 \thanks{
	        Zeeshan Akhtar and K. Rajawat are with the Department of Electrical Engineering,
	        Indian Institute of Technology Kanpur, Kanpur 208016, India (e-mail: zeeshan@iitk.ac.in,
	        ketan@iitk.ac.in). }}
\begin{document}
\maketitle

\begin{abstract}
 This paper considers stochastic convex optimization problems with two sets of constraints: (a) deterministic constraints on the domain of the optimization variable, which are difficult to project onto; and (b) deterministic or stochastic constraints that admit efficient projection. Problems of this form arise frequently in the context of semidefinite programming as well as when various NP-hard problems are solved approximately via semidefinite relaxation. Since projection onto the first set of constraints is difficult, it becomes necessary to explore projection-free algorithms, such as the stochastic Frank-Wolfe (FW) algorithm. On the other hand, the second set of constraints cannot be handled in the same way, and must be incorporated as an indicator function within the objective function, thereby complicating the application of FW methods. Similar problems have been studied before; however, they suffer from slow convergence rates. This work, equipped with momentum based gradient tracking technique, guarantees fast convergence rates on par with the best-known rates for problems without the second set of constraints. Zeroth-order variants of the proposed algorithms are also developed and again improve upon the state-of-the-art rate results. We further propose the novel trimmed FW variants that enjoy the same convergence rates as their classical counterparts, but are empirically shown to require significantly fewer calls to the linear minimization oracle speeding up the overall algorithm. The efficacy of the proposed algorithms is tested on relevant applications of sparse matrix estimation, clustering via semidefinite relaxation, and uniform sparsest cut problem. 
\end{abstract}

\section{Introduction}\label{sec:intro}
We consider problems of the form
\begin{table}[h!]\label{prob_table}
	\centering
	\scalebox{1.2}{
		\begin{tabular}{c|c}
			\colorbox{MistyRose1}{$\pd$}  & \colorbox{LightCyan1}{$\ps$}\\ 
			\hline \\
			\begin{tabular}[c]{@{}c@{}}\vspace{1.5mm} $\min_{\bbx \in \cc}  \mbE[f(\bbx,\xi)]$\\ $\text{s.t.}$ $\G\x \in \cX$\end{tabular} & \begin{tabular}[c]{@{}c@{}}\vspace{1.5mm} $\min_{\bbx \in \cc}  \mbE[f(\bbx,\xi)]$\\ $\text{s.t.}$ $\G(\xi)\x \in \cX(\xi)$ \scriptsize{almost surely},
			\end{tabular} 
	\end{tabular} }	
\end{table}

\noindent where $\cc\subseteq\mbR^{m}$ is a convex and compact set, $\xi$ is a random variable with unknown distribution, and $f(\cdot, \xi):\Rn^m \rightarrow \Rn$ is a smooth and convex function. The matrices $\G, \G(\xi) \in \Rn^{n \times m}$ are arbitrary, and the sets $\cX,\cX(\xi) \subseteq \Rn^n$ are convex. We write $\G(\xi)$ and $\cX(\xi)$ as explicit functions of $\xi$ in order to emphasize their stochastic nature and require the constraint $\G(\xi)\x \in \cX(\xi)$ to be satisfied with probability one. While such a requirement is stronger than satisfying the constraints in expectation, it is also impossible to check. Therefore, we will instead seek to find an $\x$ whose expected distance from the feasible set can be made arbitrarily small.

Affine-constrained stochastic optimization problems of the form in $\pd$ arise in a number of areas, such as control theory, multiple kernel learning, blind deconvolution, matrix learning, and communications \cite{shapiro2014lectures,ahmed2013blind}. Efficient algorithms for solving these problems make use of the stochastic gradient $\nabla f(\x,\xi)$ at every iteration, usually obtained from one or a few data points. In many cases, such a stochastic gradient can be calculated efficiently and exactly, giving rise to the first-order stochastic gradient algorithms. On the other hand, in some problems, such as those arising in simulation-based optimization \cite{conn2009introduction}, one must contend with only the function values $f(\x,\xi)$, which are then used to approximate the required stochastic gradients $\nabla f(\x,\xi)$. Algorithms using such gradient estimators are referred to as gradient-free or zeroth-order (ZO) algorithms, and find applications in parameter estimation and classification of black-box systems\cite{chen2017zoo}. Finally, deterministic constraints such as those in $\pd$ can be handled in different ways, such as through the use of projection \cite{nemirovski2009robust}, homotopy \cite{tran2018smooth}, or duality theory \cite{yurtsever2019conditional2a}.  

The problem $\ps$ is less standard, and involves stochastic constraints, which must be satisfied with probability one. Such a formulation is useful in settings where the constraints are not known in advance or are too many to be processed in one batch. For instance, the constraints in some online learning problems, such as online portfolio optimization \cite{fercoq2019almost} and online compressed sensing \cite{garrigues2008homotopy}, are revealed in a sequential fashion. Likewise, the number of constraints in the semidefinite relaxed versions of $k$-means clustering \cite{peng2007approximating}, uniform sparsest cut \cite{arora2009expander}, and maximum a posterior estimation \cite{huang2014scalable} problems, is very large, ranging from $\O(m^2)$ to $\O(m^3)$. In both such cases, efficient algorithms for solving $\ps$ must rely on using only a minibatch of the constraints at every iteration. Again, both first-order and ZO algorithms are well-motivated for solving $\ps$.

In this work, we focus on problems of the form $\pd$ and $\ps$, where the projection onto the set $\cc$ is computationally expensive, but projection onto the set $\cX$ is simpler. Such templates are prevalent
	among semidefinite programming (SDP) problems, where $\cc$ is the positive semidefinite (PSD) cone, while $\cX$ incorporates simple affine or box constraints. Examples include relaxations of combinatorial optimization problems such as quadratic assignment \cite{zhao1998semidefinite}, maximum cut \cite{goemans1995improved}, etc. SDPs are also encountered frequently in machine learning problems such as for certifying the robustness of neural networks \cite{raghunathan2018semidefinite}, unsupervised  clustering and embedding \cite{kulis2007fast}, and Lipschitz constant estimation \cite{latorre2020lipschitz}. However, solving such large-scale SDPs is almost always difficult, even when utilizing first- or zeroth-order stochastic methods due to the presence of the PSD cone constraint. Projection onto the PSD cone requires carrying out a singular value decomposition, which is a prohibitively expensive operation in large-scale settings. 

Let us motivate the issue by illustrating it on the SDP formulation of the \emph{$k$-means clustering} problem \cite{peng2007approximating}:
	\begin{align}
	&\min_{\X \in \cc} \frac{1}{|\Omega|}\sum_{(i,j)\in\Omega}[\M]_{ij}[\X]_{ij}\label{obj_p2}\\
	&\text{s.t.}\; \X\textbf{1}=\textbf{1},\;[\X]_{ij}\geq 0,\;(i,j)\in \Ec\label{cons_p2},
	\end{align}
	where $\M \in \mathbb{R}^{N\times N}$ is the pairwise distance matrix, \textbf{1} denotes the  $N\times 1$ vector of all ones, $\cc := \{\X \succeq 0, \text{tr}(\X)\leq K \}$, $\Ec := \{(i,j)\mid 1\leq i, j \leq N\}$, and $\Omega\subseteq \Ec$. A possible solution to problem \eqref{obj_p2} is to incorporate the additional constraint \eqref{cons_p2} into the objective and then use a projection-based algorithm. However, projection onto the set $\cc$ will require the eigenvalue decomposition of an $N \times N$ matrix, thus incurring an $\O(N^3)$ complexity per iteration.

Stochastic Frank-Wolfe (FW) algorithms, also known as conditional gradient methods (CGM), avoid the complicated projection operation, and have been applied to solve $\pd$ and $\ps$ in \cite{locatello2019stochastic} and \cite{vladarean2020conditional}, respectively. In the FW class of algorithms, the projection step is replaced with a linear minimization step, which in many cases, is significantly cheaper.
For instance, when $\cc$ is a PSD cone, the linear minimization step involves finding only the largest singular vector, which can be efficiently calculated by using the Lanczos method \cite{jaggi2013revisiting}. 

The affine constraints in $\pd$ and $\ps$ present another challenge, as they cannot generally be incorporated within $\cc$. For instance, in the $k$-means clustering problem \eqref{obj_p2}, incorporating the constraints \eqref{cons_p2} into $\cc$ will no longer allow us to use the efficient Lanczos method. The works in \cite{locatello2019stochastic} and \cite{vladarean2020conditional} proposed incorporating the constraints within the objective function using indicator functions, and subsequently applying homotopy and Nesterov smoothing techniques \cite{nesterov2005smooth} to make them differentiable. 

Finally, the stochastic nature of $\pd$ and $\ps$ obviates the use of the classical FW algorithm. The stochastic FW (SFW) algorithms, unlike the stochastic projected gradient counterparts, require the use of increasingly accurate gradient estimates, i.e., gradient estimates whose bias \emph{and} variance decrease with iterations. Typically, such gradient estimates are constructed using mega-batches whose size increases with iterations, but result in impractical algorithms. More recently, one-sample-per-iteration variants of SFW have become popular, where a gradient tracking approach is utilized to maintain good quality gradient estimates; the same is employed in \cite{locatello2019stochastic} and \cite{vladarean2020conditional}. However, a key issue with the gradient tracking approach in these works is that it results in suboptimal convergence rates of the SFW algorithm. 

Building on the techniques introduced in  \cite{locatello2019stochastic}-\cite{vladarean2020conditional}, we put forth improved first-order algorithms for solving $\pd$ and $\ps$. The proposed algorithms utilize a superior momentum-based gradient tracking approach and newer convergence proofs that allow us to obtain better convergence rates. While the momentum-based gradient tracking is well-known in the context of non-convex optimization and has been applied to vanilla SFW, its application to constrained problems such as in $\pd$ and $\ps$ throws up technical challenges that must be addressed. Specifically, the presence of the affine constraints, which are incorporated within the objective via smoothed penalty functions, makes the momentum-based gradient tracking process complicated and requires careful choice of step-size, which now becomes coupled with the smoothing parameter.  The optimality gap of the proposed MOmentum-based STochastic FW (MOST-FW) algorithm for solving $\pd$ decays at the rate of $\O(k^{-1/2})$, at par with the rate for the standard SFW case \cite{zhang2020one}-\cite{xie2020efficient}. The corresponding algorithm for $\ps$, namely \ms\!\!, achieves a rate of $\O(k^{-1/4})$, which is again better than the $\O(k^{-1/6})$ rate achieved in \cite{vladarean2020conditional}. We also propose the corresponding ZO variants for \md and \ms\!\!, useful for scenarios where the stochastic gradient of the objective function is not readily available, and must be estimated using ZO information. 

Existing deterministic and stochastic FW algorithms, including the proposed ones, generally require one call to the linear minimization oracle (LMO) per iteration. Consequently, the wall-clock time of large-scale problems is largely dominated by the time required to solve the linear minimization sub-problems \cite{braun2017lazifying,kerdreux2018frank,mhammedi2021efficient}. In order to alleviate the issue, we put forth trimmed-FW variants, wherein LMO calls are made only when the observed stochastic gradient causes sufficient change in the gradient estimate. Otherwise, the LMO call is skipped, and the previously available LMO output is re-utilized to carry out the update. The trimming technique is novel in the context of FW algorithms, and gives rise to the \textbf{T}rimmed-\md and \textbf{T}rimmed-\ms algorithms. Although intuitive as a technique, we also establish that the trimming process does not hurt the overall rate of convergence of the proposed algorithms. From the experiments, however, we observe that the trimmed-FW variants need to make only a fraction of calls to the LMO as compared to their regular (non-trimmed) counterparts. 
%
%
\subsection{Related Work}
We review some of the related work in the context of FW algorithms, smoothing, affine-constrained optimization, momentum techniques, and ZO methods. 
%
\subsubsection{FW for smooth functions} Though the FW method was first proposed in 1956, it has only been recently used for solving large-scale optimization problems in a projection-free manner \cite{jaggi2013revisiting}. Subsequently, online and stochastic variants of the FW method for minimizing smooth convex functions have been widely studied \cite{hazan2012projection,hazan2016variance}. Many of the early variants of the stochastic FW algorithms utilized a double loop structure, wherein the inner loop utilizes a single mini-batch of stochastic gradients per iteration, while the outer loop is used to update some of the algorithm parameters. Equivalently, these algorithms can be seen as using ``checkpoints'' for updating various algorithm parameters, with the interval between checkpoints increasing polynomially \cite{shen2019complexities} or even exponentially \cite{hazan2016variance}. Such mega-batches result in infrequent update of algorithm parameters, which is not desirable in practice \cite{defazio2018ineffectiveness}. Modern stochastic FW algorithms require only a single stochastic gradient as well as a single linear minimization step per iteration \cite{akhtar2021momentum, mokhtari2020stochastic,zhang2020one,xie2020efficient}. Of these, works in \cite{zhang2020one,xie2020efficient} were the first to utilize a momentum-based gradient tracking routine, achieving the state-of-the-art convergence rate of $\O(k^{-1/2})$. 

\subsubsection{FW for non-smooth functions}
Complications arise when extending these results to constrained problems, such as $\pd$ and $\ps$. Unlike proximal methods, the constraints cannot be handled using an indicator penalty function, since the resulting objective would be non-smooth. On the other hand, the stochastic FW variants discussed so far can handle only smooth objective functions. For instance, \cite{akhtar2021conservative} deals with additional expectation constraints using penalty reformulation but requires the functional constraints to be smooth. For non-smooth but Lipschitz continuous objectives, it may still be possible to apply stochastic FW through the use of Nesterov's smoothing \cite{lan2016conditional}. The indicator function is, however, not Lipschitz continuous, and hence the constrained formulations in $\pd$ and $\ps$ are not amenable to the techniques proposed in \cite{lan2016conditional,lan2017conditional}. Another approach, proposed in \cite{lu2020generalized}, entails incorporating the constraints within the linear minimization step. However, such  inclusion may significantly increase the complexity of carrying out the linear minimization, since the special structure present in $\cc$ may be lost.
\begin{table*}[]
	\centering\scalebox{1.1}{
		\begin{tabular}{|c|c|c|c|c|c|c|c|}
			\hline
			Reference                                                & \begin{tabular}[c]{@{}c@{}}Objective\\ Type\end{tabular} & \begin{tabular}[c]{@{}c@{}}Additional\\ Affine Constraint\end{tabular} & \begin{tabular}[c]{@{}c@{}}Affine\\ Constraint Type\end{tabular} & \begin{tabular}[c]{@{}c@{}}Optimality \\ Gap\end{tabular} & \begin{tabular}[c]{@{}c@{}}Constraint\\ Feasibility\end{tabular} & \begin{tabular}[c]{@{}c@{}}Query\\ Size\end{tabular}           & Oracle \\ \hline
			SFW \cite{mokhtari2020stochastic}                                                    & Stochastic                                               & \xmark                                                                      & --                                                             & $\mathcal{O}(k^{-1/3})$                                   & --                                                             & $\mathcal{O}(1)$ & SFO    \\ \hline
			1-SFW \cite{zhang2020one}                                                    & Stochastic                                               & \xmark                                                                      & --                                                             & $\mathcal{O}(k^{-1/2})$                                   & --                                                             & $\mathcal{O}(1)$ & SFO   
			
			\\ \hline
			ORGFW \cite{xie2020efficient}                                                    & Stochastic                                               & \xmark                                                                      & --                                                             & $\mathcal{O}(k^{-1/2})$                                   & --                                                             & $\mathcal{O}(1)$ & SFO   
			
			\\ \hline
			
			HFW\cite{yurtsever2018conditional}      & Deterministic                                            & \cmark                                                                     & Deterministic                                                    & $\mathcal{O}(k^{-1/2})$                                   & $\mathcal{O}(k^{-1/2})$                                          & --             & FO    \\ \hline
			SHCGM\cite{locatello2019stochastic}      & Stochastic                                               & \cmark                                                                     & Deterministic                                                    & $\mathcal{O}(k^{-1/3})$                                   & $\mathcal{O}(k^{-5/12})$                                         & $\mathcal{O}(1)$ & SFO    \\ \hline
			H1-SFW\cite{vladarean2020conditional} & Stochastic                                               & \cmark                                                                      & Stochastic                                                       & $\mathcal{O}(k^{-1/6})$                                   & $\mathcal{O}(k^{-1/6})$                                          & $\mathcal{O}(1)$ & SFO    \\ \hline
			\rowcolor{MistyRose1}\textbf{\md}                        & Stochastic                                               & \cmark                                                                      & Deterministic                                                    & $\mathcal{O}(k^{-1/2})$                                   & $\mathcal{O}(k^{-1/2})$                                          & $\mathcal{O}(1)$ & SFO    \\ \hline
			\rowcolor{LightCyan1}\textbf{\ms}                    & Stochastic                                               & \cmark                                                                      & Stochastic                                                       & $\mathcal{O}(k^{-1/4})$                                  & $\mathcal{O}(k^{-1/4})$                                         & $\mathcal{O}(1)$ & SFO    
			\\ [0.75ex] \hline \hline
			ZO-FW \cite{sahu2019towards}                                                     & Deterministic                                               & \xmark                                                                      & --                                                             & $\mathcal{O}(k^{-1})$                                   & --                                                             & $\mathcal{O}(m)$ & ZO    \\ \hline 
			ZO-SFW \cite{sahu2019towards}                                                     &  Stochastic                                              & \xmark                                                                      & --                                                             & $\mathcal{O}(k^{-1/3})$                                   & --                                                             & $\mathcal{O}(m)$ & SZO      \\ \hline
			\rowcolor{MistyRose1}\textbf{\md}                        & Stochastic                                               & \cmark                                                                     & Deterministic                                                    & $\mathcal{O}(k^{-1/2})$                                   & $\mathcal{O}(k^{-1/2})$                                          & $\mathcal{O}(m)$ & SZO    \\ \hline
			\rowcolor{LightCyan1}\textbf{\ms}                    & Stochastic                                               & \cmark                                                                     & Stochastic                                                       & $\mathcal{O}(k^{-1/4})$                                  & $\mathcal{O}(k^{-1/4})$
			& $\mathcal{O}(m)$ & SZO    \\ \hline	\end{tabular}}
	\caption{Summary of related works. For zeroth order oracle, query size indicates the number of  function evaluations required in estimating an $m$-dimensional gradient. For first-order oracle, it denotes the number of stochastic gradient samples required at each iteration.}
	\label{table1}
\end{table*}
The deterministic counterpart of $\pd$ was first considered in \cite{yurtsever2018conditional}, and combines the ideas of homotopy and smoothing to handle the non-smooth component of the objective. The idea there was to replace the non-smooth component with its smooth approximation. The error due to this approximation is controlled by decreasing the smoothing parameter (and hence tightening the approximation) at an appropriate rate. The problem in $\pd$ was first considered in \cite{locatello2019stochastic}, which again used the homotopy and smoothing ideas from \cite{yurtsever2018conditional}, but used gradient tracking estimator to handle the stochastic gradient noise. The stochastic homotopy CGM (SHCGM) of \cite{locatello2019stochastic} attains a convergence rate of $\O(k^{-1/3})$ for the optimality gap and $\O(k^{-5/12})$ for the constraint violation. As we shall establish later, the proposed MOST-FW algorithm also uses homotopy and smoothing, but combines it with momentum-based gradient tracking, thus achieving a rate of $\O(k^{-1/2})$ for both, optimality gap and constraint violation. 

Stochastically constrained stochastic optimization problems of the form in $\ps$ have been well-studied in the context of proximal stochastic algorithms \cite{patrascu2017nonasymptotic} and proximal online algorithms \cite{fercoq2019almost}. The works in \cite{wang2015random} and \cite{ patrascu2017nonasymptotic} can also accommodate large number of constraints, but cannot handle the additional difficult-to-project set $\cc$. Motivated by these works, the development of a projection-free algorithm for solving $\ps$ was pursued in \cite{vladarean2020conditional}, building upon the techniques from \cite{locatello2019stochastic}. The H1-SFW algorithm in \cite{vladarean2020conditional} achieves a convergence rate of $\O(k^{-1/6})$ for the optimality gap and constraint violation. In comparison, the \ms algorithm proposed here achieves the rate of $\O(k^{-1/4})$ for the same. 

\subsubsection{Variance-reduced FW} Variance-reduced algorithms have been widely used in the context of FW methods; examples include SCGS \cite{lan2016conditional}, SVRF \cite{hazan2016variance}, STORC \cite{hazan2016variance}, SAGA-FW \cite{ reddi2016stochastic}, and  SPIDER-FW \cite{yurtsever2019conditional}. However, none of these approaches can handle affine constraints as in $\pd$ or $\ps$. A variance-reduced approach called  H-SPIDER-FW is proposed in \cite{vladarean2020conditional}, but requires stochastic gradient batch-sizes that increase exponentially with the iteration index. The proposed \ms algorithm converges at the same rate but works with a standard mini-batch of stochastic gradients, akin to the state-of-the-art FW algorithms in the standard setting.  Recently, a  momentum based \textit{deterministic} accelerated-FW (AFW) algorithm was proposed in \cite{li2021momentum} by mimicking the steps of accelerated gradient
methods. Although it covers many important cases, it is not applicable to general constraint set $\mathcal{C}$. Unlike AFW, propose stochastic algorithms that perform momentum-based tracking over the gradient with the target to reduce the gradient approximation noise.
%
%
\subsubsection{FW under zeroth order oracle} Projection-free ZO algorithms have been studied in \cite{balasubramanian2018zeroth,sahu2019towards,huang2020accelerated}. Of these, the work in  \cite{balasubramanian2018zeroth} proposed a  stochastic zeroth-order (SZO) FW algorithm using the gradient estimation technique of \cite{nesterov2017random}, but the number of samples of directional derivatives required at each epoch scaled linearly with the iteration index as well as the problem dimension. An improved variant was proposed in  \cite{sahu2019towards}, which adopted the gradient averaging technique from \cite{mokhtari2020stochastic}, achieving a convergence rate of $\O(k^{-1/3})$ with a query size of $\O(m)$. Recently, a ZO projection-free algorithm was proposed in \cite{huang2020accelerated}, and resulted in state-of-the-art rates for the non-convex case. To the best of our knowledge, there are no ZO algorithms for solving $\pd$ or $\ps$.
\subsection{FW with reduced LMO calls}
There have been few attempts at reducing the number of LMO calls of conditional gradient methods, all in deterministic settings. Among these, the oldest one is the lazy conditional gradient approach  designed for deterministic FW that replaces the linear optimization oracle by a (weak) separation oracle which approximately solves a separation problem. One recent direction to improve the LMO complexity is to use a randomized linear oracle \cite{braun2017lazifying,kerdreux2018frank,frandi2014complexity}, wherein linear minimization  is performed only over a random sample of the original atomic domain. However, the effectiveness of such a technique depends on whether a specified fraction of the constraint set can be efficiently subsampled. Another relevant recent approach is that of using the set-membership oracle instead of the LMO \cite{mhammedi2021efficient}. The approach does not generally improve upon the FW algorithm, because for many problems, linear minimization may actually be cheaper than determining the set membership, even approximately. In summary, no existing works have attempted to improve the LMO complexity by conditionally skipping the LMO calls. Such a trimming process is novel in the context of FW algorithms, but bears resemblance to the \emph{censoring} idea in wireless sensor networks and more recently in distributed optimization \cite{li2021communication}.
\subsection{Contributions}
In this work, we develop stochastic first- and zeroth-order FW algorithms for solving $\pd$ and $\ps$. The proposed algorithms uses momentum-based gradient tracking and novel proof techniques to provide stronger guarantees than the existing
FW variants \cite{locatello2019stochastic,vladarean2020conditional},  while still using one mini-batch per iteration. Our contributions are summarized as follows:
\begin{itemize}
	\item We propose \md to solve $\pd$ and show that it achieves the state-of-the-art convergence rate of $\mathcal{O}(k^{-1/2})$ for both optimality gap and constraint feasibility. Remarkably, the state-of-the-art convergence rate of SFW for solving the standard version (without additional constraint set $\cX$) of $\pd$ is also $\O(k^{-1/2})$ \cite{zhang2020one,xie2020efficient}. 
	
	\item We further propose \ms to solve $\ps$ that achieves a convergence rate of $\mathcal{O}(k^{-1/4})$ for both optimality gap and constraint feasibility. The proposed \ms algorithm is a fully stochastic version of \md and follows the same idea of momentum-based gradient tracking; however, now it tracks the gradient of both objective and smoothed affine constraints. 
	\item For the first time, we propose ZO methods for solving the above problems. Interestingly, the proposed ZO versions have the same iteration complexity as their first-order counterparts. The obtained rates even \emph{improve} over the state-of-the-art ZO algorithms for solving the set-constrained version of $\pd$ or $\ps$ \cite{sahu2019towards}.
	\item We further propose trimmed variants of both the algorithms called T-\md and T-\ms by employing a novel trimming technique to improve LMO complexity while maintaining identical convergence rates up to small constant factors. The key idea is to skip the LMO calls when the observed stochastic gradients are not sufficiently new. 
	\item We provide numerical evidence of the superiority of the proposed approaches on various relevant applications, namely, sparse matrix estimation, clustering via semidefinite relaxation, and uniform sparsest cut problem. We also demonstrate the computational superiority of the trimmed variants by extensive numerical comparisons  with their non-trimmed version. In all cases, we report a significant reduction in the total number of LMO calls.
\end{itemize}

A complete comparison of the proposed algorithm with  related works is provided in Table \ref{table1}. We only compare with algorithms that make use of a standard mini-batch of stochastic gradients per iteration, and do not include the algorithms relying on mega-batches, as those are generally considered impractical \cite{defazio2018ineffectiveness}.  

The rest of this paper is organized as follows. We begin with reviewing the notions of smoothing, linear minimization oracle, and the coordinate-wise gradient estimation techniques before proceeding to discuss the proposed algorithms of \md and \ms in Sec. \ref{sec:Algo_Dev}. Various theoretical results are presented in Sec. \ref{sec:convergence}, while the numerical validation is provided in Sec. \ref{sec:application}. Finally, Section \ref{sec:conclusion} concludes the paper.

\textbf{Notation:} A matrix (vector) is denoted by uppercase (lowercase) letters in bold font. The $(i,j)$-th element of a matrix $\X$ is denoted by $[\X]_{ij}$. The notation $\norm{\cdot}$ refers to the norm, which when applied to vectors, represents the Euclidean norm and when applied to matrices, the spectral norm. The $\ell_1$ and Frobenius norms are denoted by $\norm{\cdot}_1$ and $\norm{\cdot}_F$, respectively. The inner product is represented by $\ip{\cdot,\cdot}$. The distance between the point $\x$ and the set $\cX$ is denoted by $\cD_{\cX}(\x):=\text{inf}_{\u\in \cX}\norm{\u-\x}$.

\section{Algorithm Development}\label{sec:Algo_Dev}
In this section we develop the \md and \ms algorithms (and their trimmed version T-\md and T-\ms) for solving the problems $\pd$ and $\ps$, respectively. We begin with discussing some preliminaries. The performance of the proposed algorithms will be characterized in terms of their oracle complexities. Depending on the problem at hand, we allow two possible choices of the oracle:
\begin{itemize}
	\item The \emph{Stochastic First-Order (SFO)} oracle, which provides $\nabla f(\x,\xi)$ for a given $\x$; and
	\item The \emph{Stochastic Zeroth-Order (SZO)} oracle provides $f(\x,\xi)$ for a given $\x$. 
\end{itemize}
When using the SZO oracle, the stochastic gradient is estimated using the so-called coordinate-wise gradient estimator (CGE) \cite{liu2018zeroth,ji2019improved}:  
\begin{align}\label{coo:grad_apprx}
\nt f (\x,\xi)=\sum_{i=1}^{m}\frac{f(\x+\rho\u_i,\xi)-f(\x-\rho\u_i,\xi)}{2\rho}\u_i,
\end{align}
where $\rho$ is the element-wise smoothing parameter, and $\u_i\in \mathbb{R}^m$ is a standard basis vector with $[\u_i]_j=1$ if $j=i$, otherwise zero. 

Observe here that the calculation of $\nt f(\x,\xi)$ requires $2m$ calls to the SZO oracle. For the sake of brevity, we will henceforth use $g(\x,\xi)$ to denote the stochastic gradient, with the understanding that $g(\x,\xi)$ may either be $\n f(\x,\xi)$ or $\nt f(\x,\xi)$, depending on the oracle being used. 

As stated earlier, we are interested in settings where the projection over $\cX$ is easy but the projection over $\cc$ is difficult. In particular, we will require access to a linear minimization oracle (LMO), that provides the solution to the optimization problem $\min_{\u \in \cc} \ip{\u,\x}$ for a given $\x$. The number of calls to the LMO will also be a performance metric. 

Towards using the SFW framework to solve these problems, let us define the indicator function
\begin{align}
\one_{\cX}(\x)=  \begin{cases} 0 & \x\in \cX \\
\infty & \x\notin \cX,
\end{cases}
\end{align}
which allows to write the problems $\pd$ and $\ps$ compactly as
\begin{align}\label{prob1}
\min_{\x\in \cc}F(\x) &:= \mbE[f(\x,\xi)] + \one_{\cX}(\G\x)  \\
\text{and} ~~~ 	\min_{\x\in \cc}\hat{F}(\x) &:= \mbE[f(\x,\xi) + \one_{\cX_{\xi}}(\G(\xi)\x)], \label{prob2}
\end{align}
respectively. Observe here that only the affine constraints have been incorporated within the indicator function while the constraint $\x \in \cc$ is retained as is. 

A projection free algorithm, relying instead on the LMO described earlier, along the lines of \cite{hazan2012projection,hazan2016variance,mokhtari2020stochastic,xie2020efficient} may now be applied to \eqref{prob1}-\eqref{prob2}. Note however that the analysis of SFW in these works requires the objective function to be smooth, which is not the case here due to the presence of indicator functions. To this end, we must use a smooth approximation of the objective functions. We adopt the Nesterov's smoothing technique \cite{yurtsever2018conditional}, which entails replacing the indicator function $\one_{\cX}(\y)$ in \eqref{prob1}-\eqref{prob2} with its smooth approximation \cite{beck2017first}:
\begin{align}
h_\mu(\y,\cX) =\frac{1}{2\mu}\cD^2_{\cX}(\y) = \frac{1}{2\mu}\norm{\y - \Pi_{\cX}(\y)}^2_2 \label{smooth}
\end{align}
where $\mu > 0$ is an algorithm parameter, $\Pi_{\cX}(\y):=\argmin_{\u \in \cX}\norm{\u-\y}^2$  is the projection operator, and $\cD_{\cX}(\y):= \norm{\y - \Pi_{\cX}(\y)}_2$. More generally, the scaled \emph{squared set-distance} function $h_{\mu}$ is the Moreau envelope of $\one_{\cX}$, and therefore convex as well as  $\frac{1}{\mu}$-smooth \cite{beck2017first}. It is remarked that the idea of using a smooth approximation of the objective, so as to allow the use of projection free algorithms (such as SFW) is well known \cite{fercoq2019almost,tran2018smooth}. Typically, the parameter $\mu$ must be carefully tuned so as to ensure that the approximation error, arising from the use of $h_{\mu}$ instead of $\one_{\cX}$, remains less than or equal to the optimality gap at a given iteration. 
\subsection{Performance Metrics} 
We will analyze the performance of the proposed algorithms in terms of the following parameters
\begin{itemize}
	\item Average optimality gap
	$\mbE[f(\x,\xi)]-f(\x^\star)$; and
	\item Average constraint violation give by 
	\begin{align}
	D(\x) := \begin{cases}
	\mbE[\cD_{\cX}(\G\x)] & \text{for $(\pd)$} \\
	\mbE[\cD_{\cX_\xi}(\G(\xi)\x)] & \text{for $(\ps)$} 
	\end{cases}	
	\end{align}
\end{itemize}
A point $\tilde{\x}$ is said to be $(\eps,\delta)$-optimal if it satisfies $\mbE[f(\tilde{\x})]-f(\x^{\star})\leq \eps$ and $D(\x) \leq \delta$. The algorithms will be designed such that the constraint $\x \in \cc$ will always be satisfied for all iterates. The goal will be to obtain the SFO/SZO complexities of finding an $(\eps,\delta)$-optimal solution to $(\pd)$ and $(\ps)$. It is remarked that for the proposed algorithms, the LMO complexity, which counts the number of calls to the LMO oracle, is the same as the SFO/SZO complexity as both \md and \ms solves only a single linear minimization problem at each iteration. Having discussed the preliminaries, we are now ready to detail the proposed algorithms. 
\begin{algorithm}[t]
	\caption{\colorbox{MistyRose1}{\md}}
	\label{alg:1}
	\begin{algorithmic}
		\State {\bfseries Initialization:} $\x_0,\x_1\in \cc, \y_0$, parameters $\{\gamma_k,\eta_{k}, \mu_k, \rho_k\}$
		\For{$k=1,2,\ldots$}
		\State Compute: $\y_k=(1-\gamma_k)\y_{k-1} + \gamma_kg(\x_k,\xi_k)$\\
		\hspace{3cm}$+ (1-\gamma_k)(g(\x_{k},\xi_k)-g(\x_{k-1},\xi_k))$\\
		where,
		\begin{align*}
		g(\x,\xi)=\begin{cases*}
		\nt f (\x,\xi) &\textit{for SZO oracle}\\
		\n f(\x,\xi) &\textit{for SFO oracle}
		\end{cases*}
		\end{align*}			
		\State Compute: $\w_k=\y_k+ \mu_k^{-1}\G^{T}\left(\G\x_k-\Pi_{\cX}(\G\x_k)\right)$
		\State Compute: $\z_k=\argmin_{\z\in \cc}\ip{\z,\w_k}$
		\State Update: $\x_{k+1}= \x_k+\eta_k(\z_k-\x_k)$		
		\EndFor 
	\end{algorithmic}
\end{algorithm}
\subsection{\colorbox{MistyRose1}{\md}}
We begin with writing the smoothed approximation of $\pd$, which takes the form:
\begin{align}\label{prob1_smoooth}
\min_{\x\in \cc}F_{\mu}(\x):=\Ex{f(\x,\xi)}+h_{\mu}(\G\x,\cX) \tag{$\pdm$}.
\end{align}
For the sake of brevity, let us denote $F_{\mu}(\x,\xi):= f(\x,\xi) + h_{\mu}(\G\x,\cX)$, so that $F_{\mu}(\x) = \Ex{F_{\mu}(\x,\xi)}$. 

As mentioned earlier, the problem in $\pd$ was first considered in \cite{locatello2019stochastic}, where the first order SFW variant called SHCGM was proposed for the case when $\n f(\x,\xi)$ was available. The SHCGM approach relies on three key steps: (a) tracking $\n \E {f(\x_k,\xi)}$ through the recursive update rule
\begin{align}\label{old_estimator}
\y_k=(1-\gamma_k)\y_{k-1}+\gamma_k\n f(\x_k,\xi_k),
\end{align}
(b) calling the LMO to solve the problem
\begin{align}\label{lmo_alg1}
\z_k=\argmin_{\z\in \cc}\ip{\z,\y_k+\n_{\x} h_{\mu_k}(\G\x_k,\cX)},
\end{align}
and (c) carrying out the update 
\begin{align}\label{update_alg1}
\x_{k+1} &= \x_k+\eta_k(\z_k-\x_k).
\end{align}
It can be seen that the gradient of the smooth approximation $h_{\mu}$ is well-defined and given by
\begin{align}\label{grad_cons_funct}
\n_{\x} h_{\mu}(\G\x,\cX)&=\G^T\n h_{\mu}(\G\x,\cX)=\frac{1}{\mu}\G^{T}(\G\x-\Pi_{\cX}(\G\x)).
\end{align}

We emphasize that one cannot directly use the stochastic gradient $\nabla F_{\mu_k}(\x_k,\xi_k)$ in the place of $\y_k+\nabla_{\x} h_{\mu_k}(\G\x_k,\cX)$ in  \eqref{lmo_alg1}, as the resulting algorithm does not converge. This is because the variance of the stochastic gradient $\nabla f(\x_k,\xi_k)$ does not go to zero with $k$. The situation can be remedied by using the average of several iid samples of $\n f(\x_k,\xi)$, though such an approach yields poor SFO complexity \cite{hazan2016variance}. A more sophisticated variance-reduction approach was proposed in \cite{hazan2016variance}, which achieved a better SFO complexity but still required prohibitively large batches of samples per iteration. Finally, tracking-updates in \eqref{old_estimator} were proposed in \cite{mokhtari2020stochastic}, wherein $\y_k$ serves as a biased estimator for $\Ex{\n f(\x_k,\xi)}$ but has a much lower variance, and hence achieves the same SFO complexity but using a single sample (or a small mini-batch) per iteration. Similar gradient-tracking updates have since been used in other SFW variants \cite{sahu2019towards,locatello2019stochastic,vladarean2020conditional}.

The proposed MOST-FW algorithm introduces two key innovations over these variants. First, to allow for the SZO oracle, we replace $\n f(\x,\xi)$ with $g(\x,\xi)$, which could either be $\n f(\x,\xi)$ or $\nt f(\x,\xi)$. Second, we make use of a momentum-based gradient tracker instead of that in \eqref{old_estimator}, which takes the form:
\begin{align}\label{estimator}
\y_k &=(1-\gamma_k)\y_{k-1} + \gamma_kg(\x_k,\xi_k)+ (1-\gamma_k)(g(\x_{k},\xi_k)-g(\x_{k-1},\xi_k)).
\end{align}
\begin{algorithm}[t]
	\caption{\colorbox{LightCyan1}{\ms}}
	\label{alg:2}
	\begin{algorithmic}
		\State {\bfseries Initialization:} $\x_0,\x_1\in \cc, \y_0$, parameters $\{\gamma_k,\eta_{k}, \mu_k, \rho_k\}$
		\For{$k=1,2,\ldots$}
		\State Compute: $\y_k=(1-\gamma_k)\y_{k-1} + \gamma_k\check{g}_k(\x_k,\xi_k)$\\
		\hspace{2cm}$+ (1-\gamma_k)(\check{g}_k(\x_{k},\xi_k)-\check{g}_{k-1}(\x_{k-1},\xi_k))$\\
		where, $\check{g}_k(\x,\xi) = g(\x,\xi) + \n h_{\mu_k}(\G(\xi)\x,\cX_\xi)$ and
		\begin{align*}
		g(\x,\xi)=\begin{cases*}
		\nt f (\x,\xi) &\textit{for SZO oracle}\\
		\n f(\x,\xi) &\textit{for SFO oracle}
		\end{cases*}
		\end{align*}			
		\State Compute: $\z_k=\argmin_{\z\in \cc}\ip{\z,\y_k}$
		\State Update: $\x_{k+1}= \x_k+\eta_k(\z_k-\x_k)$		
		\EndFor 
	\end{algorithmic}
\end{algorithm}
The momentum term in \eqref{estimator} was first introduced in \cite{cutkosky2019momentum} in the context of classical stochastic gradient descent, and has been shown to improve the gradient tracking performance. In the present case, we will establish that the momentum-based updates in \eqref{estimator} yield improved oracle complexity bounds. The updates in  \eqref{lmo_alg1}-\eqref{update_alg1} remain the same, and the full algorithm is summarized in Algorithm \ref{alg:1}. Recall that $\rho_k$ is the smoothing parameter associated with the gradient estimator. The choice of various algorithm parameters $\{\gamma_k,\eta_k,\mu_k,\rho_k\}$ will be specified later.

\subsection{\colorbox{LightCyan1}{\ms}}
The derivation of the \ms algorithm for solving $\ps$ follows along similar lines. We begin by first writing down the smoothed approximation of \eqref{prob2} as
\begin{align}\label{prob2_smoooth}
\min_{\x\in \cc} \hat{F}_{\mu}(\x):=\mbE[f(\x,\xi)+h_{\mu}(\G(\xi)\x,\cX_\xi)] \tag{$\mathcal{P}_{\mu}^{+}$}
\end{align}
where $h_{\mu}$ is the smooth approximation of $\one_{\cX_\xi}$ defined in \eqref{smooth}. For the sake of brevity, we denote $\n \hat{F}_{\mu}(\x,\xi):= \n f(\x,\xi) + \n_{\x} h_{\mu}(\G(\xi)\x,\cX_{\xi})$.
Observe here that the gradient of $h_{\mu}$, given by 
\begin{align}\label{grad_Hh}
\n_{\x} h_{\mu}(\G(\xi)\x,\cX_\xi)&=\G(\xi)^T\n h_{\mu}(\G(\xi)\x,\cX_\xi)\\&=\frac{1}{\mu}\G^\T(\xi)(\G(\xi)\x-\Pi_{\cX_\xi}(\G(\xi)\x))\nonumber
\end{align} 
is also random and must also be tracked. 

To this end, let us define
\begin{align}
\check{g}_k(\x,\xi) = g(\x,\xi) + \n_{\x} h_{\mu_k}(\G(\xi)\x,\cX_\xi)
\end{align}
where, as earlier, $g(\x,\xi)$ may either be $\n f(\x,\xi)$ or $\nt f(\x,\xi)$, depending on the oracle available. Then the gradient tracking update takes the form:
\begin{align}\label{estimator2}
\y_k &=(1-\gamma_k)\y_{k-1} + \gamma_k\check{g}_k(\x_k,\xi_k)+ (1-\gamma_k)(\check{g}_k(\x_{k},\xi_k)-\check{g}_{k-1}(\x_{k-1},\xi_k)).
\end{align}
Observe that different from \eqref{prob1_smoooth}, the smoothed component of the objective  in \eqref{prob2_smoooth} is also random and hence its gradient must also be tracked via \eqref{estimator2}. Next, the LMO is called \eqref{lmo_alg2} and finally we carry out the updates \eqref{update_alg2}
\begin{align}\label{lmo_alg2}
\z_k&=\argmin_{\z\in \cc}\ip{\z,\y_k},\\
\x_{k+1} &= \x_k+\eta_k(\z_k-\x_k).\label{update_alg2}
\end{align}
The full algorithm is summarized in Algorithm \ref{alg:2}. Before proceeding to introduce a novel trimmed variant of our proposed algorithms, we provide a remark below to better understand the use of momentum technique in our work.

\noindent \textit{Remark 1:}  There have been recent attempts to improve the performance  of deterministic Frank-Wolfe algorithms using Nesterov momentum \cite{li2021momentum}. However, faster rates for deterministic FW has been achieved under very specific circumstances when the constraint set $\mathcal{C}$ is either a polytope \cite{kerdreux2021projection}, active strongly convex set \cite{garber2015faster}, or active $\ell_p$-norm ball \cite{li2021momentum}. Different from Nesterov momentum, the goal of momentum in stochastic optimization \eqref{estimator} is to reduce the gradient approximation noise. Our use of this particular momentum technique in the FW context results in key Lemmas \ref{track1} and \ref{track2}, which yield improved rates. The analysis in these lemmas is new and instrumental as they provide us the flexibility to set both the $\gamma_k$ and $\eta_k$ decaying at the same rate, unlike the gradient tracking strategy \eqref{old_estimator}. Of particular note is the usage of momentum to track not only the objective gradient but also the constraint in the MOST-FW$^+$ algorithm, thanks to its versatility.   
\begin{algorithm}[t]
	\caption{T-\md / T-\ms}
	\label{alg:3}
	\begin{algorithmic}
		\State {\bfseries Initialization:}  parameters $\{\gamma_k,\eta_{k},\tau_k, \mu_k, \rho_k\}$ and $\x_0,\;\x_1\;\in\;\cc$, $ \y_0,\v_0$ 
		\For{$k=1,2,\ldots$}
		\State Obtain: 
		\begin{align*}
		\s_k=\begin{cases*}
		\w_k &\textit{for \colorbox{MistyRose1}{\md}}\\
		\y_k &\textit{for \colorbox{LightCyan1}{\ms}}
		\end{cases*}
		\end{align*}			
		\State  Trimming:
		\If{$\norm{\s_k-\v_{k-1}}\geq \tau_k$ \textbf{or} $k=1$}
		\State Set $\v_k=\s_k$
		\State Compute: $\z_k=\argmin_{\z\in \cc}\ip{\z,\v_k}$
		\Else
		\State $\z_k=\z_{k-1}$ and set $\v_k=\v_{k-1}$
		\EndIf
		\State Update: $\x_{k+1}= \x_k+\eta_k(\z_k-\x_k)$		
		\EndFor 
	\end{algorithmic}
\end{algorithm}
\subsection{Trimmed \md and \ms}
This section presents a novel trimmed version of both the proposed algorithm.  The initial steps of the algorithm remains the same, i.e., start with obtaining the tracked gradient $\s_k = \w_k$ for \md and $\s_k = \y_k$ for \ms. However, before proceeding to solve the linear minimization problem (or call to the LMO), we first compare $\s_k$ with its old copy $\v_{k-1}$  and solve  \eqref{lmo_alg1} by calling LMO only when $\s_k$ provides sufficient new information. The sufficiency is determined by comparing  the norm difference $\norm{\s_k-\v_{k-1}}$ with a threshold  $\tau_k$. In nutshell, a new call to LMO is only made if  $\norm{\s_k-\v_{k-1}}\geq \tau_k$, where $\norm{\cdot}$ is the norm induced by the inner product used in the linear minimization step. Finally, we carry out the updates as \eqref{update_alg1}. We will show that for  carefully designed  $\tau_k$, this simple thresholding strategy  can scale down	the number of LMO calls significantly while still converging at almost the same rate. The full algorithm is summarized in Algorithm \ref{alg:3}.

\vspace*{0.5cm}

\section{Convergence Analysis}\label{sec:convergence}
In this section, we study the convergence rate of the proposed \md and \ms algorithms. We begin with stating the assumptions required for the analysis: 
\begin{assumption}\label{smoothness}
	\normalfont
	(\emph{Smoothness}) The objective function $f(\cdot,\xi)$ is $L$-smooth on $\cc$, i.e., $\|\n f(\x,\xi)-\n f(\y, \xi)\|\leq L\|\x-\y\|$ for all $\x$, $\y \in \cc$.
\end{assumption}
\begin{assumption}\label{compact}
	\normalfont (\emph{Compact domain})
	The convex set $\cc$ is compact, so that  $\norm{\x-\y}\leq D$ for all  $\x$, $\y \in \cc$.
\end{assumption}
\begin{assumption}\label{boundedf}
	\normalfont
	(\emph{Bounded variance}): The variance of the stochastic gradients $\n f(\x,\xi)$ is bounded by $\sigma^2$, i.e., 
	\begin{align}
	\E{\norm{\n f(\x,\xi) - \n f(\x)}^2}\leq \sigma^2.
	\end{align}
	where $f(\x):=\Ex{f(\x,\xi)}$. 
\end{assumption}
\begin{assumption}\label{spectral norm}
	\normalfont (\emph{Bounded spectral norm})
	The spectral norm of the linear operators $\G$ and $\G(\xi)$ are bounded, i.e.,
	\begin{align}
	\norm{\G(\xi)}^2 &\leq L_G < \infty\\
	\norm{\G}^2 &\leq L_G < \infty
	\end{align}
\end{assumption}
\begin{assumption}\label{slater}
	\normalfont (\emph{Slater's Condition})
	Slater's condition holds for $\pd$ and $\ps$. In other words, for $\pd$, 
	\begin{align}
	\text{relint}(\cc \times \cX)\cap\{(\x,\b)\in \mathbb{R}^m\times \mathbb{R}^n:\G\x=\b\}\neq \emptyset.
	\end{align}
	and for $\ps$, 
	Let $C:\mathcal{H}\rightarrow \mathbb{R}\cup \{\infty\},\; C(\G\x):=\mbE[\delta_{\cX_{\xi}}(\G(\xi)\x)],$ with the linear operator $\G:\mathbb{R}^m\rightarrow \mathcal{H}$ defined as $\G(\x)(\xi):=\G(\xi)\x, \forall\; \x$, we require that
	\begin{align}\label{sri1}
	0\in sri(dom(C)-\G(dom(f)))
	\end{align} 
	here, \textit{sri} stands for strong relative interior \cite{bauschke2011convex} and defined as
	\begin{align}\label{sri2}
	sriC=\{\x \in C|cone(C-\x)=span(C-\x)\}
	\end{align}
\end{assumption}
Assumptions \ref{smoothness}, \ref{compact}, and \ref{boundedf} are standard in the context of SFW algorithms. Observe that Assumption \ref{smoothness} also implies that $f(\x)$ is $L$-smooth. Further we have from \cite[Lemma 3]{ji2019improved} that 
\begin{align}\label{coro_ZOGrad}
\norm{\nt f(\x,\xi) - \n f(\x,\xi)} \leq \sqrt{m}L\rho,
\end{align}
which in turn, implies that $\norm{\nt f(\x) - \n f(\x)} \leq \sqrt{m}L\rho$.

Assumptions \ref{smoothness} and \ref{boundedf} imply that the variance of $\nt f(\x,\xi)$ is bounded since 
\begin{align}\label{up2}
&\Ex{\norm{\nt f(\x) - \nt f(\x,\xi)}^2} = \mbE\left[\|\nt f(\x)+\n f(\x)\right.\\
&\hspace{1cm}\left.+\n f(\x,\xi)  -\n f(\x) + \nt f(\x,\xi) - \n f(\x,\xi)\|^2\right]\nonumber\\
& \leq 3\Ex{\norm{\nt f(\x)-\n f(\x)}^2} + 3\Ex{\norm{\n f(\x)-\n f(\x,\xi)}^2} \nonumber\\
&\hspace{1cm}+3\Ex{\norm{\n f(\x,\xi)-\nt f(\x,\xi)}^2} \leq 3\sigma^2+6mL^2\rho^2\nonumber
\end{align}
Assumption \ref{spectral norm} is intuitive and required to ensure that the composite objectives $F_\mu$ (cf. \eqref{prob1_smoooth})  and $\hat{F}_{\mu}$ (cf. \eqref{prob2_smoooth}) are smooth; see also \cite{fercoq2019almost,vladarean2020conditional}. Finally, Slater's condition in Assumption \ref{slater} implies that for the convex problems $\pd$ and $\ps$, strong duality holds and that the dual optimum variable, denoted by $\lam^\star$, is bounded and ensures  the constraint set to have a non-empty strong relative interior.

We begin with a brief outline of the proof. We will decrease $\mu$ with each iteration so that the decision variable converges to the original solution as the algorithm proceeds. However, we will observe that the variance of the stochastic gradient $\n h_{\mu}(\G(\xi)\x,\cX_\xi)$ is inversely proportional to the smoothing parameter $\mu$ (c.f. \eqref{var_g}) and hence reducing $\mu$ results in a proportional increase in the gradient approximation noise. The same issue was also encountered in  \cite{vladarean2020conditional} because of which their algorithm (H1-SFW) settled at the convergence rate of $\mathcal{O}(k^{-1/6})$. In this work, we will show that by careful selection of smoothness parameter $\mu$ and using the momentum-based gradient tracking given in \eqref{estimator2}, our algorithm can deal with this issue of gradient approximation noise in a much better way than H1-SFW and  helps in achieving better results. 
\subsection{Convergence Analysis of \colorbox{MistyRose1}{\md}}
We begin with establishing a key lemma characterizing the evolution of the tracking error $\Ex{\y_k - g(\x_k)}$, where $g(\x_k):= \Ex{g(\x_k,\xi)}$. To allow us to present the results in a unified manner, we adopt the convention that $\rho_k = 0$ when using the SFO oracle, and $\rho_k > 0$ when using the SZO oracle. Since a zero value of $\rho_k$ does not make sense in \eqref{coo:grad_apprx}, the usage should be clear from the context. 
\begin{lemma}\label{track1}
	(a) Under Assumptions \ref{smoothness}-\ref{boundedf}, the iterates generated by \md satisfy:
	\begin{align}
	&\Ex{\norm{\y_k - g(\x_k)}^2} \leq (1-\gamma_k)^2\Ex{\norm{\y_{k-1}-g(\x_{k-1})}^2}+ 6\gamma_k^2\sigma^2 + 24m L^2\rho^2_{k-1}+ 6\eta_{k-1}^2L^2D^2 \label{mainlem}
	\end{align}
	(b) For the choice $\gamma_k = \tfrac{1}{k}$, $\eta_k = \tfrac{2}{k+1}$, and $\rho_k \leq \frac{D}{\sqrt{m}(k+1)}$, it holds that 
	\begin{align}
	\Ex{\norm{\y_k - \n f(\x_k)}^2} \leq \frac{16(3\sigma^2 + 25L^2D^2)}{k}.
	\end{align}
\end{lemma}
\begin{IEEEproof}
	We start by subtracting $g(\x_k)$ from both sides of \eqref{estimator} and introducing $g(\x_{k-1})$ on the right to obtain
	\begin{align}
	\y_k - g(\x_k) &=(1-\gamma_k)(\y_{k-1} - g(\x_{k-1})) \nonumber\\
	&\quad + (g(\x_k,\xi_k) - g(\x_k))  + (1-\gamma_k)\big(g(\x_{k-1}) - g(\x_{k-1},\xi_k)\big).\label{ttem1z}
	\end{align}
	From the definition of $g$, it holds that, $\Ex{g(\x) - g(\x,\xi)} = 0$ for all $\x$, implying that the last two terms in \eqref{ttem1z} are zero mean. Hence, taking squared-norm on both sides and taking expectation with respect to $\xi_k$, we obtain
	\begin{align}
	\Ek{\norm{\y_k - g(\x_k)}^2} &= (1-\gamma_k)^2 \norm{\y_{k-1} - g(\x_{k-1})}^2 \nonumber\\
	& \quad+ \mbE_k\Big[\|(g(\x_k,\xi_k) - g(\x_k)) \left. + (1-\gamma_k)\big(g(\x_{k-1}) - g(\x_{k-1},\xi_k)\|^2 \right]\label{yk_1stz}
	\end{align}
	where $\mbE_k$ denotes the expectation with respect to the random variable $\xi_k$, while keeping everything else fixed.  Next, let us consider the second term and bound it separately. Defining 
	\begin{align}
	\sX_k &:= \gamma_k(g(\x_k,\xi_k) - g(\x_k))\\
	\sY_k &:= (1-\gamma_k)(g(\x_k,\xi_k) - g(\x_{k-1},\xi_k)),
	\end{align}
	we  see that the second term in \eqref{yk_1stz} is given by $\En{\sX_k + \sY_k - \mbE_k\sY_k}$ and can be bounded by using the inequality $\En{\sX_k + \sY_k - \mbE_k\sY_k} \leq 2\En{\sX_k} + 2\En{\sY_k}$. 
	
	To obtain bounds on $\En{\sX_k}$ and $\En{\sY_k}$, we consider the SFO and SZO cases separately. 
	
	\noindent \textbf{SFO oracle:} In this case, $g(\x_k,\xi_k) = \n f(\x_k,\xi_k)$, so it follows that
	\begin{align}
	\En{\sX_k}\!=\!\gamma_k^2\Ex{\norm{\n f(\x_k)\! -\! \n f(\x_k,\xi_k)}^2}\! \leq \gamma_k^2\sigma^2 \label{xk1}
	\end{align}
	from Assumption \ref{boundedf}. Likewise, from Assumption \ref{smoothness} we have
	\begin{align}
	\En{\sY_k} &= (1-\gamma_k)^2\Ek{\norm{\n f(\x_k,\xi_k) - \n f(\x_{k-1},\xi_k)}^2} \nonumber\\
	&\leq L^2\norm{\x_k - \x_{k-1}}^2 \nonumber\\
	&\leq L^2\eta_{k-1}^2D^2\label{yk1}
	\end{align}
	where we have also used the update \eqref{update_alg1}, dropped the factor $(1-\gamma_k)^2$, and used Assumption \ref{compact}. 
	
	\noindent \textbf{SZO oracle:} In this case, $g(\x_k,\xi_k) = \nt f(\x_k,\xi_k)$, so that
	\begin{align}
	\En{\sX_k} &= \gamma_k^2\Ex{\norm{\nt f(\x_k) - \nt f(\x_k,\xi_k)}^2} \nonumber\\
	&\leq 3\gamma_k^2(\sigma^2 + 2mL^2\rho_k^2)\nonumber\\
	&\leq 3\gamma_k^2\sigma^2 + 6mL^2\rho_k^2\label{xk2}
	\end{align}
	from \eqref{up2} and since $\gamma_k \leq 1$. Next, expanding $\sY_k$ and using \eqref{coro_ZOGrad}-\eqref{up2}, we obtain
	\begin{align}\label{yk2}
	\En{\sY_k} &= (1-\gamma_k)^2\Ek{\norm{\nt f(\x_{k-1},\xi_k)-\nt f(\x_{k},\xi_k)}^2}\nonumber\\
	&=(1-\gamma_k)^2\mbE_k\|\nt f(\x_{k-1},\xi_k)-\n f(\x_{k-1},\xi_k) +\n f(\x_{k-1},\xi_k) -\n f(\x_{k},\xi_k)+\n f(\x_{k},\xi_k)-\nt f(\x_{k},\xi_k)\|^2\nonumber\\&\leq 3\mbE_k\|\nt f(\x_{k-1},\xi_k)-\n f(\x_{k-1},\xi_k)\|^2+ 3\mbE_k\|\n f(\x_{k},\xi_k)-\nt f(\x_{k},\xi_k)\|^2+ 3\mbE_k\|\n f(\x_{k-1},\xi_k)-\n f(\x_{k},\xi_k)\|^2\nonumber\\
	&\leq 3m L^2\rho^2_k + 3mL^2\rho_{k-1}^2 + 3L^2\|\x_k-\x_{k-1}\|^2\nonumber\\
	&\leq 3m L^2\rho^2_k + 3mL^2\rho_{k-1}^2	+ 3\eta_{k-1}^2L^2D^2.
	\end{align}
	where we again used \eqref{update_alg1}, Assumption \ref{compact}, and dropped the factor $(1-\gamma_k)^2$. Finally, if we ensure that $\rho_k \leq \rho_{k-1}$, then the right-hand side of \eqref{yk2} can be written as $6mL^2\rho_{k-1}^2	+ 3\eta_{k-1}^2L^2D^2$.
	
	Since the bounds for the SFO and SZO cases differ only in constant factors and in terms depending on $\rho_k$, they can be unified as
	\begin{align}
	\Ek{\norm{\sX_k + \sY_k - \Ek{\sY_k}}^2} \leq 6\gamma_k^2\sigma^2 + 24m L^2\rho^2_{k-1}+ 6\eta_{k-1}^2L^2D^2 \label{tempol2z}
	\end{align}
	where recall our convention that $\rho_k = 0$ for SFO case. Substituting \eqref{tempol2z} into \eqref{yk_1stz}, we obtain the required result. 
	
	\noindent (b) Substituting  in \eqref{mainlem}, we obtain
	\begin{align}
	&\Ex{\norm{\y_k - g(\x_k)}^2} \\
	& \leq \left(1-\tfrac{1}{k}\right)^2\Ex{\norm{\y_{k-1}-g(\x_{k-1})}^2}  + 6(\sigma^2 + 8L^2D^2)\frac{1}{k^2}\nonumber
	\end{align}
	Subsequently, using Lemma \ref{simplelemma} with $a = 2$ and $A = 6(\sigma^2 + 8L^2D^2)$, we obtain
	\begin{align}\label{sm1}
	\Ex{\norm{\y_k - g(\x_k)}^2} \leq \frac{24(\sigma^2 + 8L^2D^2)}{k}
	\end{align} 
	For SFO case, we get the desired bound simply by setting $g(\x_k)=\n f(\x_k)$ in \eqref{sm1}. For SZO case, we start with $\Ex{\norm{\y_k - \n f(\x_k)}^2}$, introduce $\nt f(\x_k)$, set $g(\x_k)=\nt f(\x_k)$ in \eqref{sm1}, and recall that $\norm{\n f(\x_k) - \nt f(\x_k)}^2 \leq mL^2\rho_k^2$, so as to yield:
	\begin{align}
	\Ex{\norm{\y_k - \n f(\x_k)}^2} &\leq 2\Ex{\norm{\y_k - g(\x_k)}^2 + \norm{\n f(\x_k) - \nt f(\x_k)}^2} \\
	&\leq  \frac{16(3\sigma^2 + 25L^2D^2)}{k}.\label{sm2}
	\end{align} 
	Observe that the bound in \eqref{sm1} is less than that in \eqref{sm2}, and differs only in constant factors. Therefore, for the sake of simplicity, we use \eqref{sm2} as the unified bound. 
\end{IEEEproof}

Bound for Lemma \ref{track1} establishes that the auxiliary variable $\y_k$ tracks the actual gradient $\n f(\x_k)$ with average error bounded as $\O\left(\frac{\sigma^2+L^2D^2}{k}\right)$.

\noindent\textit{Remark 2}: The gradient tracking approach defined in \eqref{estimator} provides a modification of the approach in \eqref{old_estimator} by adding a correction term  $(1-\gamma_k)(g(\x_k,\xi_k)-g(\x_{k-1},\xi_k))$. This correction term plays a crucial role in variance reduction by exploiting the smoothness of $f(\cdot,\xi)$. To see this, let us inspect the tracking error $E_k=\mbE\norm{\y_k-g(\x_k)}^2$ for both the cases, which yield bounds:
\begin{align}
&E_k\leq (1-\gamma_k)E_{k-1}+\mathcal{O}(\gamma_k^2)+\mathcal{O}(\eta_k^2/\gamma_k^2),\label{mom1}\\&\label{mom2}
E_k\leq (1-\gamma_k)E_{k-1}+\mathcal{O}(\gamma_k^2)+\mathcal{O}(\eta_k^2),
\end{align}
respectively. Observe here that the error dynamics differs only in the last term, where for \eqref{old_estimator}, we have a $\mathcal{O}(\eta_k^2/\gamma_k^2)$ term while for \eqref{estimator}, we only have an $\mathcal{O}(\eta_k^2)$ term. This means that if we want $E_k$ to decrease with $k$, we would have to ensure for \eqref{old_estimator} that the term $\frac{\eta_k^2}{\gamma_k^2} = \mathcal{O}(\gamma_k) \rightarrow 0$, while no such restriction applies to \eqref{estimator}. It turns out that this flexibility is the key to achieving a better rate for the proposed algorithms; in fact we set $\eta_k = \gamma_k = \mathcal{O}(k^{-1})$ for MOST-FW and $\eta_k = \gamma_k = \mathcal{O}(k^{-1/2})$ for MOST-FW$+$ to obtain the required rates. This particular choice of step-size parameters would not have been possible had we used the classical gradient tracking estimator \eqref{old_estimator}. 
Next, we have the following lemma regarding the evolution of the smoothed objective function. 
\begin{lemma}\label{lemma:basic_F(x)_1st}
	(a) Under Assumptions \ref{smoothness}, \ref{compact}, \ref{spectral norm}, and \ref{slater}, and for $\mu_k \geq \mu_{k-1}(1-\eta_k)$, it holds that
	\begin{align}\label{smgaprec}
	&F_{\mu_k}(\x_{k+1})-f(\x^{\star})\leq (1-\eta_k)\left(F_{\mu_{k-1}}(\x_{k})-f(\x^{\star})\right)+\frac{\eta_k^2}{2}\left(L+\frac{L_G}{\mu_k}\right)D^2+\eta_kD\norm{\n f(\x_k)-\y_k}
	\end{align}
	(b) Under Assumptions \ref{smoothness}-\ref{slater}, and for $\eta_k = \frac{2}{k+1}$, $\gamma_k = \frac{1}{k}$, $\rho_k \leq \frac{D}{\sqrt{m}(k+1)}$, and $\mu_k = \frac{\mu_c }{\sqrt{k}}$, the smoothed function gap is bounded by 
	\begin{align}
	\mbE[F_{\mu_k}(\x_{k+1})]\!-\!f(\x^{\star}) \leq \frac{8\sqrt{3}\sigma D\!+\!(41L\!+\!L_G\mu_c^{-1})D^2}{\sqrt{k}}
	\end{align}
\end{lemma}
\textbf{Proof}: The proof of Lemma \ref{lemma:basic_F(x)_1st} is provided in   Sec.\ref{proof-smoothgap} of the supplementary material
and is largely similar to proof of Theorem 9 in \cite{locatello2019stochastic}. For a constant $\mu_c$, Lemma \ref{lemma:basic_F(x)_1st} establishes that the smoothed gap evolves as $\O\left(\frac{\sigma D + LD^2}{\sqrt{k}}\right)$, and will allow us to infer the optimality gap and the constraint violation in the subsequent theorem.  It is remarked here that $\rho_k$ is not a tuning parameter. Ideally, we would like $\rho_k$ to be as small as possible, so as to minimize the approximation error. In practice though, a very small value of $\rho_k$ may not be viable due to computational issues. The results developed here make sense provided that $\rho_k$ is not too large so as to impact the rate of convergence. Therefore, the statement of the results require an explicit upper bound on $\rho_k$ as a function of $k$. 

\begin{thm}\label{theorem1}
	Under Assumptions \ref{smoothness}-\ref{slater} and for $\eta_k = \frac{2}{k+1}$, $\gamma_k = \frac{1}{k}$, $\rho_k \leq \frac{D}{\sqrt{m}(k+1)}$, and $\mu_k = \frac{\mu_c}{\sqrt{k}}$, we have:
	\begin{enumerate}[(a)]
		\item Optimality gap:
		\begin{align}
		\mbE[f(\x_{k})]-f(\x^{\star})\!\leq\! \frac{8\sqrt{3}\sigma D\!+\!(41L\!+\!L_G\mu_c^{-1})D^2}{\sqrt{k}}
		\end{align}
		\item Constraint violation:
		\begin{align}
		\mbE[\cD_{\cX}(\G\x_{k})]\leq  \frac{Q}{\sqrt{k}}
		\end{align}
		where $Q:=2\mu_c\norm{\lam^{\star}}+ 6\sqrt{\sigma D \mu_c} + 10D\sqrt{L\mu_c}+2D\sqrt{L_G}$. 	
	\end{enumerate}	
\end{thm}
\textbf{Proof}: The proof is provided in Sec. \ref{Proof:th1} of the supplementary material. Theorem \ref{theorem1} establishes an upper bound
on the expected suboptimality and constraint feasibility for the iterates generated by \md and shows that they converge to zero at least at the rate of $\mathcal{O}(k^{-1/2})$. These results can be directly translated to oracle complexities of finding an $(\epsilon,\delta)$-optimal solution to $\pd$. Defining $\Upsilon := \sigma D + (L+L_G)D^2$, it can be seen that the SFO oracle complexity is given by $\O\left(\max\{\frac{\Upsilon^2}{\epsilon^2},\frac{\norm{\lam^\star}^2+\Upsilon}{\delta^2}\}\right)$ while the SZO complexity is $\O\left(m\max\{\frac{\Upsilon^2}{\epsilon^2},\frac{\norm{\lam^\star}^2+\Upsilon}{\delta^2}\}\right)$, since gradient calculation requires $2m$ calls to the SZO oracle. This concludes our analysis of the \md algorithm.
\subsection{Convergence Analysis of \colorbox{LightCyan1}{\ms}}
We now consider the \ms algorithm for solving $\ps$. As in \md, we begin with bounding the mean square tracking error $\mbE\|\ffb  -\y_k\|^2$.
\begin{lemma}\label{track2}
	(a) Under Assumptions \ref{smoothness}-\ref{spectral norm}, the iterates generated by \ms satisfy:
	\begin{align}
	\Ek{\norm{\y_k - \check{g}_k(\x_k)}^2} &\leq (1-\gamma_k)^2 \norm{\y_{k-1} - \check{g}_{k-1}(\x_{k-1})}^2 \nonumber\\
	&+12\gamma_k^2\sigma^2 + 60m L^2\rho^2_{k-1}+ 18\eta_{k-1}^2L^2D^2 
	\label{mainlem2}\\& +\frac{4\gamma_k^2L_G^2D^2}{\mu^2_k}+6L_G^2D^2\left(\frac{1}{\mu_k}-\frac{1}{\mu_{k-1}}\right)^2
	+\frac{24\eta_{k-1}^2L_G^2D^2}{\mu_{k-1}^2}\nonumber
	\end{align}
	(b) For the choice $\gamma_k=\frac{1}{k},\;$ $\rho_k\leq\frac{D}{\sqrt{m}(k+1)},$ $ \eta_k=\frac{2}{k+1}$, and  $\mu_k =\frac{\mu_c}{(k+1)^{1/4}}$, it holds that 
	\begin{align}
	\mbE\|\ffb  -\y_k\|^2\leq \frac{96(\sigma^2+12L^2D^2+9L_G^2\mu_c^{-2}D^2)}{\sqrt{k}}.
	\end{align}
\end{lemma}
\textbf{Proof:}  The proof of Lemma \ref{track2} is provided in Appendix \ref{Proof:track2}. Bound for Lemma \ref{track2} establishes that for a constant $\mu_c$, the auxiliary variable $\y_k$ tracks the actual gradient $\ffb$ with average error bounded as $\O\left(\frac{\sigma^2+(L^2+L_G^2)D^2}{\sqrt{k}}\right)$. 
Next,  we present the following lemma on the evolution of the smoothed objective function. 
\begin{lemma}\label{lemma:F-smooth}
	Under Assumptions \ref{smoothness}-\ref{slater},  $\rho_k\leq\frac{D}{\sqrt{m}(k+1)},$  $\gamma_k=\frac{1}{k}$, and $\mu_k =\frac{\mu_c}{(k+1)^{1/4}}$, the smoothed function gap is bounded by	
	\begin{align}\label{F-Smooth}
	\mbE[\hat{F}_{\mu_k}(\x_{k+1})]\!-\!f(\x^\star)
	\!\leq\! \frac{8(4\sqrt{6}\sigma D+(35L+31L_G\mu_c^{-1})D^2)}{k^\frac{1}{4}}
	\end{align}
\end{lemma}
\textbf{Proof}: The proof of Lemma \eqref{lemma:F-smooth} is provided in   Sec.\ref{proof:lemma:F-smooth} of the supplementary material. For a constant $\mu_c$, Lemma \ref{lemma:F-smooth} establishes that the smoothed gap evolves as $\O\left(\frac{\sigma D + (L+L_G)D^2}{k^{\frac{1}{4}}}\right)$, and will allow us to infer the optimality gap and the constraint violation in the subsequent theorem. 

\begin{thm}\label{theorem2}
	Under Assumptions \ref{smoothness}-\ref{slater} and for $\eta_k = \frac{2}{k+1}$, $\gamma_k = \frac{1}{k}$, $\rho_k \leq \frac{D}{\sqrt{m}(k+1)}$, and $\mu_k = \frac{\mu_c}{(k+1)^{1/4}}$, we have:
	\begin{enumerate}[(a)]
		\item Optimality gap:
		\begin{align}
		\mbE[f(\x_{k})]-f(\x^{\star})\leq \frac{8(4\sqrt{6}\sigma D+(35L+31L_G\mu_c^{-1})D^2)}{k^{\frac{1}{4}}}
		\end{align}
		\item Constraint violation:
		\begin{align}\label{tem22}
		\mbE[\cD_{\cX(\xi)}(\G(\xi)\x_{k+1})]\leq\frac{Q^{+}}{k^{\frac{1}{4}}},
		\end{align}	 
		where \small{$$Q^{+}:=2\left(\mu_c \norm{\lam^{\star}}+9\sqrt{\sigma  D \mu_c}+17D\sqrt{L\mu_c}+16D\sqrt{L_G}\right).$$}
	\end{enumerate}	
\end{thm}
\textbf{Proof}: The proof  is provided in Sec. \ref{proof:theorem2} of the supplementary material. Theorem \ref{theorem2} establishes an upper bound
on the expected suboptimality and constraint feasibility for the iterates generated by \ms and shows that they converge to zero at least at the rate of $\mathcal{O}(k^{1/4})$. These results can be directly translated to oracle complexities of finding an $(\epsilon,\delta)$-optimal solution to $\ps$. Defining $\Upsilon := \sigma D + (L+L_G)D^2$,  it can be seen that the SFO oracle complexity is given by $\O\left(\max\{\frac{\Upsilon^4}{\epsilon^4},\frac{\norm{\lam^\star}^4+\Upsilon^2}{\delta^4}\}\right)$ while the SZO complexity is $\O\left(m\max\{\frac{\Upsilon^4}{\epsilon^4},\frac{\norm{\lam^\star}^4+\Upsilon^2}{\delta^4}\}\right)$, since gradient calculation requires $2m$ calls to the SZO oracle. Note that although the variance reduced algorithm H-SPIDER-FW \cite{vladarean2020conditional} achieved the \textit{sfo} complexity of $\mathcal{O}(\eps^{-4},\delta^{-4})$, it requires evaluating a batch of gradients at each iteration with batch size changing with iteration as $\mathcal{O}(2^k)$. This restricts the use of H-SPIDER-FW in online setting where at each iteration we have access only to a single stochastic gradient.

\begin{figure}[t!]
	\centering
	\begin{subfigure}{0.5\columnwidth}		\includegraphics[width=\linewidth, height = 0.7\linewidth]{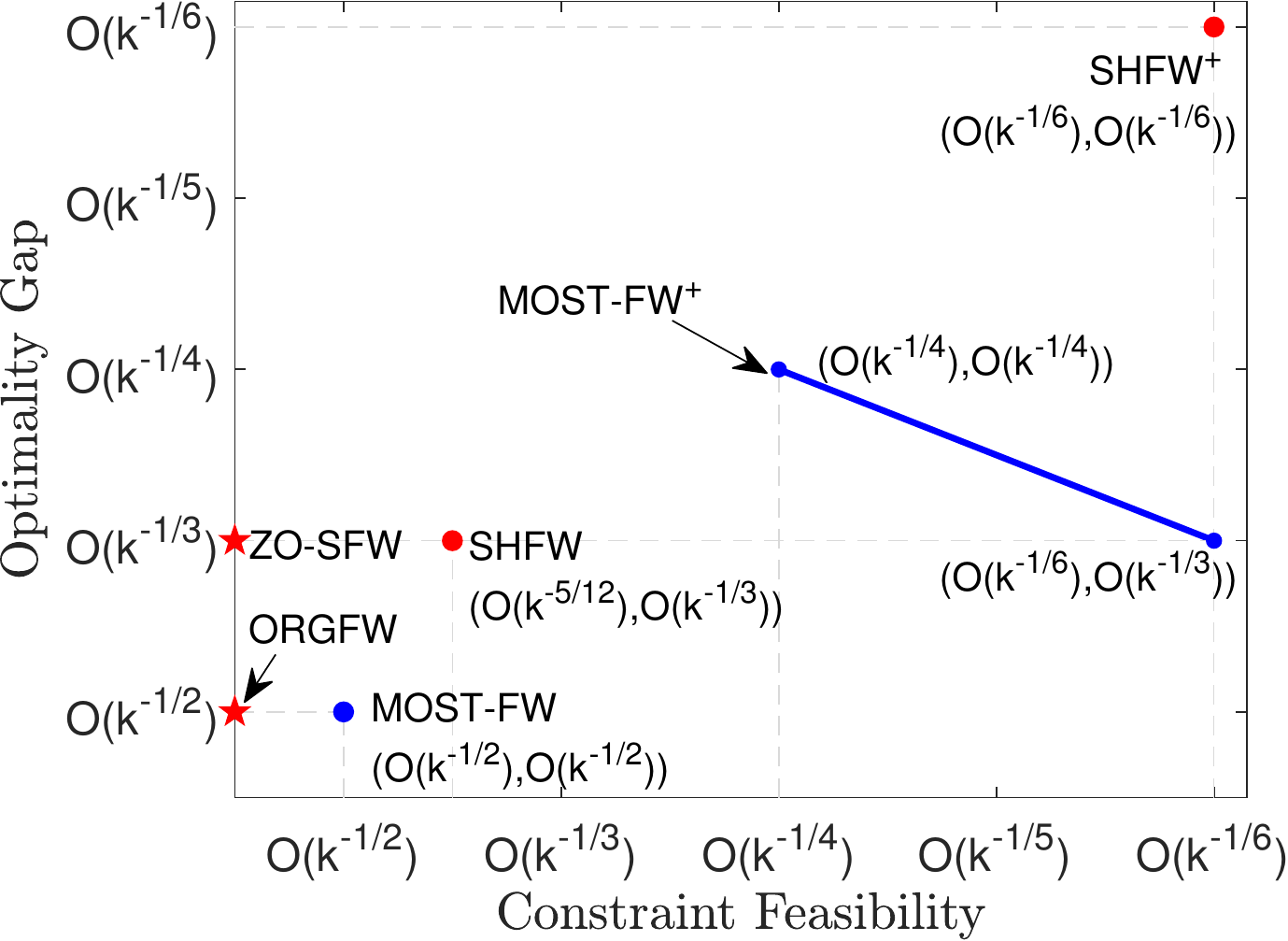}
		
	\end{subfigure} 	
	\caption{Comparison of convergence rates for different algorithms. Red color denotes existing state-of-the-art algorithms while blue color represents the proposed algorithms. Star is used to represent FW-algorithms for problems without additional constraints.}  
	\label{fig:objvscv}
\end{figure}

\noindent \textit{Remark 3:}  We can explore the dependency of convergence results on the smoothness parameter $\mu_k$ by setting $\mu_k = \frac{\mu_c}{(k+1)^{b}}$ with $b\in(1/6,1/4)$ in Theorem \ref{theorem2}. For this choice of $\mu_k$, we obtain an optimality gap of $\mathcal{O}(k^{b-\frac{1}{2}})$ and a constraint feasibility of $\mathcal{O}(k^{-\min\{b,\frac{1}{4}\}})$. Varying $b$ in the range $(1/6,1/4)$ one can observe the trade-off between optimality gap and constraint feasibility. This trade-off can be better understood from Fig.\ref{fig:objvscv}, which compares the convergence rates of all the discussed algorithms. For instance, on setting $b=1/6$, for the same order of constraint feasibility ($\mathcal{O}(k^{-1/6})$), the proposed algorithm \ms shows an improved optimality gap of $\mathcal{O}(k^{-1/3})$ over  $\mathcal{O}(k^{-1/6})$ of H1-SFW. Similarly on setting $b=1/4$, \ms  attains ($\mathcal{O}(k^{-1/4})$) optimality gap and while the constraint feasibility gets improved to $\mathcal{O}(k^{-1/4})$ compared to  $\mathcal{O}(k^{-1/6})$ of H1-SFW. 

\noindent \textit{Remark 4:}  Note that CGE in \eqref{coo:grad_apprx} is only a representative of various ZO gradient estimators. One can also use other estimators as long as it is ensured that the numerical gradient is asymptotically unbiased with $\rho_t\rightarrow 0$. For instance, we can use I-RDSA \cite{sahu2019towards} where only $\hat{m} < m$ directions are sampled at each iteration. Although such an approach reduces the SZO call  from $2m$ to $2\hat{m}$ at each iteration, it makes the final rate  proportional to $\frac{m}{\hat{m}}$. Hence, the total calls to SZO oracle to achieve $\eps$-optimal solution remains the same as in our case at $\mathcal{O}(m\eps^{-2})$ for $\mathcal{P}$ and $\mathcal{O}(m\eps^{-4})$ for $\mathcal{P}^+$.

\subsection{Convergence Analysis of T-\md and T-\ms}
We begin the analysis of trimmed variant with establishing a key lemma characterizing the error introduced due to trimming operation and then provide the convergence rates.
\begin{lemma}\label{cens-error}
	For Algorithm \ref{alg:3}, the error  introduced is upper bounded by the threshold $\tau_k$, i.e.,
	\begin{align}
	e_k:=\norm{\s_k-\v_k}< \tau_k \forall k \geq 0
	\end{align}
\end{lemma} 
\noindent\textbf{Proof}: According to the trimming strategy, we have  $\v_k=\s_k$ if $\norm{\s_k-\v_{k-1}}\geq \tau_k$ and $\v_k=\v_{k-1}$ if  $\norm{\s_k-\v_{k-1}}< \tau_k$. Therefore, in both cases, we have $e_k=\norm{\s_k-\v_k}< \tau_k$.


\begin{thm}\label{theorem3}
	Under Assumptions \ref{smoothness}-\ref{slater}, with $\eta_k = \frac{2}{k+1}$, $\gamma_k = \frac{1}{k}$, and $\rho_k \leq \frac{D}{\sqrt{m}(k+1)}$,  
	
	\noindent(a) for T-\md with  $\tau_k=\frac{\tau_0}{(k+1)^{1/2}}$ and $\mu_k = \frac{\mu_c }{\sqrt{k}}$, we have the following bounds
	
	\begin{align}
	\mbE[f(\x_{k})]-f(\x^{\star})\leq \mathcal{O}\left(\frac{1}{\sqrt{k}}\right)+\frac{8\tau_0D}{\sqrt{k}},
	\end{align}
	\begin{align}
	\mbE[\cD_{\cX}(\G\x_{k})]\leq  \mathcal{O}\left(\frac{1}{\sqrt{k}}\right)+\frac{4\sqrt{\tau_0D}}{\sqrt{k}},
	\end{align}
	\noindent(b) and for T-\ms with  $\tau_k=\frac{\tau_0}{(k+1)^{1/4}}$ and $\mu_k =\frac{\mu_c}{(k+1)^{1/4}}$, we have the following bounds
	\begin{align}
	\mbE[f(\x_{k})]-f(\x^{\star})\leq \mathcal{O}\left(\frac{1}{k^\frac{1}{4}}\right)+\frac{8\tau_0D}{k^\frac{1}{4}}
	\end{align}
	\begin{align}
	\mbE[\cD_{\cX}(\G\x_{k})]\leq  \mathcal{O}\left(\frac{1}{k^\frac{1}{4}}\right)+\frac{6\sqrt{\tau_0D}}{k^\frac{1}{4}}.
	\end{align}
	
\end{thm}
\textbf{Proof}: The proof  is provided in Sec. \ref{proof:theorem3} of the supplementary material. 
Theorem \ref{theorem3} establishes an upper bound
on the expected suboptimality and constraint feasibility for the iterates generated by T-\md (part(a)) and T-\ms (part(b)). Note that, as compared to the non-trimmed version, these convergence rates are reduced by the factor of $\tau_0$. Hence, the SFO oracle complexity for T-\md and T-\ms will be $\O\left(\frac{(1+\tau_0)^2}{\epsilon^2}\right)$ and $\O\left(\frac{(1+\tau_0)^4}{\epsilon^4}\right)$, respectively. Although it is not possible to directly quantify the LMO complexity for trimmed case, we expect it to be a fraction of SFO complexity. In fact, experimentally (in Section \ref{sec:application}) we observe a significant reduction in LMO calls (more than $30\%$) when compared  to the non-trimmed version while maintaining almost similar accuracy. 
\begin{figure}
	\centering
	\setcounter{subfigure}{0}
	\begin{subfigure}{0.35\columnwidth}		\includegraphics[width=\linewidth, height = 0.7\linewidth]{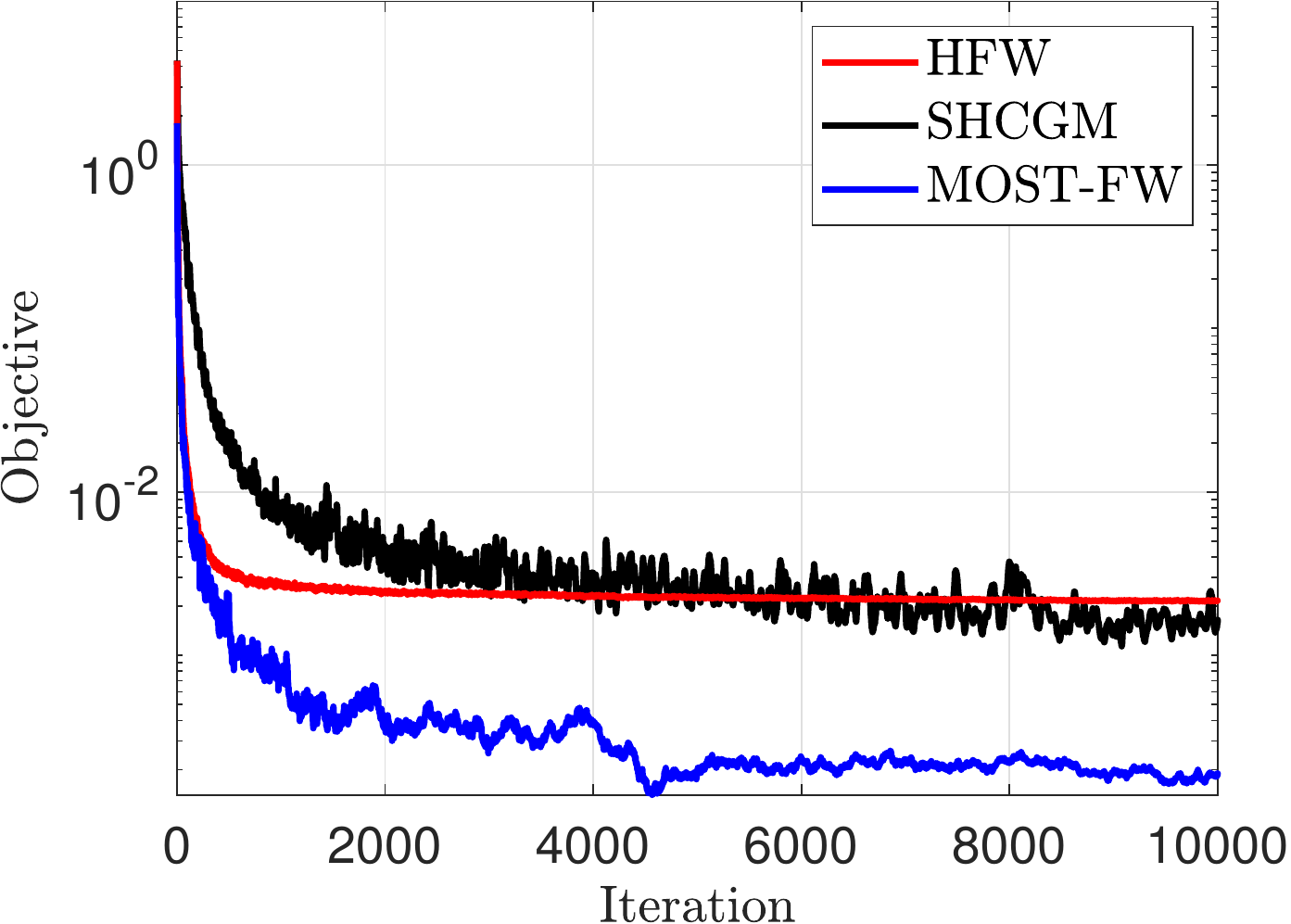}
		\caption{Objective}
		\label{ex2obj} 
	\end{subfigure}
	\begin{subfigure}{0.35\columnwidth}
		\includegraphics[width=\linewidth, height = 0.7\linewidth]{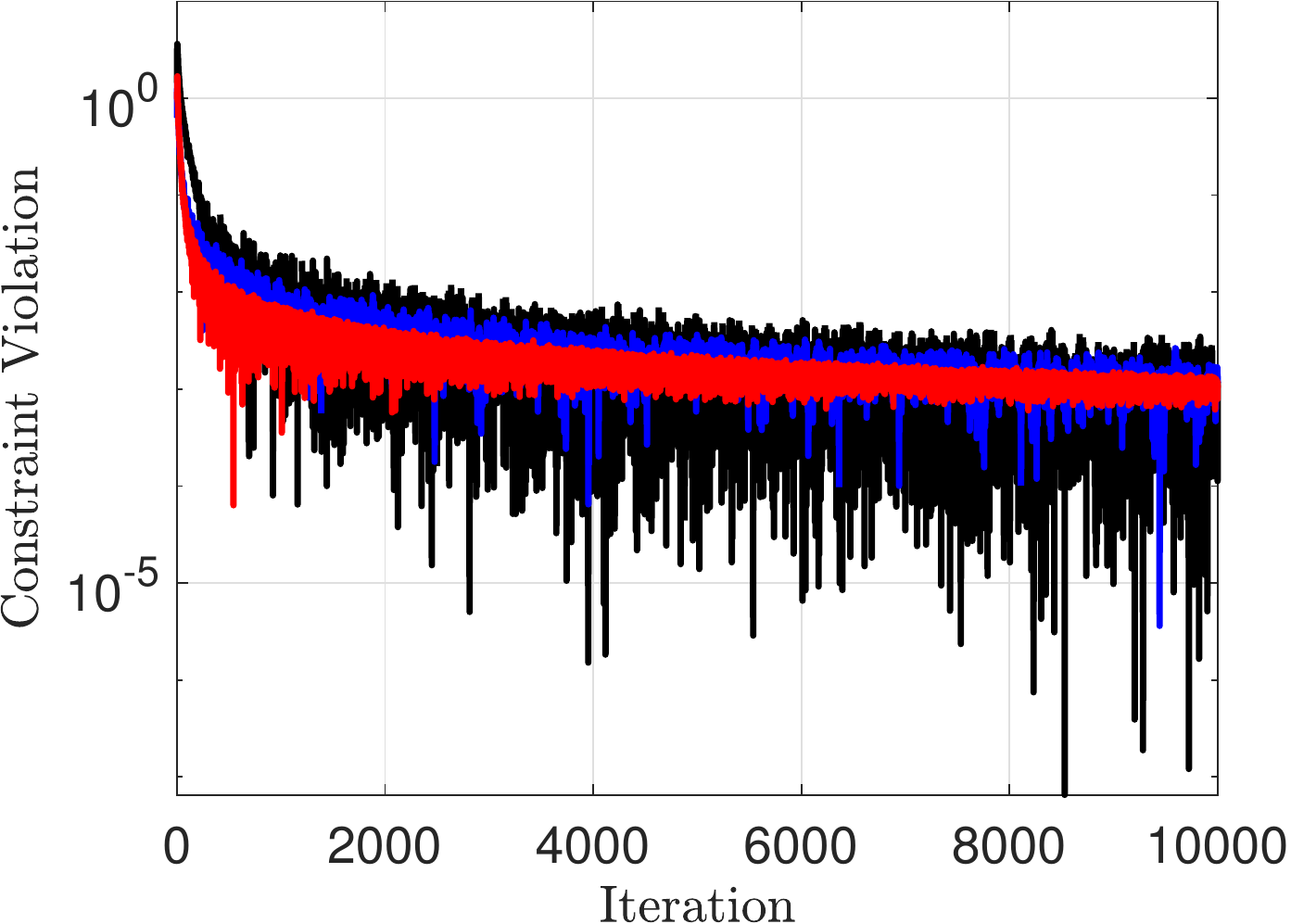}
		\caption{Constraint Violation}
		\label{ex2cv1} 
	\end{subfigure}
	\begin{subfigure}{0.35\columnwidth}		\includegraphics[width=\linewidth, height = 0.7\linewidth]{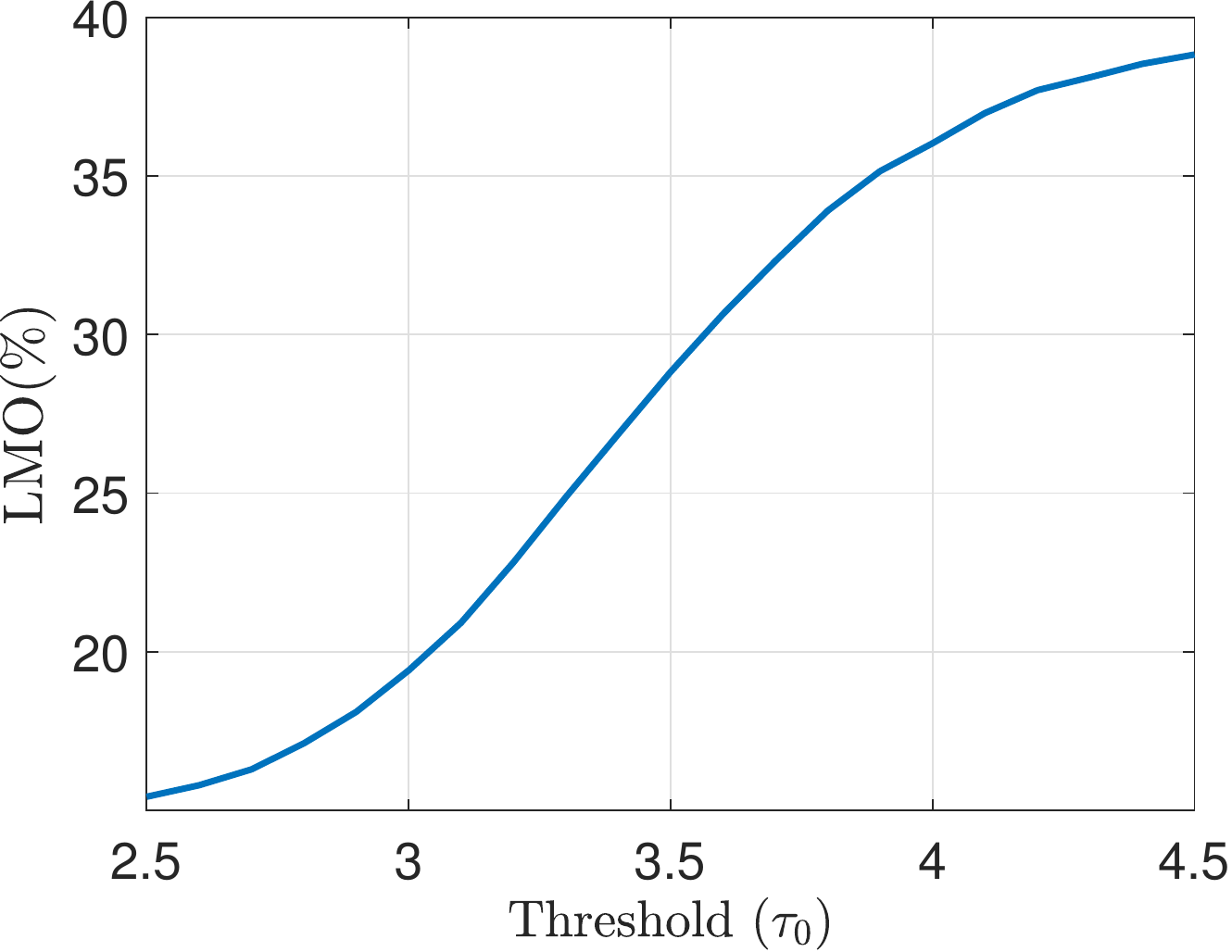}
		\caption{T-MOST-FW}
		\label{LO1_1} 
	\end{subfigure}
	\begin{subfigure}{0.35\columnwidth}
		\includegraphics[width=\linewidth, height = 0.7\linewidth]{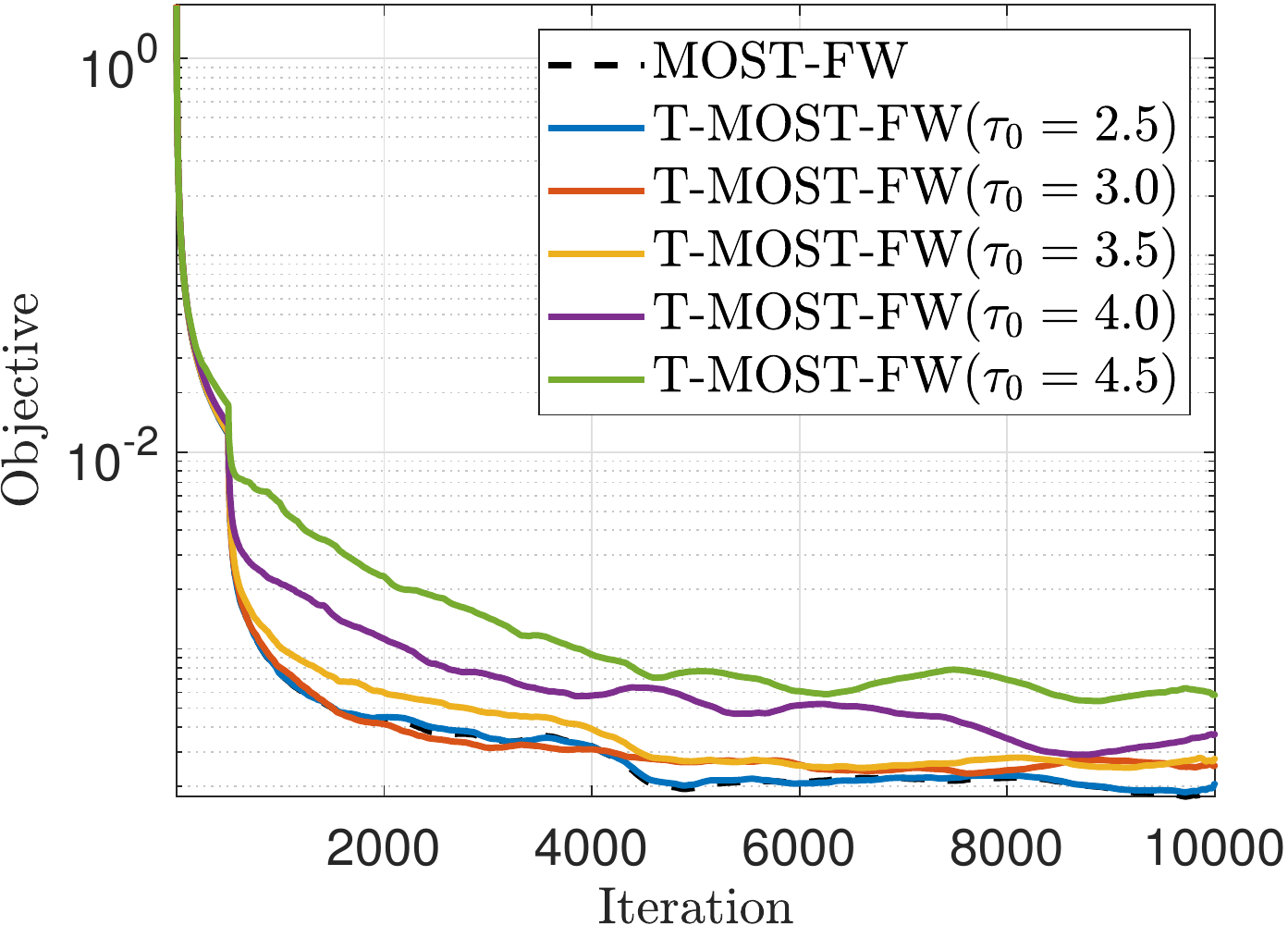}
		\caption{Objective}
		\label{OT1_1} 
	\end{subfigure}
	\caption{Sparse matrix estimation (Sec.\ref{cov-mat}): (a)-(b): Comparison of the proposed algorithm MOST-FW with HFW and SHCGM, (c): percentage of LMO (compared to MOST-FW) skipped  and (d) convergence performance of  T-MOST-FW at different threshold values ($\tau_0$).}
	\vspace{-0mm}
	\label{fig:matrix_est}
\end{figure}

\section{Numerical Experiments}\label{sec:application} 
In this section, 
we compare the performance of proposed algorithms \md and \ms with that of  HFW\cite{yurtsever2018conditional}, SHCGM \cite{locatello2019stochastic} and H1-SFW \cite{vladarean2020conditional}. We set $\gamma_k = \frac{1}{k}$, $\eta_k = \frac{2}{k+1}$, and $\mu_k=\frac{\mu_c}{\sqrt{k}}$ for \md and $\gamma_k = \frac{1}{k}$, $\eta_k = \frac{2}{k+1}$, and $\mu_k=\frac{\mu_c}{k^{1/4}}$ for \ms, as dictated by theory.. Further, we set $\tau_k$ as $\frac{\tau_0}{(k+1)^{1/2}}$ and $\frac{\tau_0}{(k+1)^{1/4}}$ for T-\md and T-\ms, respectively.  The parameter $\mu_c$ and $\tau_0$ are tuned to yield the best possible performance. 
We plot the objective and constraint violation for different applications to highlight the improved convergence bounds obtained theoretically.  Further, to understand the effect of the trimming scheme, we  plot the total number of times the trimmed variant of the algorithms skipped the LMO calls (in $\%$) for different values of threshold $\tau_0$ as compared to the total LMO calls required by the non-trimmed version ($\tau_0= 0$) of the algorithm. We also  demonstrate the
potential advantage of the trimmed version over the original methods by plotting the evolution of optimality gap at different threshold values.  

\subsection{Sparse matrix estimation}\label{cov-mat}
We consider the problem of estimating a sparse covariance matrix, given independent samples of a Gaussian random vector, observed in a sequential fashion \cite{richard2012estimation}. The problem arises in a number of disciplines, such as physics,  finance, machine learning, genome sequencing, and portfolio optimization \cite{zhao2014robust,deshmukh2020improved}. The problem can be formulated as
\begin{align}
&\min_{\X \in \cc} f(\X):=\mbE\norm{\X-\w \w^\T}^2_F \label{obj_p1}\\
&\text{s.t.}\; \norm{\vect{\X}}_1\leq \alpha,\label{cons_p1}
\end{align}
where $\cc := \{\X \succeq 0, \text{tr}(\X)\leq K \}$ and the expectation is over random vector $\w$. For the experiments, we consider the setup in  \cite{locatello2019stochastic}. The underlying covariance matrix $\W \in \Rn^{1000 \times 1000}$ is generated as $\W = \sum_{i=1}^{10}\boldsymbol{\psi}_i\boldsymbol{\psi}_i^\mathsf{T}$ where the entries of $\boldsymbol{\psi}_i$ are drawn uniformly at random from $[-1, 1]$. At each $t$, we observe $\w_t \sim \mathcal{N}(\mathbf{0},\W)$ and the goal is to estimate $\W$ by solving \eqref{obj_p1}-\eqref{cons_p1}, which matches the template in \md. Observe that the projection over $\cc$ is difficult, but projection over the $\ell_1$ norm ball is easy. For this problem, we set $\alpha = \text{tr}(\W)$ and $K = \norm{\vect{\W}}_1$. 
We compare the performance of the deterministic HFW, its online version SHCGM, and the proposed MOST-FW algorithms, the latter two (being stochastic) uses a mini-batch of 200 data points. 
\begin{figure*}[t!]
	\centering
	\setcounter{subfigure}{0}
	\begin{subfigure}{0.32\columnwidth}		\includegraphics[width=\linewidth, height = 0.7\linewidth]{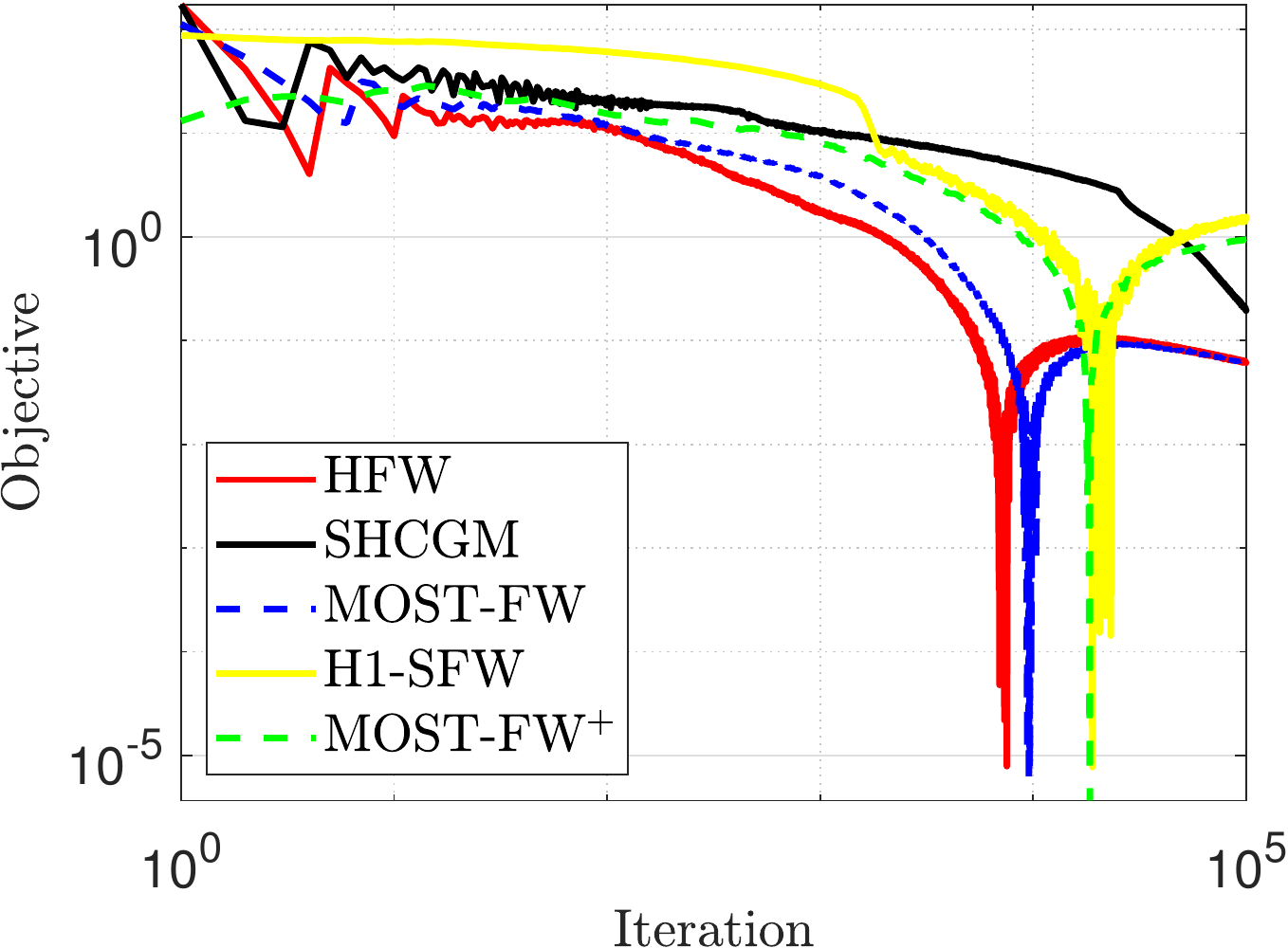}
		\caption{Objective}
		\label{ex1obj}
	\end{subfigure} 	
	\begin{subfigure}{0.32\columnwidth}		\includegraphics[width=\linewidth, height = 0.7\linewidth]{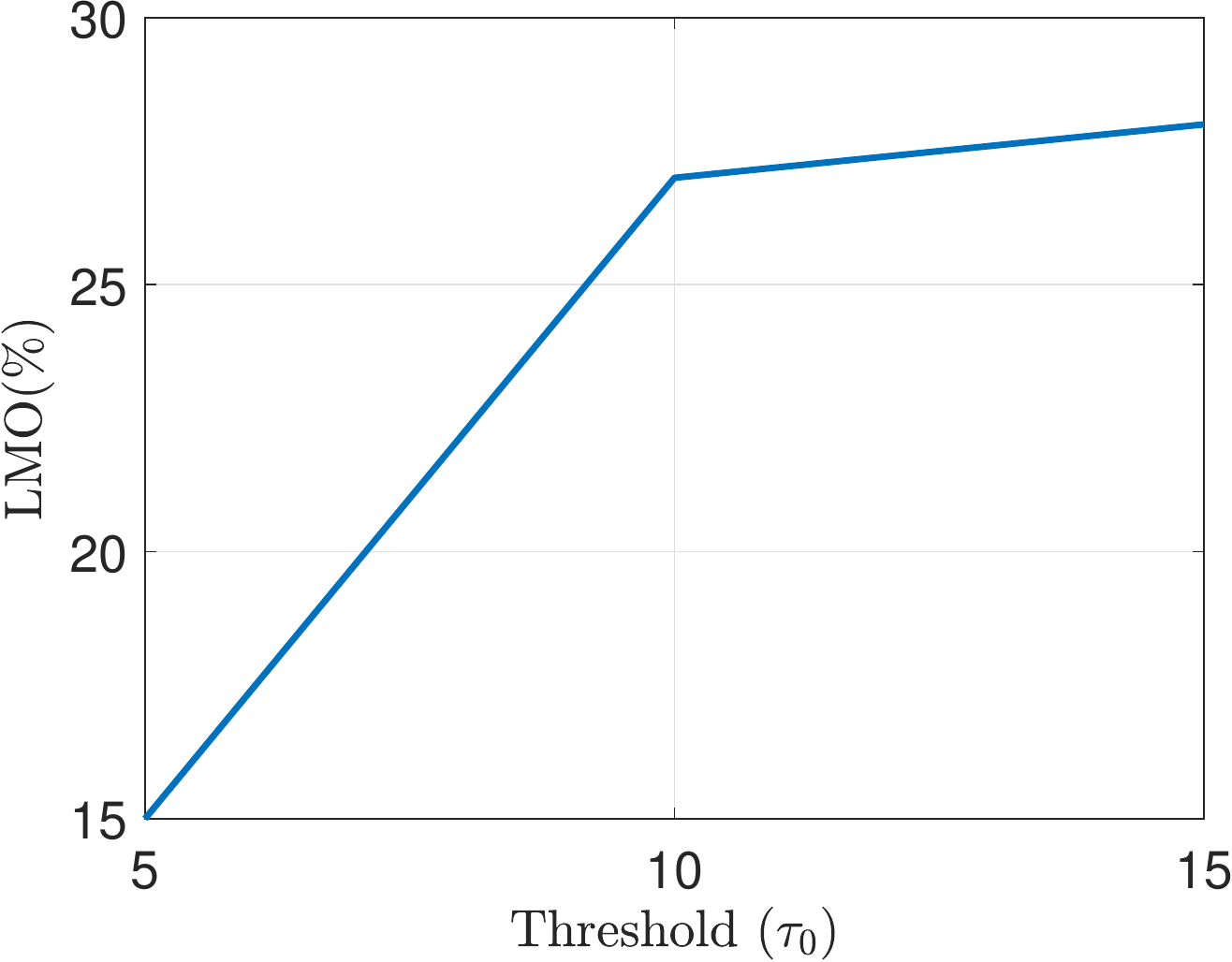}
		\caption{T-MOST-FW}
		\label{LO1_2}
	\end{subfigure} 	
	\begin{subfigure}{0.32\columnwidth}
		\includegraphics[width=\linewidth, height = 0.7\linewidth]{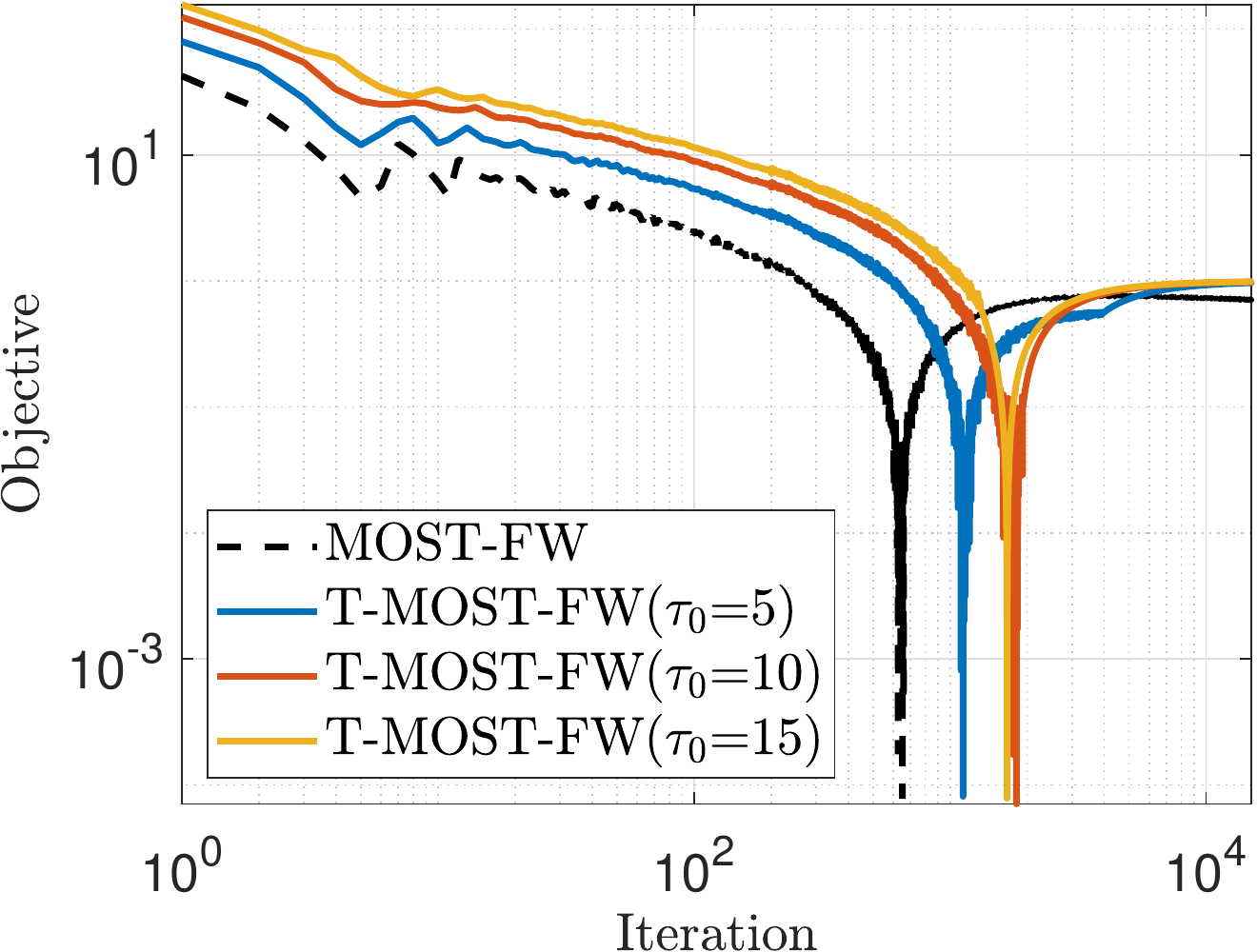}
		\caption{Objective}
		\label{OT1_2}
	\end{subfigure}
	\begin{subfigure}{0.32\columnwidth}
		\includegraphics[width=\linewidth, height = 0.7\linewidth]{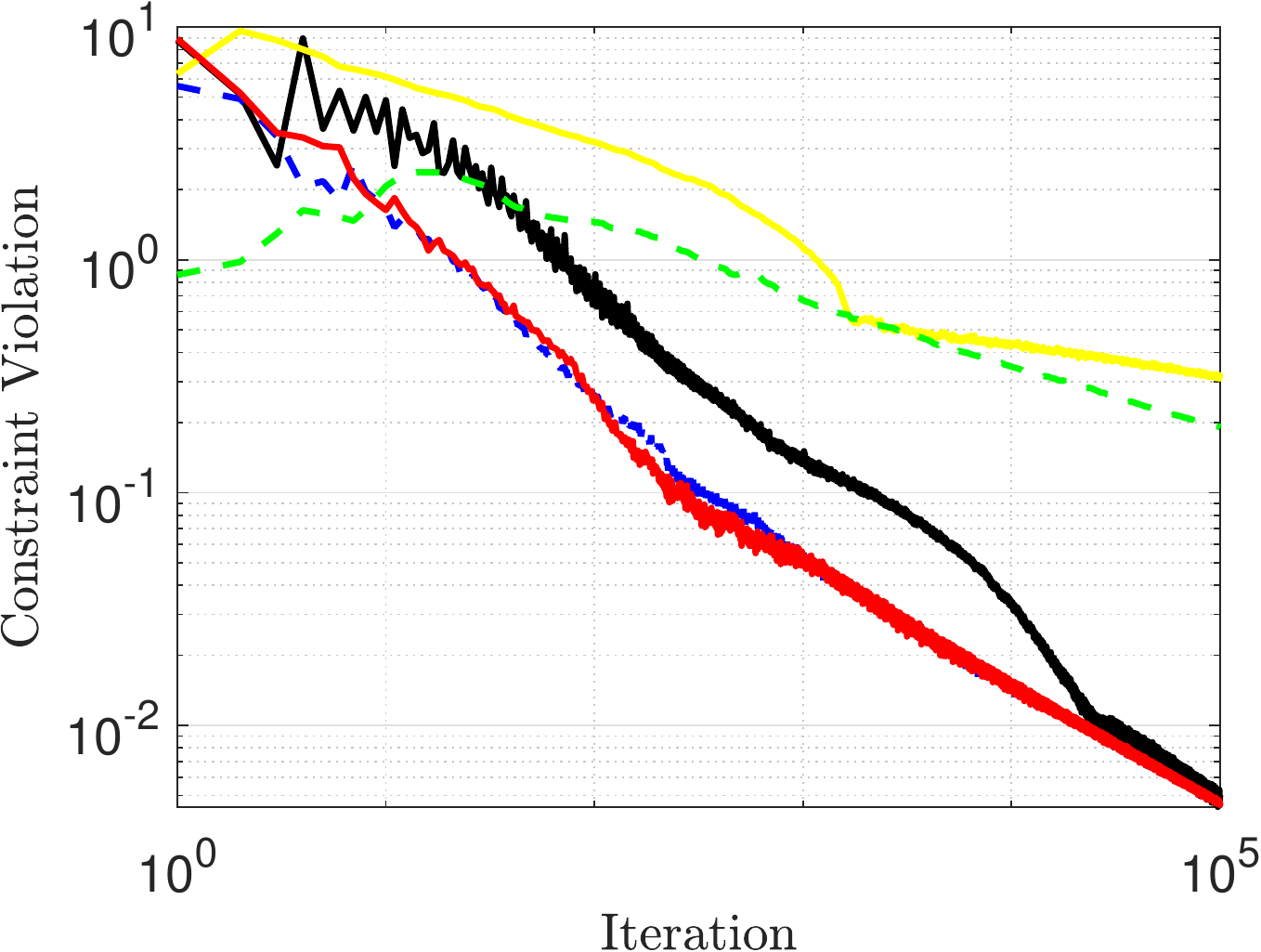}
		\caption{Constraint Violation}
		\label{ex1cv1}
	\end{subfigure}
	\begin{subfigure}{0.32\columnwidth}		\includegraphics[width=\linewidth, height = 0.7\linewidth]{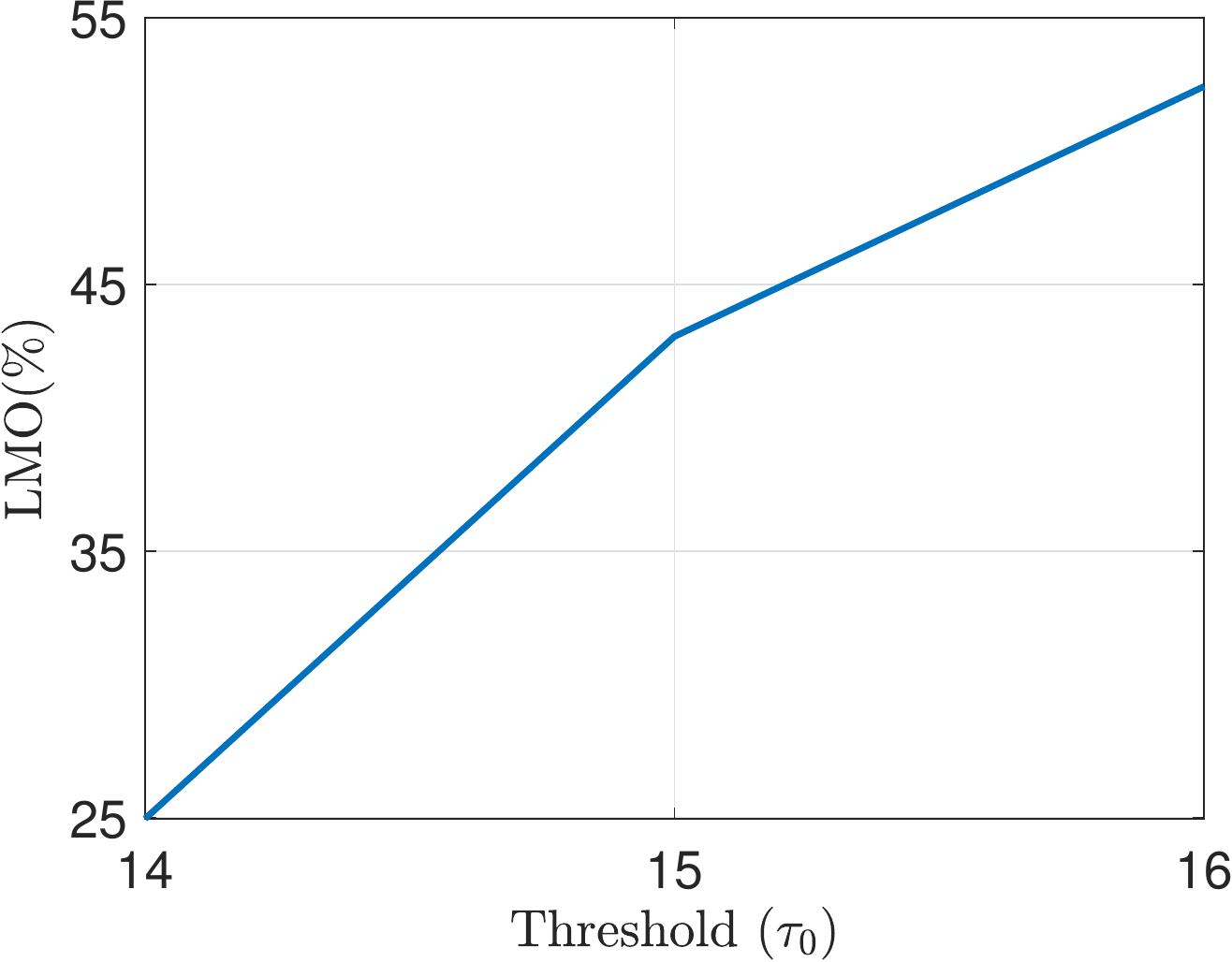}
		\caption{T-MOST-FW$^+$}
		\label{LO2_2}
	\end{subfigure} 	
	\begin{subfigure}{0.32\columnwidth}
		\includegraphics[width=1.08\linewidth, height = 0.8\linewidth]{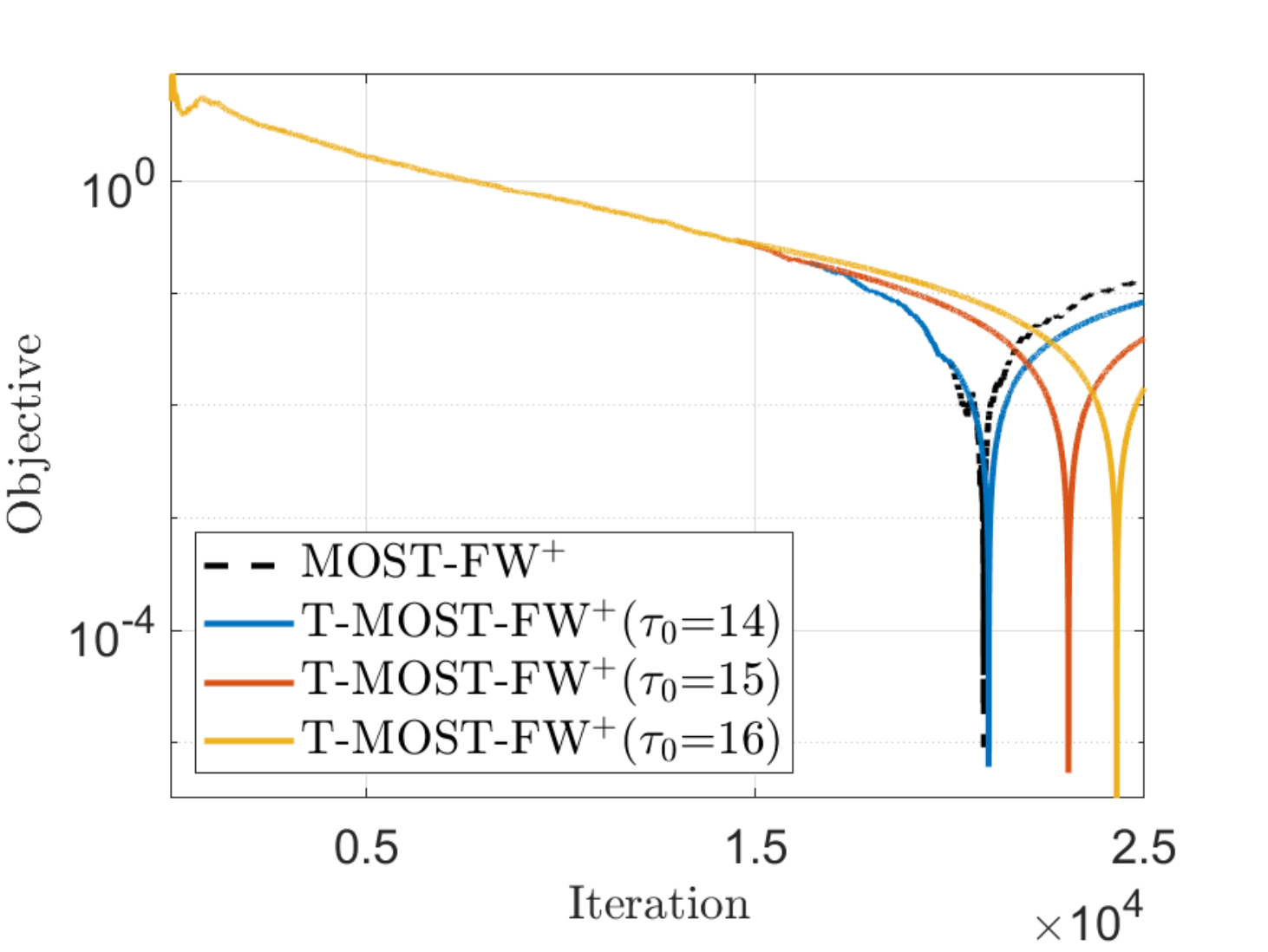}
		\caption{Objective}
		\label{OT2_2}
	\end{subfigure}
	\caption{K-means clustering (Sec.\ref{clustering}): Comparison of (a) Optimality Gap (d) Constraint violation of the proposed algorithms.  Percentage of LMO calls skipped at different threshold values  (b) by T-MOST-FW (compared to MOST-FW) and (e) by T-MOST-FW$^+$ (compared to MOST-FW$^+$).  Convergence performance (c) for T-MOST-FW and (f) for T-MOST-FW$^+$ at different threshold values.  }
	\label{fig:kmeans_trim}
\end{figure*}

We set $\mu_c=1$ for MOST-FW while the parameters of HFW and SHCGM are kept the same as in \cite{yurtsever2018conditional},\cite{locatello2019stochastic}, respectively. We run all three algorithms for $10^4$ iterations and analyze their performance in terms of  objective convergence as $\norm{\X-\W}^2_F/\norm{\W}^2_F$ and constraint violation as $\max(\norm{\X}_1-\alpha,0)/\alpha$. 

The convergence plots for this experiment are shown in Fig.\ref{fig:matrix_est}. HFW, being a deterministic algorithm, converges fast but gets saturated at some accuracy as it uses the same datapoints at all the iterations. Advantages of \md over both the compared algorithms are evident from the plots of optimality gap Fig.\ref{ex2obj} and constraint violation Fig.\ref{ex2cv1}, which converges much faster and to a better accuracy level. 
Next, we analyze the performance of T-\md at different threshold values $\tau_0$ and plot the variation of LMO($\%$) in Fig.\ref{LO1_1}.  As expected, the algorithm frequently skips calls to the LMO with the increase in the threshold value. We proceed to study the effect of trimming on the convergence in Fig.\ref{OT1_1}  Observe that though trimming affects the convergence, for properly chosen threshold we can get comparable  performance with that of the non-trimmed version. For instance, at $\tau_0=3.5$, T-\md skips around $28\%$ of LMO calls while maintaining almost the same performance as \md.
\subsection{Clustering via semidefinite relaxation}\label{clustering}
\begin{figure*}
	\centering
	\setcounter{subfigure}{0}
	\begin{subfigure}{0.32\columnwidth}		\includegraphics[width=\linewidth, height = 0.7\linewidth]{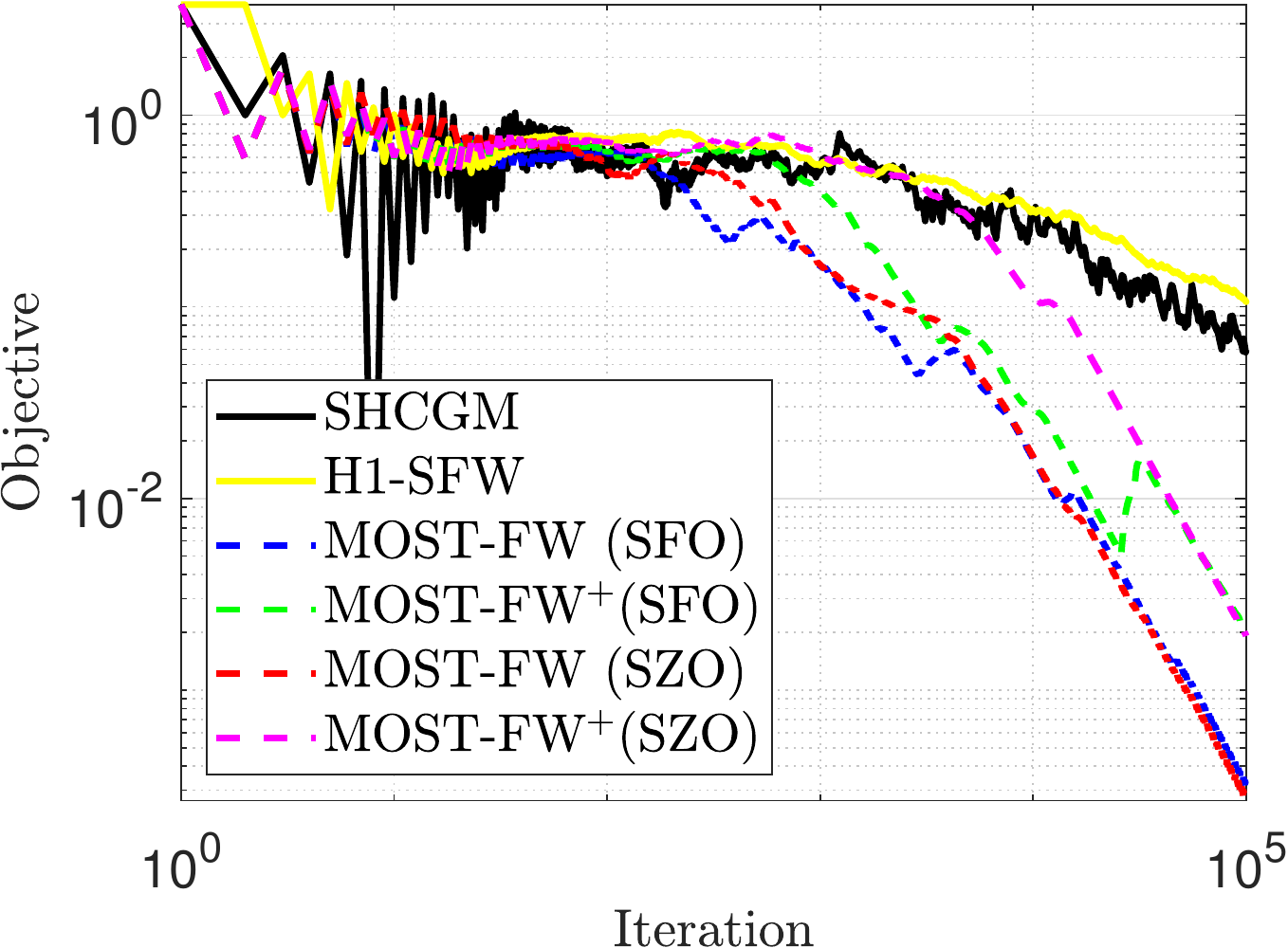}
		\caption{Objective}
		\label{scut_obj}
	\end{subfigure}
	\begin{subfigure}{0.32\columnwidth}		\includegraphics[width=\linewidth, height = 0.7\linewidth]{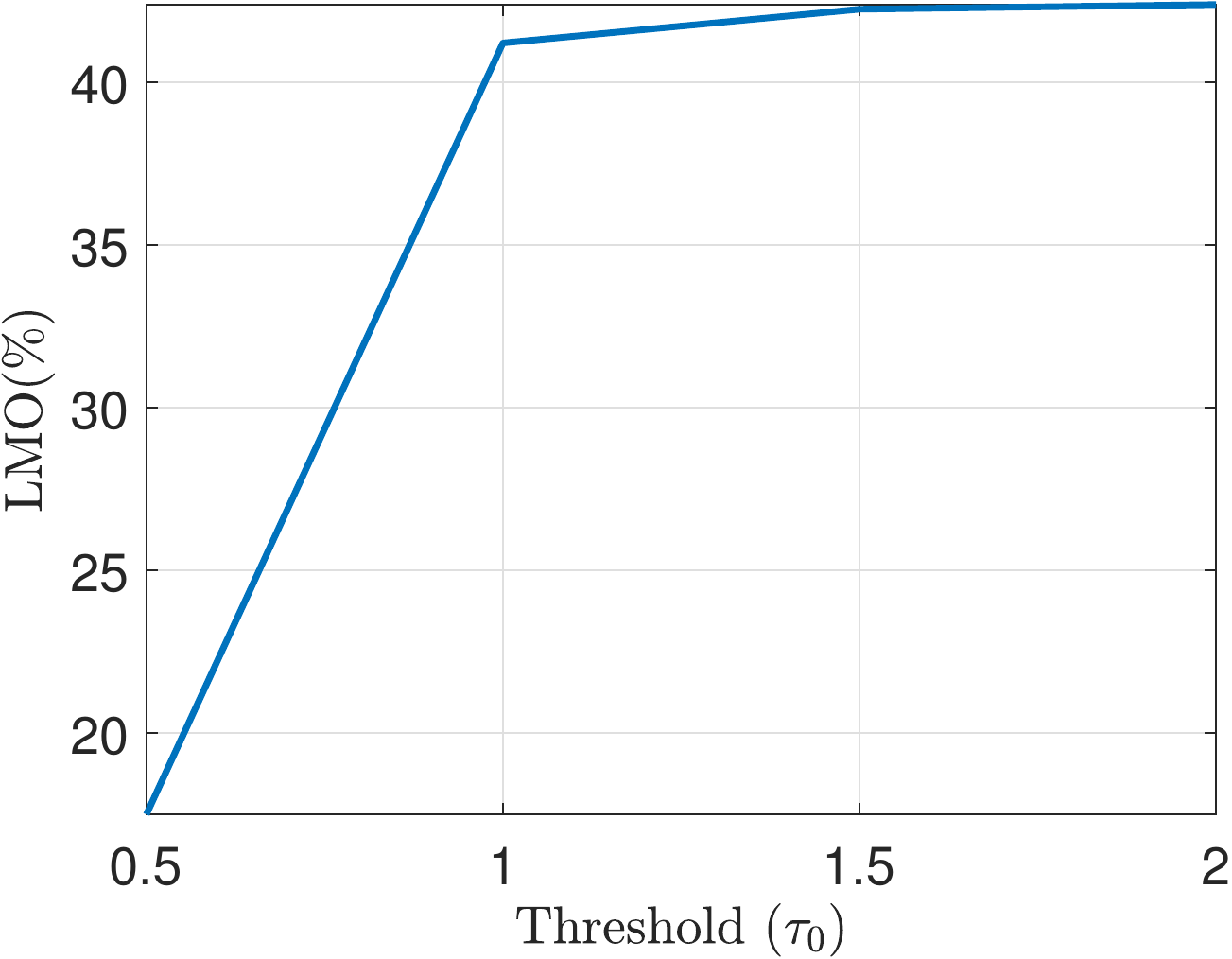}
		\caption{T-MOST-FW}
		\label{LO1_3}
	\end{subfigure}
	\begin{subfigure}{0.32\columnwidth}
		\includegraphics[width=\linewidth, height = 0.7\linewidth]{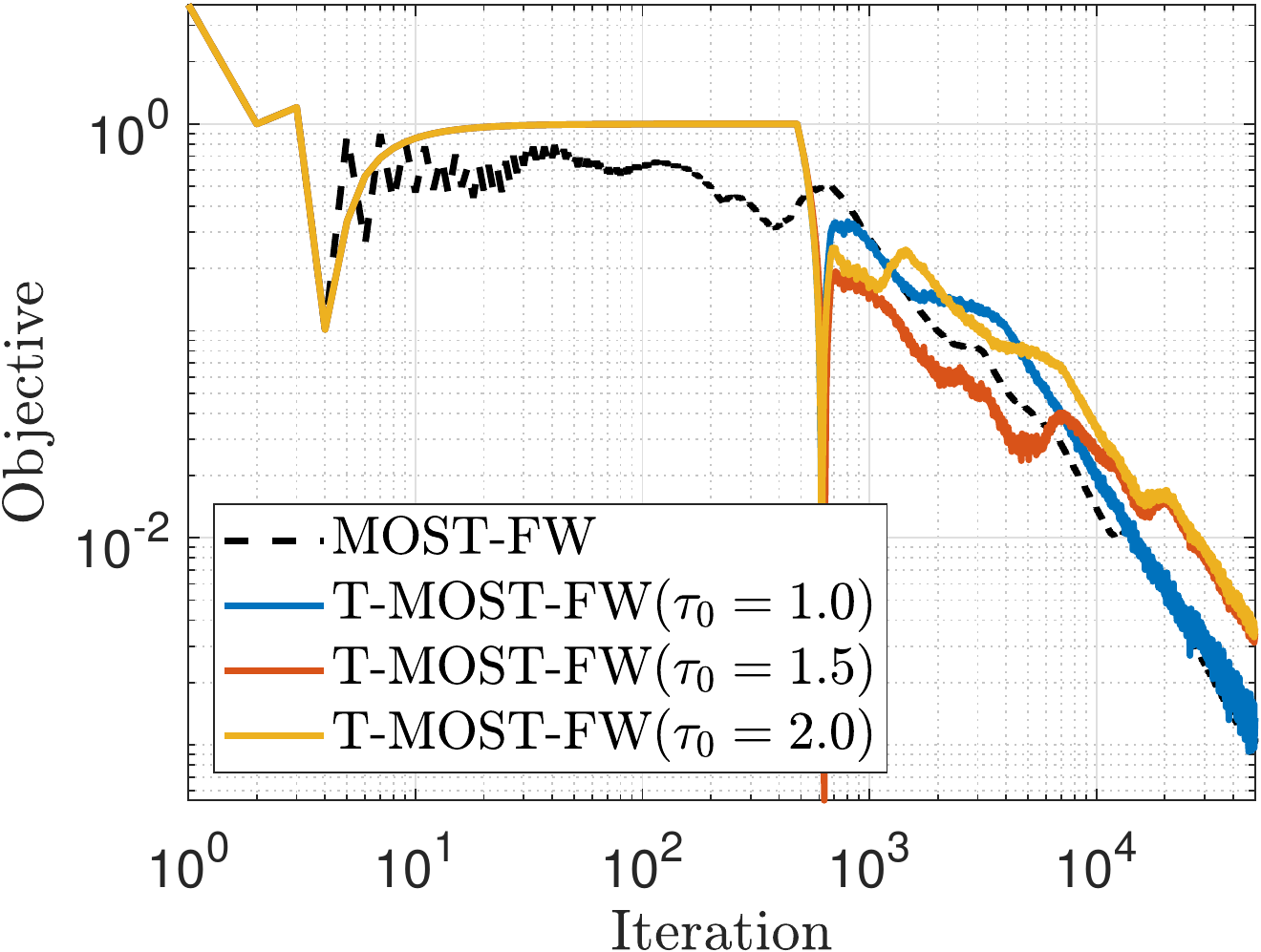}
		\caption{Objective}
		\label{OT1_3}
	\end{subfigure}
	\begin{subfigure}{0.32\columnwidth}
		\includegraphics[width=\linewidth, height = 0.7\linewidth]{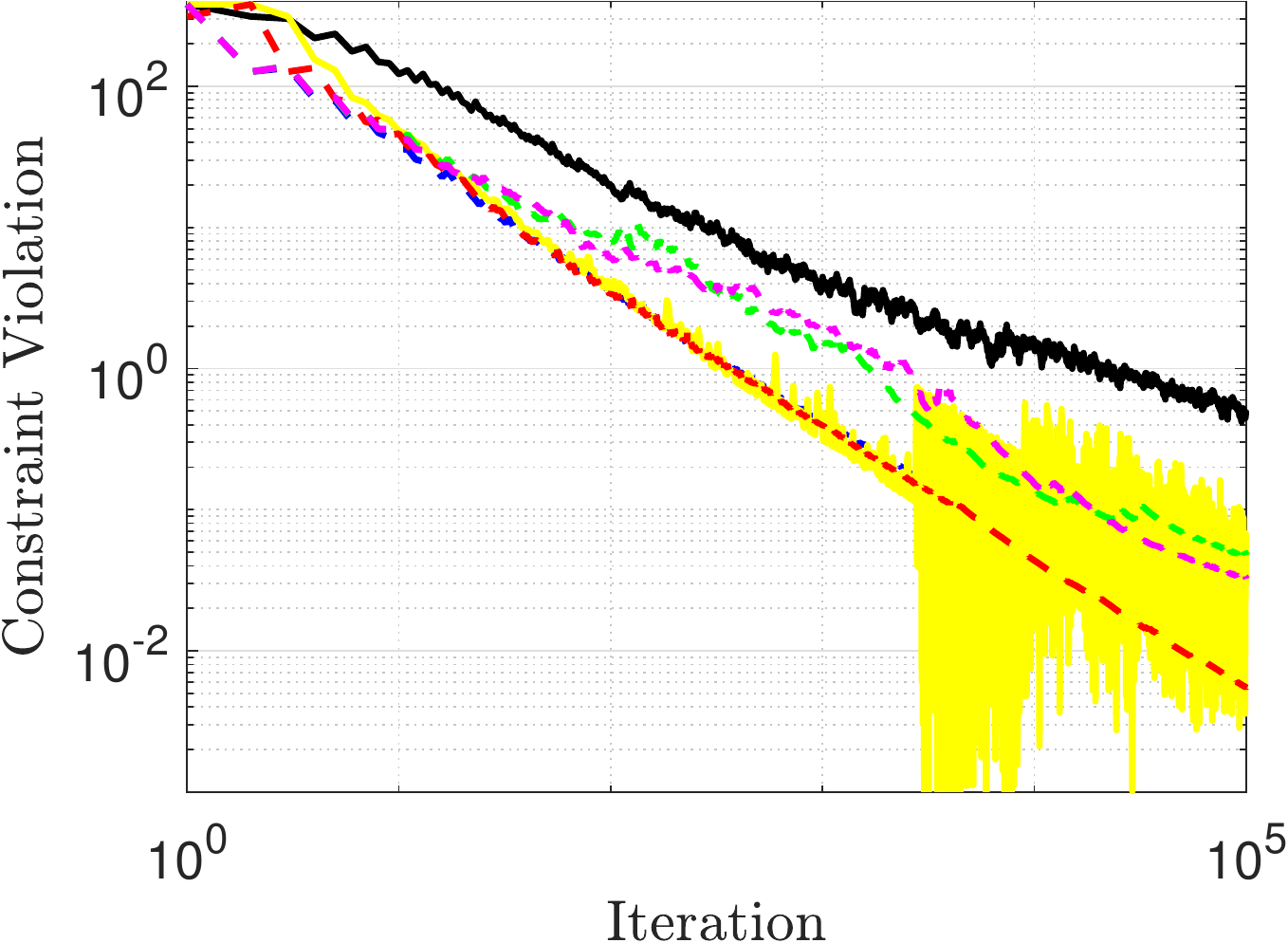}
		\caption{Constraint Violation}
		\label{scut_cv1}
	\end{subfigure}
	\begin{subfigure}{0.32\columnwidth}		\includegraphics[width=\linewidth, height = 0.7\linewidth]{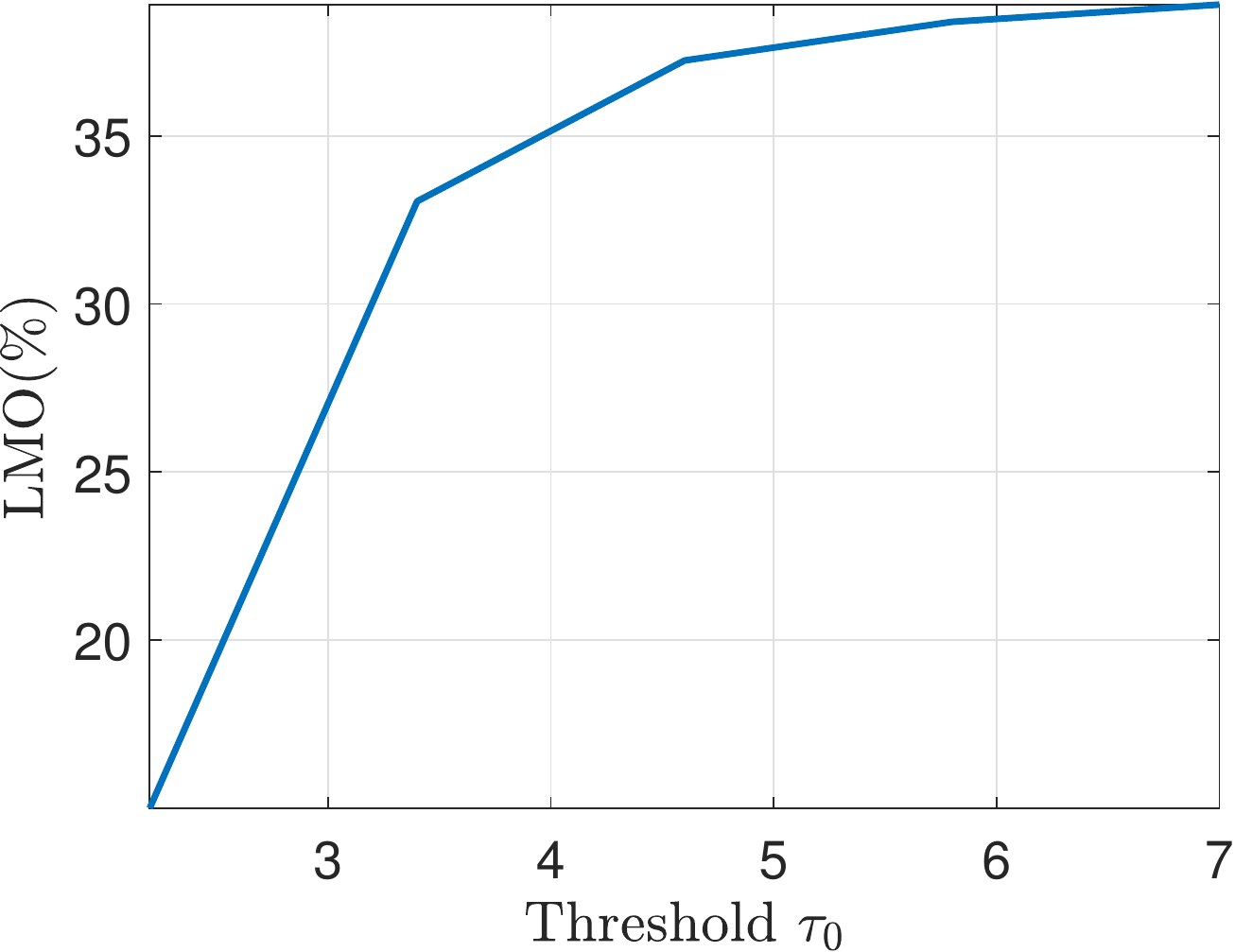}
		\caption{T-MOST-FW$^+$}
		\label{LO2_3}
	\end{subfigure}
	\begin{subfigure}{0.32\columnwidth}
		\includegraphics[width=\linewidth, height = 0.7\linewidth]{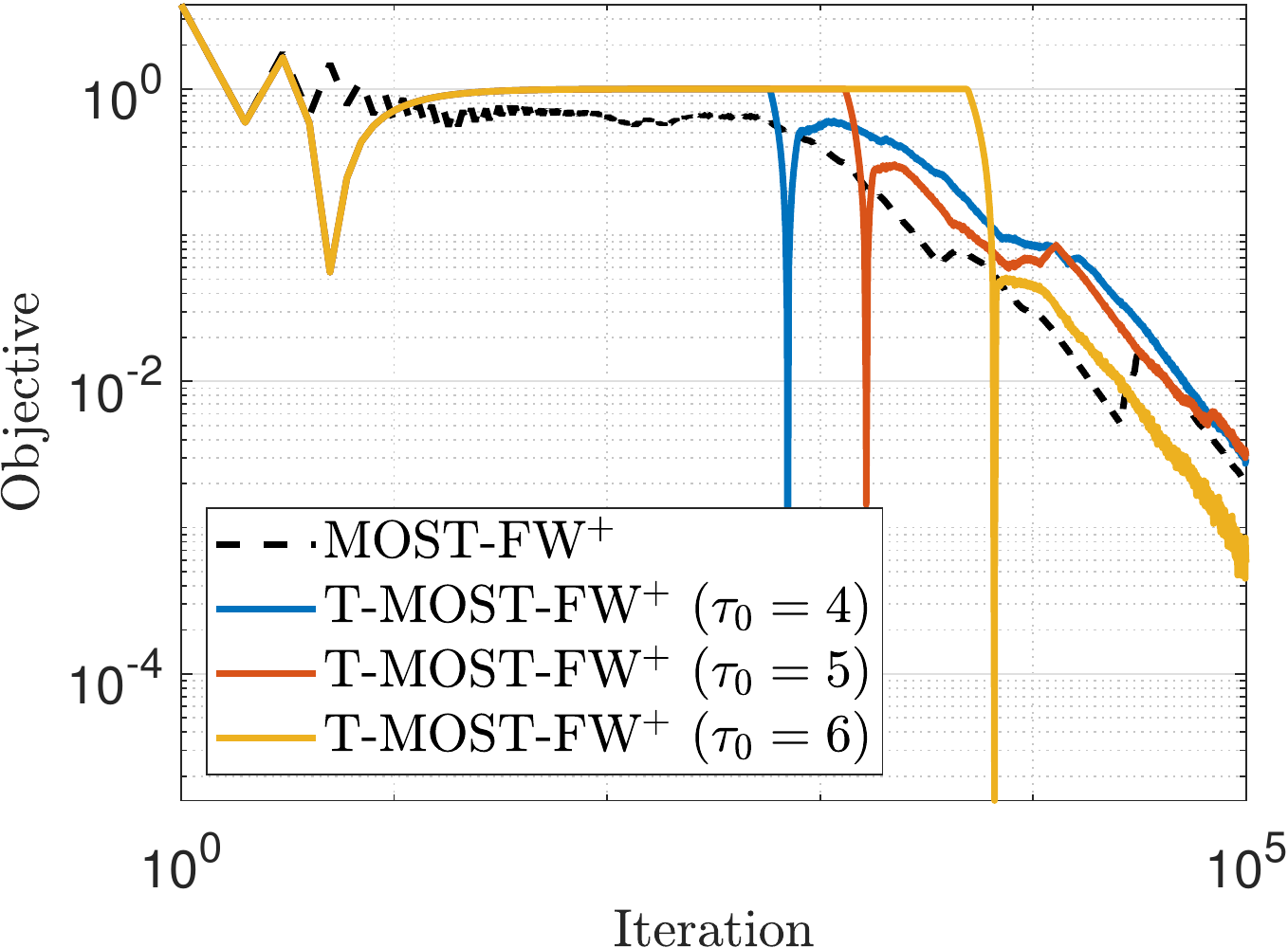}
		\caption{Objective}
		\label{OT2_3}
	\end{subfigure}
	\caption{Uniform Sparsest Cut (Sec.\ref{sparse-cut}): Comparison of (a) Optimality Gap (d) Constraint violation of the proposed algorithms. Percentage of LMO calls skipped at different threshold values  (b) by T-MOST-FW (compared to MOST-FW) and (e) by T-MOST-FW$^+$ (compared to MOST-FW$^+$).  Convergence performance (c) for T-MOST-FW and (f) for T-MOST-FW$^+$ at different threshold values.}
	\vspace{-0mm}
	\label{fig:scut_trim}
\end{figure*}
In this experiment, we  will consider the $K$-means clustering problem defined in \eqref{obj_p2}-\eqref{cons_p2} which is also considered in the previous related works
\cite{mixon2016clustering,yurtsever2018conditional,locatello2019stochastic}. 

Observe that \eqref{obj_p2}-\eqref{cons_p2} is of the form in $\pd$ when only the entries of $\M$ corresponding to a randomly selected subset $\Omega_t \subset \Omega$ are revealed at every iteration. The $K$-means problem can therefore be solved using the SHCGM algorithm from \cite{locatello2019stochastic} and the proposed MOST-FW algorithm. Further, if a randomly selected subset of the constraints are imposed at every iteration, we obtain the formulation in $\ps$, allowing us to use the H1-SFW algorithm \cite{vladarean2020conditional} and the proposed \ms algorithm. It is remarked that the projection over a subset constraints is simpler, and therefore the formulation in $\ps$ is better suited to larger scale problems. 

For the experiment, we use $N=1000$ samples from MNIST dataset \cite{lecun2010mnist} with $K=10$ clusters. For \md, \ms, SHCGM, and H1-SFW, we estimate objective gradient using only $1\%$ of the data at each iteration while HFW being deterministic in nature, uses full data. In addition, for \ms and H1-SFW, we use only $1\%$ of the randomly  sampled constraints at each iteration.  
We tune $\mu_c$ to achieve the best performance for all the algorithms and finally set $\mu_c=10$ for \md and $\mu_c=2.75$ for \ms while for HFW, SHCGM and H1-SFW we use the same parameters as in \cite{yurtsever2018conditional,locatello2019stochastic,vladarean2020conditional}, respectively. We run all the algorithms for $10^5$ iterations and analyze their performance in terms of objective convergence as $|f(\X)-f^{\star}|/|f^{\star}|$ where $f^{\star}=f(\X^{\star})$ is derived from the ground truth\footnote{ The ground truth is obtained using kmeans-sdp of \cite{mixon2016clustering}; see also \href{https://github.com/solevillar/kmeans_sdp}{https://github.com/solevillar/kmeans\_sdp}} and constraint violation as the sum of violation of both the constraints of \eqref{cons_p2}, i.e., $\norm{\X\textbf{1}-\textbf{1}}/\norm{\textbf{1}}$ and $\norm{\X-\Pi_{[\X]_{i,j}\geq 0}(\X)}_F$. 
The convergence results are shown in Fig. \ref{ex1obj} and \ref{ex1cv1}. It can be observed that \md shows the best practical performance, very similar to the deterministic algorithm HFW, although it uses only $1\%$ of the data at each iteration. Also, \ms shows slightly better performance compared to H1-SFW. These results support our theoretical findings.

Next, to analyze the performance of trimmed variants, in our experiment, we use $N=2000$ samples while keeping all other settings the same as above. The experimental results are shown in Fig.\ref{fig:kmeans_trim} (b-c) and Fig.\ref{fig:kmeans_trim} (e-f). Observe that both the trimmed algorithms skip more LMO calls with the increase in the threshold value. Although T-\md skips only around $15\%$ of LMO calls at $\tau_0=5$, we observed slightly improved performance for T-\ms that skips  $25\%$ of LMO calls (at $\tau_0=14$) while maintaining almost the same convergence performance as \ms.  
\subsection{Uniform Sparsest Cut Problem}\label{sparse-cut}
Consider the problem of graph partitioning where the target is to divide the graph into two or more sub-graphs by cutting the smallest number of edges. The problem arises in the design of divide-and-conquer algorithms for a number of other problems including, communications in distributed
networks,  cluster analysis and, machine learning \cite{chatziafratis2018hierarchical, bonsma2012complexity}.


Consider a graph $\mathcal{G}=(\mathcal{V},\mathcal{E})$, where $\mathcal{V}$ is a set of $d$ vertices, and $\mathcal{E}$ is a set of edges. Let $\mathcal{A}, \hat{\mathcal{A}} \subset \mathcal{V}$ denote two disjoint sets of vertices and let $E(\mathcal{A},\hat{\mathcal{A}})$ denote the number of edges between these sub-graphs. The uniform sparsest cut problem considered in \cite{arora2009expander} aims to find the cut $(\mathcal{A},\hat{\mathcal{A}})$ which minimizes 
\begin{align}
\frac{|E(\mathcal{A},\hat{\mathcal{A}})|}{|\mathcal{A}||\hat{\mathcal{A}}|},
\end{align}
where the $\abs{\cdot}$ operator returns the cardinality of the corresponding set. However, finding such a partition is NP-hard, and a number of approximation algorithms for solving it exist   \cite{leighton1999multicommodity, shmoys1997cut}. Of particular interest is the $\mathcal{O}(\sqrt{\log(d)})$-approximation proposed in \cite{arora2009expander}, that relies on embedding the nodes of the graph onto the $d$-dimensional space. The embedding can be obtained by solving the SDP  \cite{vladarean2020conditional} 
\begin{align}\label{sparsest cut}
&\min_{\X \in \cc}\frac{1}{d^2}\sum_{i,j}[\mathbf{L}]_{ij}[\X]_{ij}\\
&\text{s.t.}\; d\; \text{Tr}(\X)-\text{Tr}(\textbf{1}^{d\times d}\X)=\frac{d^2}{2},\label{cons_p31}\\
&\quad\; \X^i_j+\X^j_k-\X^i_k-\X^j_j\leq 0,\; \forall \; i,j,k\in \mathcal{V},\label{cons_p32}
\end{align}
where $\mathbf{L}$ is the Laplacian of $\mathcal{G}$ and $\cc := \{\X \succeq 0, \text{tr}(\X)\leq d \}$. The problem can be formulated as $\pd$  and solved using \md by considering only a subset of the summands in the objective function at every iteration. Further, if only a subset of the constraints are considered at every iteration, the problem is cast as $\ps$ and can be solved using \ms. 

We will adopt the experimental setup of \cite{vladarean2020conditional} and  use \textit{mammalia-primate-association-13} graph from the Network Repository dataset  \cite{nr} which has 25 nodes and 181 edges. The SDP dimension for this graph is $\X\in \mathbb{R}^{25\times 25}$ while the number of constraints is around $6.9\times10^3.$  We set the batch size to $5\%$ for all algorithms and additionally use only $5\%$ of the randomly sampled constraints at each iteration for H1-SFW and \ms. Finally, we set $\mu_c=1.5\;\text{and}\;1$ for \md and \ms, respectively. For SHCGM and H1-SFW, we use the same parameters as in \cite{locatello2019stochastic} and \cite{vladarean2020conditional}, respectively.  

We also test the performance of the ZO methods by assuming that only the loss function values are available but not its gradients. We set $\rho_k=\frac{2}{\sqrt{m}(k+1)}$ with $m=d^2=625$, and the approximate gradients are obtained using \eqref{coo:grad_apprx}.

We run all the algorithms for $10^5$ iterations and plot the convergence results in Fig.\ref{scut_obj} and Fig.\ref{scut_cv1}. Both the proposed algorithms outperform SHCGM and H1-SFW. Also, the ZO variants of the proposed algorithms show almost similar performance. However, these require $2m$ SZO calls per iteration. 
Next, we plot the performance variation with trimming threshold $\tau$ for the trimmed variant of both the proposed algorithms. The center plots (Fig.\ref{LO1_3} and Fig.\ref{LO2_3}) show the variation in LMO ($\%$) while the right plots (Fig.\ref{OT1_3} and Fig.\ref{OT2_3}) show the effect of trimming on the objective for few  threshold values. Observe that, for both the cases, the trimmed version provides a significant reduction in total LMO calls. In particular, at $\tau_0=1$, T-\md skips around $40\%$ of LMO calls while at $\tau_0=5$, T-\ms skips around $37\%$ of LMO calls while maintaining a similar performance compared to their non-trimmed version.

\section{Conclusion} \label{sec:conclusion}
This work puts forth projection-free algorithms \md and \ms to solve constrained stochastic optimization algorithms. The problems considered here contain two sets of constraints. The first set of constraints are deterministic and difficult to project onto, motivated from the semidefinite cone constraints that arise in semidefinite programming. The second set of constraints are easy to project onto, but may be large in number or stochastic, and can model the very large number of constraints that arise in context of semidefinite relaxation. We use the stochastic Frank-Wolfe (FW) method to develop projection-free algorithms for solving these problems. The second set of constraints is incorporated using an indicator function within the objective and Nesterov's smoothing is applied to simplify the application of the FW method. Different from existing FW methods to solve similar problems, we utilize a momentum-based gradient tracker that results in improved convergence rates, at par with the set-constrained problems. We also develop zeroth-order algorithms for solving the same problems, yielding gradient-free and projection-free algorithms for solving the same problems while again achieving state-of-the-art convergence rates.  We have also proposed a variant of both the algorithms employing a novel trimming strategy that reduces the number of times linear minimization problem  required to be solved over the  domain to reach a certain accuracy.   
Finally, the performance of the different algorithms is validated numerically on the problems of sparse matrix estimation, clustering, and sparsest cut. Also, to highlight   the usefulness of the trimming scheme, in experiments, we plot the reduction in LMO calls at different threshold and studied the effect of trimming (LMO calls) on the optimality gap. Results showed that the trimmed versions provide significant reduction in total LMO calls while maintaining a similar convergence performance as the original non-trimmed methods.
\appendices
\section{ Proof of Lemma \ref{track2}}\label{Proof:track2}
To begin with, using the definition of the projection operation as well as Assumptions \ref{compact} and \ref{spectral norm}, we state the following result that will be used repeatedly:
\begin{align}
&\Ek{\norm{\G(\xi_k)^\mathsf{T}(\G(\xi_k)\x_k-\Pi_{\cX_{\xi_k}}(\G(\xi_k)\x_k))}^2} \nonumber\\
&\leq \Ek{\norm{\G(\xi_k)}^2}\mbE_k\norm{\G(\xi_k)\x_k-\Pi_{\cX_{\xi_k}}(\G(\xi_k)\x_k)}^2 \\
&\leq L_G^2 \Ek{\norm{\G(\xi_k)\x_k-\G(\xi_k)\x^{\star}}^2} \leq L_G^2D^2 \label{gbound}
\end{align}
(a)	Starting with definition of $\y_k$ from \eqref{estimator2}, observing that $\Ex{\check{g}_k(\x) - \check{g}_k(\x,\xi)} = \Ex{\check{g}_{k-1}(\x) - \check{g}_{k-1}(\x,\xi)} = 0$, and proceeding similarly as in \eqref{ttem1z}-\eqref{yk_1stz}, we can write
\begin{align}
\Ek{\norm{\y_k - \check{g}_k(\x_k)}^2}& = (1-\gamma_k)^2 \norm{\y_{k-1} - \check{g}_{k-1}(\x_{k-1})}^2 \nonumber\\
& + \mbE_k\Big[\|(\check{g}_k(\x_k,\xi_k) - \check{g}_k(\x_k)) \left. + (1-\gamma_k)\big(\check{g}_{k-1}(\x_{k-1}) - \check{g}_{k-1}(\x_{k-1},\xi_k)\|^2 \right]\label{yk_2}
\end{align}
Defining 
\begin{align}
\ssX_k &:= \gamma_k(\check{g}_k(\x_k,\xi_k) - \check{g}_k(\x_k))\\
\ssY_k &:= (1-\gamma_k)(\check{g}_k(\x_k,\xi_k) - \check{g}_{k-1}(\x_{k-1},\xi_k)),\label{def_Y}
\end{align}
and using the inequality $\En{\ssX_k + \ssY_k - \mbE_k\ssY_k} \leq 2\En{\ssX_k} + 2\En{\ssY_k}$, we obtain
\begin{align}
\En{\y_k - \check{g}_k(\x_k)}& \leq (1-\gamma_k)^2 \norm{\y_{k-1} - \check{g}_{k-1}(\x_{k-1})}^2 +2\En{\ssX_k} + 2\En{\ssY_k}\label{yk_3}.
\end{align}
Let $H_{\mu}(\G\x):=\mbE[h_{\mu}(\G(\xi)\x,\cX_\xi)]$, so that $\n_{\x} H_{\mu}(\G\x)=\mbE[\n_{\x} h_{\mu}(\G(\xi)\x,\cX_\xi)]$ and $\En{\ssX_k}$ can be bounded as
\begin{align}\label{rhs2ykb}
&\En{\ssX_k}= \gamma_k^2\En{\check{g}_k(\x_k,\xi_k) - \check{g}_k(\x_k)}\nonumber\\&\leq 2\En{g(\x_k)-g(\x_k,\xi_k)}+2\Ek{\norm{\n_{\x} h_{\mu_k}(\G(\xi_k)\x_k,\cX_{\xi_k})-\n_{\x} H_{\mu_k}(\G\x_k)}^2}.
\end{align}
We can further bound the last term of \eqref{rhs2ykb} using the non-negativity of the variance and from \eqref{grad_Hh} as
\begin{align}
&\Ek{\norm{\n_{\x} h_{\mu_k}(\G(\xi_k)\x_k,\cX_{\xi_k})-\n_{\x} H_{\mu_k}(\G\x_k)}^2}\nonumber\\
& \leq \Ek{\norm{\n_{\x} h_{\mu_k}(\G(\xi_k)\x_k,\cX_{\xi_k})}^2}\nonumber\\
& = \frac{1}{\mu^2_k}\Ek{\norm{\G(\xi_k)^\mathsf{T}(\G(\xi_k)\x_k-\Pi_{\cX_{\xi_k}}(\G(\xi_k)\x_k))}^2}\label{tempi1}\\
& \leq \frac{L_G^2D^2}{\mu^2_k},\label{var_g}
\end{align}	
where \eqref{var_g} follows from \eqref{gbound}. Substituting, we obtain
\begin{align}\label{X2}
\En{\ssX_k}&\!\leq\! 2\gamma_k^2\mbE_k[\norm{g(\x_k)\!-\!g(\x_k,\xi_k)}^2]\!+\!\frac{2\gamma_k^2L_G^2D^2}{\mu^2_k}.
\end{align}
Now, taking norm square of $\ssY_k$ defined in \eqref{def_Y}, dropping the factor $(1-\gamma_k)^2$ and introducing $\n_{\x} h_{\mu_{k-1}}(\G(\xi_k)\x_k,\cX_{\xi_k})$ we obtain
\begin{align}
&\En{\ssY_k} =\Ek{\norm{\check{g}_k(\x_k,\xi_k) - \check{g}_{k-1}(\x_{k-1},\xi_k)}^2} \label{t11}\\
&\leq 3\Ek{\norm{g(\x_k,\xi_k) - g(\x_{k-1},\xi_k)}^2} + 3B_k + 3C_k
\end{align}
where,
\begin{align}
B_k&:= \mbE_k \big[\|\n_{\x} h_{\mu_k}(\G(\xi_k)\x_k,\cX_{\xi_k})-\n_{\x} h_{\mu_{k-1}}(\G(\xi_{k})\x_{k},\cX_{\xi_k})\|^2\big] \\
C_k&:=\mbE_k \big[\|\n_{\x} h_{\mu_{k-1}}(\G(\xi_k)\x_k,\cX_{\xi_k}) -\n_{\x} h_{\mu_{k-1}}(\G(\xi_{k})\x_{k-1},\cX_{\xi_k})\|^2\big]
\end{align}
Here, $B_k$ can be bounded using \eqref{gbound} as
\begin{align}\label{rht1}
B_k &= \left(\tfrac{1}{\mu_k}-\tfrac{1}{\mu_{k-1}}\right)^2\times \En{\G^\T(\xi_k)(\G(\xi_k)\x_k-\Pi_{\cX_{\xi_k}}(\G(\xi_k)\x_k))}\nonumber\\
&\leq \left(\tfrac{1}{\mu_k}-\tfrac{1}{\mu_{k-1}}\right)^2 L_G^2D^2.
\end{align}
Likewise, $C_k$ can be bounded using \eqref{grad_Hh} as:
\begin{align}
C_k&=\frac{1}{\mu_{k-1}^2}\mbE_k\bigg[\|\G^\T(\xi_k)\G(\xi_k)(\x_{k}-\x_{k-1})+ \G^\T(\xi_k)\left(\Pi_{\cX_{\xi_k}}(\G(\xi_k)\x_{k-1}) -\Pi_{\cX_{\xi_k}}(\G(\xi_k)\x_{k})\right) \|^2\bigg]\nonumber\\
&\leq\frac{2}{\mu_{k-1}^2}\mbE_k\bigg[\norm{\G^\T(\xi_k)\G(\xi_k)(\x_{k}-\x_{k-1})}^2\label{tt1}+ \norm{\G^\T(\xi_k)\left(\Pi_{\cX_{\xi_k}}(\G(\xi_k)\x_{k-1}) -\Pi_{\cX_{\xi_k}}(\G(\xi_k)\x_{k})\right) }^2\bigg]
\end{align}
Here, the last term can be bounded using the non-expansiveness property of the projection operator as
\begin{align}
&\Ek{\norm{\G^\T(\xi_k)\left(\Pi_{\cX_{\xi_k}}(\G(\xi_k)\x_{k-1}) -\Pi_{\cX_{\xi_k}}(\G(\xi_k)\x_{k})\right) }^2} \nonumber\\
&\leq \mbE\norm{\G(\xi_k)}^2\mbE_k{\norm{\Pi_{\cX_{\xi_k}}(\G(\xi_k)\x_{k-1}) -\Pi_{\cX_{\xi_k}}(\G(\xi_k)\x_{k})}^2} \nonumber\\
&\leq L_G^2\mbE_k{\norm{\G(\xi_k)(\x_{k-1} -\x_{k})}^2}
\end{align}
Substituting and simplifying, we obtain
\begin{align}
C_k \leq \frac{4L_G^2}{\mu_{k-1}^2}\mbE_k\norm{\G(\xi_k)(\x_{k-1} -\x_{k})}^2 &\leq \frac{4\eta^2_{k-1}L_G^2D^2}{\mu_{k-1}^2}
\end{align}
where we have used the update equation and Assumption \ref{compact}. Thus, we have
\begin{align}
\mbE_k\|\ssY_k\|^2 &\leq 3\mbE_k\|g(\x_k,\xi_k) - g(\x_{k-1},\xi_k)\|^2+ 3L_G^2D^2 \left(\left(\tfrac{1}{\mu_k}-\tfrac{1}{\mu_{k-1}}\right)^2 + \tfrac{4\eta_{k-1}^2}{\mu_{k-1}^2}\right).
\end{align}
From \eqref{xk1} and \eqref{yk1} for SFO oracle, we have 
\begin{align}
&\gamma_k^2\Ex{\norm{ g(\x_k)\! -\!  g(\x_k,\xi_k)}^2}\! \leq \gamma_k^2\sigma^2, \label{xk11}\\
&\Ek{\norm{ g(\x_k,\xi_k) -  g(\x_{k-1},\xi_k)}^2}\leq L^2\eta_{k-1}^2D^2\label{yk11}
\end{align}
Likewise, for SZO oracle, from \eqref{xk2} and \eqref{yk2}, we have
\begin{align}
&\gamma_k^2\Ex{\norm{g(\x_k) - g(\x_k,\xi_k)}^2}\leq  3\gamma_k^2\sigma^2 + 6mL^2\rho_k^2,\label{xk22}\\
&\Ek{\norm{g(\x_{k-1},\xi_k)-g(\x_{k},\xi_k)}^2}\leq 3m L^2\rho^2_k + 3mL^2\rho_{k-1}^2	+ 3\eta_{k-1}^2L^2D^2\label{yk22} 
\end{align}
Using bounds from \eqref{xk11}-\eqref{yk22}, the results for both SFO and SZO cases can be unified as
\begin{align}\label{bt1}
\Ek{\norm{\ssX_k + \ssY_k - \Ek{\ssY_k}}^2}
& \leq 12\gamma_k^2\sigma^2 + 60m L^2\rho^2_{k-1}+ 18\eta_{k-1}^2L^2D^2 \\
& + 2L_G^2D^2\left(\tfrac{2\gamma_k^2}{\mu_k^2} + \tfrac{12\eta_{k-1}^2}{\mu_{k-1}^2} + 3\left(\tfrac{1}{\mu_k}-\tfrac{1}{\mu_{k-1}}\right)^2 \right)
\end{align}
where recall our convention that $\rho_k = 0$ for SFO case.  Substituting  \eqref{bt1} into \eqref{yk_3}, we obtain the desired expression.
%

\noindent(b) Setting  $\gamma_k=\frac{1}{k},\; \mu_k =\frac{\mu_c}{(k+1)^{1/4}},$ $\rho_k\leq\frac{D}{\sqrt{m}(k+1)},$ and $ \eta_k=\frac{2}{k+1}$ in \eqref{mainlem2}, we obtain
\begin{align}
\Ek{\norm{\y_k - \check{g}_k(\x_k)}^2} &\leq \left(1-\frac{1}{k}\right)^2 \norm{\y_{k-1} - \check{g}_{k-1}(\x_{k-1})}^2 \nonumber\\
&+\frac{12(\sigma^2+11L^2D^2)}{k^2} +\frac{(4\sqrt{2}+\tfrac{3}{8}+96)L_G^2D^2}{\mu_c^2k^{3/2}}
\label{mainlem3}
\end{align}
where we have used the fact that $\frac{(k+1)^{1/2}}{k^2}\leq \frac{\sqrt{2}}{k^{3/2}}$ and  $\tfrac{1}{\mu_k}-\tfrac{1}{\mu_{k-1}}=\tfrac{(k+1)^{1/4}-k^{1/4}}{\mu_{c}}\leq \tfrac{1/4}{\mu_c k^{(3/4)}}$.
Now, using the fact that $k^{3/2}\leq k^2$ and using Lemma \ref{simplelemma} with $a=3/2$ and $A=12(\sigma^2+11L^2D^2)+103L_G^2D^2\mu_c^{-2}$, we obtain
\begin{align}\label{bound22}
&\mbE_k\norm{\y_k - \check{g}_k(\x_k)}^2\!\leq\! \frac{48(\sigma^2\!+\!11L^2D^2\!+\!9L_G^2\mu_c^{-2}D^2)}{k^{1/2}},
\end{align}

For SFO case, we get the desired bound simply setting $\check{g}_k(\x_k)=\ffb $ in \eqref{bound22}. For SZO case, we start with $\Ex{\norm{\y_k - \ffb }^2}$, introduce $\n \tilde{F}_{\mu_k}(\x_k):= \mbE[\tilde{\n} f(\x_k,\xi) + \n_{\x} h_{\mu_k}(\G(\xi)\x_k,\cX_{\xi})]$, and obtain the bound by setting $\check{g}_k(\x_k)=\nt f(\x_k)$ in \eqref{bound22}, to yield
\begin{align}
\mbE\|\ffb  -\y_k\|^2&\leq 2\mbE\|\ffb -\ffgbz \|^2+2\mbE\|\ffgbz -\y_k\|^2
\nonumber\\&=2\mbE\|\n f(\x_k) -\nt f(\x_k)\|^2 + 2\mbE\|\ffgbz -\y_k\|^2
\nonumber\\&\leq 2\left[m\rho_k^2 L^2+\frac{48(\sigma^2+11L^2D^2+9L_G^2\mu_c^{-2}D^2)}{k^{1/2}}\right]\nonumber\\&\leq\frac{96(\sigma^2+12L^2D^2+9L_G^2\mu_c^{-2}D^2)}{k^{1/2}},\label{aa}
\end{align}
where we have used \eqref{coro_ZOGrad} and the inequality $k^{1/2}\leq k^2$. Since, the bounds for the SFO \eqref{bound22} and SZO \eqref{aa} cases differ only in constant factors, we will use \eqref{aa} as a unified bound for both the cases.

\bibliographystyle{IEEEtran}
\bibliography{ref}

\begin{thebibliography}{10}
\providecommand{\url}[1]{#1}
\csname url@samestyle\endcsname
\providecommand{\newblock}{\relax}
\providecommand{\bibinfo}[2]{#2}
\providecommand{\BIBentrySTDinterwordspacing}{\spaceskip=0pt\relax}
\providecommand{\BIBentryALTinterwordstretchfactor}{4}
\providecommand{\BIBentryALTinterwordspacing}{\spaceskip=\fontdimen2\font plus
\BIBentryALTinterwordstretchfactor\fontdimen3\font minus
  \fontdimen4\font\relax}
\providecommand{\BIBforeignlanguage}[2]{{%
\expandafter\ifx\csname l@#1\endcsname\relax
\typeout{** WARNING: IEEEtran.bst: No hyphenation pattern has been}%
\typeout{** loaded for the language `#1'. Using the pattern for}%
\typeout{** the default language instead.}%
\else
\language=\csname l@#1\endcsname
\fi
#2}}
\providecommand{\BIBdecl}{\relax}
\BIBdecl

\bibitem{shapiro2014lectures}
A.~Shapiro, D.~Dentcheva, and A.~Ruszczy{\'n}ski, \emph{Lectures on stochastic
  programming: modeling and theory}.\hskip 1em plus 0.5em minus 0.4em\relax
  SIAM, 2014.

\bibitem{ahmed2013blind}
A.~Ahmed, B.~Recht, and J.~Romberg, ``Blind deconvolution using convex
  programming,'' \emph{IEEE Transactions on Information Theory}, vol.~60,
  no.~3, pp. 1711--1732, 2013.

\bibitem{conn2009introduction}
A.~R. Conn, K.~Scheinberg, and L.~N. Vicente, \emph{Introduction to
  derivative-free optimization}.\hskip 1em plus 0.5em minus 0.4em\relax SIAM,
  2009.

\bibitem{chen2017zoo}
P.-Y. Chen, H.~Zhang, Y.~Sharma, J.~Yi, and C.-J. Hsieh, ``Zoo: Zeroth order
  optimization based black-box attacks to deep neural networks without training
  substitute models,'' in \emph{Proceedings of the 10th ACM Workshop on
  Artificial Intelligence and Security}, 2017, pp. 15--26.

\bibitem{nemirovski2009robust}
A.~Nemirovski, A.~Juditsky, G.~Lan, and A.~Shapiro, ``Robust stochastic
  approximation approach to stochastic programming,'' \emph{SIAM Journal on
  optimization}, vol.~19, no.~4, pp. 1574--1609, 2009.

\bibitem{tran2018smooth}
Q.~Tran-Dinh, O.~Fercoq, and V.~Cevher, ``A smooth primal-dual optimization
  framework for nonsmooth composite convex minimization,'' \emph{SIAM Journal
  on Optimization}, vol.~28, no.~1, pp. 96--134, 2018.

\bibitem{yurtsever2019conditional2a}
A.~Yurtsever, O.~Fercoq, and V.~Cevher, ``A conditional-gradient-based
  augmented lagrangian framework,'' in \emph{International Conference on
  Machine Learning}.\hskip 1em plus 0.5em minus 0.4em\relax PMLR, 2019, pp.
  7272--7281.

\bibitem{fercoq2019almost}
O.~Fercoq, A.~Alacaoglu, I.~Necoara, and V.~Cevher, ``Almost surely constrained
  convex optimization,'' \emph{arXiv preprint arXiv:1902.00126}, 2019.

\bibitem{garrigues2008homotopy}
P.~Garrigues and L.~Ghaoui, ``An homotopy algorithm for the lasso with online
  observations,'' \emph{Advances in neural information processing systems},
  vol.~21, pp. 489--496, 2008.

\bibitem{peng2007approximating}
J.~Peng and Y.~Wei, ``Approximating k-means-type clustering via semidefinite
  programming,'' \emph{SIAM journal on optimization}, vol.~18, no.~1, pp.
  186--205, 2007.

\bibitem{arora2009expander}
S.~Arora, S.~Rao, and U.~Vazirani, ``Expander flows, geometric embeddings and
  graph partitioning,'' \emph{Journal of the ACM (JACM)}, vol.~56, no.~2, pp.
  1--37, 2009.

\bibitem{huang2014scalable}
Q.~Huang, Y.~Chen, and L.~Guibas, ``Scalable semidefinite relaxation for
  maximum a posterior estimation,'' in \emph{International Conference on
  Machine Learning}, 2014, pp. 64--72.

\bibitem{zhao1998semidefinite}
Q.~Zhao, S.~E. Karisch, F.~Rendl, and H.~Wolkowicz, ``Semidefinite programming
  relaxations for the quadratic assignment problem,'' \emph{Journal of
  Combinatorial Optimization}, vol.~2, no.~1, pp. 71--109, 1998.

\bibitem{goemans1995improved}
M.~X. Goemans and D.~P. Williamson, ``Improved approximation algorithms for
  maximum cut and satisfiability problems using semidefinite programming,''
  \emph{Journal of the ACM (JACM)}, vol.~42, no.~6, pp. 1115--1145, 1995.

\bibitem{raghunathan2018semidefinite}
A.~Raghunathan, J.~Steinhardt, and P.~S. Liang, ``Semidefinite relaxations for
  certifying robustness to adversarial examples,'' \emph{Advances in Neural
  Information Processing Systems}, vol.~31, 2018.

\bibitem{kulis2007fast}
B.~Kulis, A.~C. Surendran, and J.~C. Platt, ``Fast low-rank semidefinite
  programming for embedding and clustering,'' in \emph{Artificial Intelligence
  and Statistics}.\hskip 1em plus 0.5em minus 0.4em\relax PMLR, 2007, pp.
  235--242.

\bibitem{latorre2020lipschitz}
F.~Latorre, P.~Rolland, and V.~Cevher, ``Lipschitz constant estimation of
  neural networks via sparse polynomial optimization,'' \emph{arXiv preprint
  arXiv:2004.08688}, 2020.

\bibitem{locatello2019stochastic}
F.~Locatello, A.~Yurtsever, O.~Fercoq, and V.~Cevher, ``Stochastic
  {Frank-Wolfe} for composite convex minimization,'' in \emph{Advances in
  Neural Information Processing Systems}, 2019, pp. 14\,269--14\,279.

\bibitem{vladarean2020conditional}
M.-L. Vladarean, A.~Alacaoglu, Y.-P. Hsieh, and V.~Cevher, ``Conditional
  gradient methods for stochastically constrained convex minimization,''
  \emph{arXiv preprint arXiv:2007.03795}, 2020.

\bibitem{jaggi2013revisiting}
M.~Jaggi, ``Revisiting {Frank-Wolfe}: Projection-free sparse convex
  optimization,'' in \emph{Proceedings of the 30th international conference on
  machine learning}, no. CONF, 2013, pp. 427--435.

\bibitem{nesterov2005smooth}
Y.~Nesterov, ``Smooth minimization of non-smooth functions,''
  \emph{Mathematical programming}, vol. 103, no.~1, pp. 127--152, 2005.

\bibitem{zhang2020one}
M.~Zhang, Z.~Shen, A.~Mokhtari, H.~Hassani, and A.~Karbasi, ``One sample
  stochastic frank-wolfe,'' in \emph{International Conference on Artificial
  Intelligence and Statistics}.\hskip 1em plus 0.5em minus 0.4em\relax PMLR,
  2020, pp. 4012--4023.

\bibitem{xie2020efficient}
J.~Xie, Z.~Shen, C.~Zhang, B.~Wang, and H.~Qian, ``Efficient projection-free
  online methods with stochastic recursive gradient.'' in \emph{AAAI}, 2020,
  pp. 6446--6453.

\bibitem{braun2017lazifying}
G.~Braun, S.~Pokutta, and D.~Zink, ``Lazifying conditional gradient
  algorithms,'' in \emph{International conference on machine learning}.\hskip
  1em plus 0.5em minus 0.4em\relax PMLR, 2017, pp. 566--575.

\bibitem{kerdreux2018frank}
T.~Kerdreux, F.~Pedregosa, and A.~d’Aspremont, ``Frank-wolfe with subsampling
  oracle,'' in \emph{International Conference on Machine Learning}.\hskip 1em
  plus 0.5em minus 0.4em\relax PMLR, 2018, pp. 2591--2600.

\bibitem{mhammedi2021efficient}
Z.~Mhammedi, ``Efficient projection-free online convex optimization with
  membership oracle,'' \emph{arXiv preprint arXiv:2111.05818}, 2021.

\bibitem{hazan2012projection}
E.~Hazan and S.~Kale, ``Projection-free online learning,'' \emph{arXiv preprint
  arXiv:1206.4657}, 2012.

\bibitem{hazan2016variance}
E.~Hazan and H.~Luo, ``Variance-reduced and projection-free stochastic
  optimization,'' in \emph{International Conference on Machine Learning}.\hskip
  1em plus 0.5em minus 0.4em\relax PMLR, 2016, pp. 1263--1271.

\bibitem{shen2019complexities}
Z.~Shen, C.~Fang, P.~Zhao, J.~Huang, and H.~Qian, ``Complexities in
  projection-free stochastic non-convex minimization,'' in \emph{The 22nd
  International Conference on Artificial Intelligence and Statistics}.\hskip
  1em plus 0.5em minus 0.4em\relax PMLR, 2019, pp. 2868--2876.

\bibitem{defazio2018ineffectiveness}
A.~Defazio and L.~Bottou, ``On the ineffectiveness of variance reduced
  optimization for deep learning,'' \emph{arXiv preprint arXiv:1812.04529},
  2018.

\bibitem{akhtar2021momentum}
Z.~Akhtar and K.~Rajawat, ``Momentum based projection free stochastic
  optimization under affine constraints,'' in \emph{2021 American Control
  Conference (ACC)}.\hskip 1em plus 0.5em minus 0.4em\relax IEEE, 2021, pp.
  2619--2624.

\bibitem{mokhtari2020stochastic}
A.~Mokhtari, H.~Hassani, and A.~Karbasi, ``Stochastic conditional gradient
  methods: From convex minimization to submodular maximization,'' \emph{Journal
  of Machine Learning Research}, vol.~21, no. 105, pp. 1--49, 2020.

\bibitem{akhtar2021conservative}
Z.~Akhtar, A.~S. Bedi, and K.~Rajawat, ``Conservative stochastic optimization
  with expectation constraints,'' \emph{IEEE Transactions on Signal
  Processing}, vol.~69, pp. 3190--3205, 2021.

\bibitem{lan2016conditional}
G.~Lan and Y.~Zhou, ``Conditional gradient sliding for convex optimization,''
  \emph{SIAM Journal on Optimization}, vol.~26, no.~2, pp. 1379--1409, 2016.

\bibitem{lan2017conditional}
G.~Lan, S.~Pokutta, Y.~Zhou, and D.~Zink, ``Conditional accelerated lazy
  stochastic gradient descent,'' \emph{arXiv preprint arXiv:1703.05840}, 2017.

\bibitem{lu2020generalized}
H.~Lu and R.~M. Freund, ``Generalized stochastic {Frank--Wolfe} algorithm with
  stochastic “substitute” gradient for structured convex optimization,''
  \emph{Mathematical Programming}, pp. 1--33, 2020.

\bibitem{yurtsever2018conditional}
A.~Yurtsever, O.~Fercoq, F.~Locatello, and V.~Cevher, ``A conditional gradient
  framework for composite convex minimization with applications to semidefinite
  programming,'' \emph{arXiv preprint arXiv:1804.08544}, 2018.

\bibitem{sahu2019towards}
A.~K. Sahu, M.~Zaheer, and S.~Kar, ``Towards gradient free and projection free
  stochastic optimization,'' in \emph{The 22nd International Conference on
  Artificial Intelligence and Statistics}.\hskip 1em plus 0.5em minus
  0.4em\relax PMLR, 2019, pp. 3468--3477.

\bibitem{patrascu2017nonasymptotic}
A.~Patrascu and I.~Necoara, ``Nonasymptotic convergence of stochastic proximal
  point methods for constrained convex optimization,'' \emph{The Journal of
  Machine Learning Research}, vol.~18, no.~1, pp. 7204--7245, 2017.

\bibitem{wang2015random}
M.~Wang, Y.~Chen, J.~Liu, and Y.~Gu, ``Random multi-constraint projection:
  Stochastic gradient methods for convex optimization with many constraints,''
  \emph{arXiv preprint arXiv:1511.03760}, 2015.

\bibitem{reddi2016stochastic}
S.~J. Reddi, S.~Sra, B.~P{\'o}czos, and A.~Smola, ``Stochastic {Frank-Wolfe}
  methods for nonconvex optimization,'' in \emph{2016 54th Annual Allerton
  Conference on Communication, Control, and Computing (Allerton)}.\hskip 1em
  plus 0.5em minus 0.4em\relax IEEE, 2016, pp. 1244--1251.

\bibitem{yurtsever2019conditional}
A.~Yurtsever, S.~Sra, and V.~Cevher, ``Conditional gradient methods via
  stochastic path-integrated differential estimator,'' in \emph{International
  Conference on Machine Learning}.\hskip 1em plus 0.5em minus 0.4em\relax PMLR,
  2019, pp. 7282--7291.

\bibitem{li2021momentum}
B.~Li, M.~Couti{\~n}o, G.~B. Giannakis, and G.~Leus, ``A momentum-guided
  frank-wolfe algorithm,'' \emph{IEEE Transactions on Signal Processing},
  vol.~69, pp. 3597--3611, 2021.

\bibitem{balasubramanian2018zeroth}
K.~Balasubramanian and S.~Ghadimi, ``Zeroth-order (non)-convex stochastic
  optimization via conditional gradient and gradient updates,'' in
  \emph{Advances in Neural Information Processing Systems}, 2018, pp.
  3455--3464.

\bibitem{huang2020accelerated}
F.~Huang, L.~Tao, and S.~Chen, ``Accelerated stochastic gradient-free and
  projection-free methods,'' in \emph{International Conference on Machine
  Learning}.\hskip 1em plus 0.5em minus 0.4em\relax PMLR, 2020, pp. 4519--4530.

\bibitem{nesterov2017random}
Y.~Nesterov and V.~Spokoiny, ``Random gradient-free minimization of convex
  functions,'' \emph{Foundations of Computational Mathematics}, vol.~17, no.~2,
  pp. 527--566, 2017.

\bibitem{frandi2014complexity}
E.~Frandi, R.~{\~N}anculef, and J.~Suykens, ``Complexity issues and
  randomization strategies in frank-wolfe algorithms for machine learning,''
  \emph{arXiv preprint arXiv:1410.4062}, 2014.

\bibitem{li2021communication}
W.~Li, Z.~Wu, T.~Chen, L.~Li, and Q.~Ling, ``Communication-censored distributed
  stochastic gradient descent,'' \emph{IEEE Transactions on Neural Networks and
  Learning Systems}, 2021.

\bibitem{liu2018zeroth}
S.~Liu, B.~Kailkhura, P.-Y. Chen, P.~Ting, S.~Chang, and L.~Amini,
  ``Zeroth-order stochastic variance reduction for nonconvex optimization,''
  \emph{Advances in Neural Information Processing Systems}, vol.~31, pp.
  3727--3737, 2018.

\bibitem{ji2019improved}
K.~Ji, Z.~Wang, Y.~Zhou, and Y.~Liang, ``Improved zeroth-order variance reduced
  algorithms and analysis for nonconvex optimization,'' \emph{arXiv preprint
  arXiv:1910.12166}, 2019.

\bibitem{beck2017first}
A.~Beck, \emph{First-order methods in optimization}.\hskip 1em plus 0.5em minus
  0.4em\relax SIAM, 2017.

\bibitem{cutkosky2019momentum}
A.~Cutkosky and F.~Orabona, ``Momentum-based variance reduction in non-convex
  sgd,'' in \emph{Advances in Neural Information Processing Systems}, 2019, pp.
  15\,236--15\,245.

\bibitem{kerdreux2021projection}
T.~Kerdreux, A.~d’Aspremont, and S.~Pokutta, ``Projection-free optimization
  on uniformly convex sets,'' in \emph{International Conference on Artificial
  Intelligence and Statistics}.\hskip 1em plus 0.5em minus 0.4em\relax PMLR,
  2021, pp. 19--27.

\bibitem{garber2015faster}
D.~Garber and E.~Hazan, ``Faster rates for the frank-wolfe method over
  strongly-convex sets,'' in \emph{International Conference on Machine
  Learning}.\hskip 1em plus 0.5em minus 0.4em\relax PMLR, 2015, pp. 541--549.

\bibitem{bauschke2011convex}
H.~H. Bauschke, P.~L. Combettes \emph{et~al.}, \emph{Convex analysis and
  monotone operator theory in Hilbert spaces}.\hskip 1em plus 0.5em minus
  0.4em\relax Springer, 2011, vol. 408.

\bibitem{richard2012estimation}
E.~Richard, P.-A. Savalle, and N.~Vayatis, ``Estimation of simultaneously
  sparse and low rank matrices,'' \emph{arXiv preprint arXiv:1206.6474}, 2012.

\bibitem{zhao2014robust}
Q.~Zhao, D.~Meng, Z.~Xu, W.~Zuo, and L.~Zhang, ``Robust principal component
  analysis with complex noise,'' in \emph{International conference on machine
  learning}, 2014, pp. 55--63.

\bibitem{deshmukh2020improved}
S.~Deshmukh and A.~Dubey, ``Improved covariance matrix estimation with an
  application in portfolio optimization,'' \emph{IEEE Signal Processing
  Letters}, 2020.

\bibitem{mixon2016clustering}
D.~G. Mixon, S.~Villar, and R.~Ward, ``Clustering subgaussian mixtures by
  semidefinite programming,'' \emph{arXiv preprint arXiv:1602.06612}, 2016.

\bibitem{lecun2010mnist}
Y.~LeCun, C.~Cortes, and C.~Burges, ``Mnist handwritten digit database. 2010,''
  \emph{URL http://yann. lecun. com/exdb/mnist}, vol.~7, p.~23, 2010.

\bibitem{chatziafratis2018hierarchical}
V.~Chatziafratis, R.~Niazadeh, and M.~Charikar, ``Hierarchical clustering with
  structural constraints,'' \emph{arXiv preprint arXiv:1805.09476}, 2018.

\bibitem{bonsma2012complexity}
P.~Bonsma, H.~Broersma, V.~Patel, and A.~Pyatkin, ``The complexity of finding
  uniform sparsest cuts in various graph classes,'' \emph{Journal of discrete
  algorithms}, vol.~14, pp. 136--149, 2012.

\bibitem{leighton1999multicommodity}
T.~Leighton and S.~Rao, ``Multicommodity max-flow min-cut theorems and their
  use in designing approximation algorithms,'' \emph{Journal of the ACM
  (JACM)}, vol.~46, no.~6, pp. 787--832, 1999.

\bibitem{shmoys1997cut}
D.~B. Shmoys, ``Cut problems and their application to divide-and-conquer,''
  \emph{Approximation algorithms for NP-hard problems}, pp. 192--235, 1997.

\bibitem{nr}
\BIBentryALTinterwordspacing
R.~A. Rossi and N.~K. Ahmed, ``The network data repository with interactive
  graph analytics and visualization,'' in \emph{AAAI}, 2015. [Online].
  Available: \url{http://networkrepository.com}
\BIBentrySTDinterwordspacing

\end{thebibliography}

\newpage
\section{Supplementary Material}
Before proceeding we replicate some technical results as following corollaries that are borrowed  from related works and  will be used in proving our theorems. We also state a preliminary lemma \ref{simplelemma} that will apply to most of the proofs.
\begin{corollary} \label{coro:borrow1}
	(\textit{Lemma 10 in \cite{tran2018smooth}}): Suppose $h_{\mu}(\x)$ be the smooth approximation  of a non-smooth convex function $h$ and $\lam^{\star}_{\mu}(\cdot)$ be the solution of the dual problem, then following inequalities holds true
	\begin{align}\label{borrow1}
	h(\x_1)\geq h_{\mu}(\x_2)+ \ip{\n h_{\mu}(\x_2),\x_1-\x_2}+\frac{\mu}{2}\|\lam^{\star}_{\mu}(\x_2)\|^2,
	\end{align}
	\begin{align}\label{borrow2}
	h_{\mu_1}(\x)\leq h_{\mu_2}(\x)+\frac{\mu_2-\mu_1}{2}\|\lam^{\star}_{\mu_1}(\x)\|^2.
	\end{align}
\end{corollary}
\begin{corollary} \label{coro:borrow3} (\textit{Theorem 2(a) in \cite{locatello2019stochastic}}): Let $(\x^{\star},\lam^{\star})$	be a saddle point of $\mathcal{L}(\x,\r,\lam^{\star})=f(\x)+\ip{\lam^{\star},\G\x-\r}$, then from Lagrange saddle point theory  following bound holds $\forall\; \x \in \cc$ and $\r \in \cX$
	\begin{align}\label{borrow3}
	\mbE[f(\x_k)]-f(\x^{\star})&\geq-\mbE[\min_{\r\in \cX}\norm{\lam^{\star}}\|\G\x-\r\|]\nonumber\\&= -\norm{\lam^{\star}}\mbE[\cD_{\cX}(\G\x_k)].
	\end{align}
\end{corollary}
\begin{corollary} \label{coro:borrow4}(Lemma B.1 in \cite{vladarean2020conditional}): Consider problem $\ps$ and its Lagrangian formulation as $$\mathcal{L}(\x,\lam):=f(\x)+\int\ip{\G(\xi)\x,\lam(\xi)}-\sup_{\y\in \cX(\xi)}\ip{\y,\lam(\xi)}\alpha(d\xi),$$
	here $\alpha$ denotes the probability measure of the random variable $\xi$. Let $(\x^{\star},\lam^{\star})$ be a saddle point of $\mathcal{L}$, then
	\begin{align}\label{borrow4}
	\int\!\!\!\cD^2_{\cX(\xi)}\!(\G(\xi)\x)\alpha(d\xi)\!\leq\! 4\mu^2\!\norm{\!\lam^{\star}\!}^2\!+\!4\mu[F_{\mu}(\x)\!-\!f(\x^{\star})].
	\end{align}
\end{corollary}

\begin{lemma}\label{simplelemma}
	Let $\{\psi_k\}_{k=1}^t$ be a sequence of numbers satisfying either of the following recursions:
	\begin{align}
	\psi_k &\leq \left(1 - \frac{1}{k}\right)^2\psi_{k-1} + \frac{A}{k^a} \label{simplerec1}\\
	\psi_k &\leq \left(1 - \frac{2}{k+1}\right)\psi_{k-1} + \frac{A}{k^a} \label{simplerec2}
	\end{align}
	for some $A > 0$ and $1 \leq a \leq 2$. Then it holds that $\psi_k \leq 4A/k^{a-1}$ for all $t \geq 1$. 
\end{lemma}
\begin{IEEEproof}[Proof of \eqref{simplerec1}]
	Multiplying both sides of \eqref{simplerec1} by $k^2$, we obtain
	\begin{align}
	k^2 \psi_k &\leq (k-1)^2\psi_{k-1} + Ak^{2-a} 
	\end{align}
	for all $k \geq 1$. Carrying out telescopic sum over $k = 1, \ldots, K$, we obtain
	\begin{align}
	\psi_K &\leq \frac{A}{K^2}\sum_{k=1}^K k^{2-a} \leq A\frac{(K+1)^{3-a}}{(3-a)K^2}  \leq \frac{4A}{K^{a-1}}
	\end{align}
	where we have used the inequality $\sum_{k=1}^{K} k^{n} \leq \frac{(K+1)^{n+1}-1}{(n+1)}$ for $n \geq 0$.
\end{IEEEproof}

\begin{IEEEproof}[Proof of \eqref{simplerec2}]
	Multiplying both sides of \eqref{simplerec2} by $k(k+1)$, we obtain
	\begin{align}
	k(k+1) \psi_k &\leq (k-1)k\psi_{k-1} + Ak^{1-a}(k+1) \\
	& \leq (k-1)k\psi_{k-1} + 2Ak^{2-a}
	\end{align}
	for all $k \geq 1$. Carrying out telescopic sum over $k = 1, \ldots, K$, we obtain
	\begin{align}
	\psi_K \!\leq\! \frac{2A}{K(K+1)}\sum_{k=1}^K k^{2-a}  \!\leq\! 2A\frac{(K+1)^{2-a}}{(3-a)K}\!\leq\! \frac{4A}{K^{a-1}}
	\end{align}
	where we have used the inequality $\sum_{k=1}^K k^{n} \leq \frac{(K+1)^{n+1}-1}{(n+1)}$ for $n \geq 0$.
\end{IEEEproof}
\section{}
\subsection{Proof of Lemma \ref{lemma:basic_F(x)_1st}}\label{proof-smoothgap}
\noindent (a) Since $h_{\mu}$ is $1/\mu$-smooth by construction and $f(\x)$ is $L$-smooth from Assumption \ref{smoothness}, it follows from Assumption \ref{spectral norm} that  $F_{\mu}(\x)$ is $\left(L+\frac{L_G}{\mu}\right)$-smooth. We can therefore write the quadratic upper bound on $F_{\mu_k}$ and use the update in \eqref{update_alg1} to obtain
\begin{align}\label{basic smoothness step}
F_{\mu_k}(\x_{k+1})&\leq F_{\mu_k}(\x_{k})+\eta_k\ip{\nf,\z_k-\x_k}+\frac{\eta_k^2}{2}\left(L+\frac{L_G}{\mu_k}\right)\|\z_k-\x_k\|^2\nonumber\\&\leq F_{\mu_k}(\x_{k})+\eta_k\ip{\nf,\z_k-\x_k}+\frac{\eta_k^2}{2}\left(L+\frac{L_G}{\mu_k}\right)D^2,
\end{align}
where we have used the compactness assumption on $\cc$ (Assumption \ref{compact}). Recall from  \eqref{grad_cons_funct} that
$\nf=\n f(\x_k)+\G^\T\n h_{\mu_k}(\G\x,\cX)$, which allows us to write the second term on the right of \eqref{basic smoothness step} as
\begin{align}
\ip{\nf,\z_k-\x_k}&
=\ip{\n f(\x_k),\z_k-\x_k}+\ip{\G^\T\n h_{\mu_k}(\G\x_k,\cX),\z_k-\x_k}\nonumber\\
&=\ip{\n f(\x_k)-\y_k,\z_k-\x_k}+\ip{\y_k+\G^\T\n h_{\mu_k}(\G\x_k,\cX),\z_k-\x_k}\label{xa1}.
\end{align}
From the definition of $\z_k$ in \eqref{lmo_alg1}, we therefore have that
\begin{align}
\ip{\nf,\z_k-\x_k}\leq \ip{\n f(\x_k)-\y_k,\z_k-\x_k}+\ip{\y_k+\G^\T\n h_{\mu_k}(\G\x_k,\cX),\x^{\star}-\x_k}.
\end{align}
Adding and subtracting $\ip{\n f(\x_k),\x^{\star}-\x_k}$, and using the convexity of $f$, we obtain
\begin{align}
\ip{\nf,\z_k-\x_k}&\leq \ip{\n f(\x_k)-\y_k,\z_k-\x^{\star}}+\ip{\n f(\x_k),\x^{\star}-\x_k}\nonumber\\&\quad+\ip{\G^\T\n h_{\mu_k}(\G\x_k,\cX),\x^{\star}-\x_k}\label{tempa}\\
&\leq\|\n f(\x_k)-\y_k\|\|\z_k-\x^{\star}\|+f(\x^{\star})- f(\x_k)\nonumber\\&\quad+\ip{\G^\T\n h_{\mu_k}(\G\x_k,\cX),\x^{\star}-\x_k}\label{tempi}\\
&\leq\|\n f(\x_k)-\y_k\|D+f(\x^{\star})- f(\x_k)\nonumber\\&\quad+\ip{\G^\T\n h_{\mu_k}(\G\x_k,\cX),\x^{\star}-\x_k},\label{grad1}
\end{align}
where we have used the Cauchy-Schwartz inequality on the first term in \eqref{tempa} and then the compactness assumption on $\cc$. 

Next, we rewrite \eqref{borrow1}  in context of our problem as
\begin{align}
\one_{\cX}(\G\x^{\star}) &\geq h_{\mu_k}(\G\x_k,\cX)+ \ip{\n h_{\mu_k}(\G \x_k,\cX),\G\x^{\star}-\G\x_k}+\frac{\mu_k}{2}\|\lam^{\star}_{\mu_k}(\G\x_k)\|^2
\end{align}
implying that
\begin{align}\label{borrow1_our}
\ip{\G^\T\n h_{\mu_k}(\x_k,\cX),\x^{\star}-\x_k}&=\ip{\n h_{\mu_k}(\x_k,\cX),\G\x^{\star}-\G\x_k}\nonumber\\&\leq 
\one_{\cX}(\G\x^{\star}) - h_{\mu_k}(\G\x_k,\cX)-\frac{\mu_k}{2}\|\lam^{\star}_{\mu_k}(\G\x_k)\|^2.
\end{align}
Using \eqref{borrow1_our} in \eqref{grad1} we therefore obtain
\begin{align}\label{half_grad}
\ip{\nf,\z_k-\x_k}
&\leq \|\n f(\x_k)-\y_k\|D+f(\x^{\star})- f(\x_k)+\one_{\cX}(\G\x^{\star})\nonumber\\&\quad-h_{\mu_k}(\G\x_k,\cX)-\frac{\mu_k}{2}\|\lam^{\star}_{\mu_k}(\G\x_k)\|^2.
\end{align}

Also, rewriting \eqref{borrow2}  in context of our problem we get
\begin{align}
h_{\mu_k}(\G\x_k,\cX)\leq h_{\mu_{k-1}}(\G\x_k,\cX)+\frac{\mu_{k-1}-\mu_k}{2}\|\lam^{\star}_{\mu_k}(\G\x_k)\|^2.\nonumber
\end{align}
Hence, we can write
\begin{align}\label{temp1}
F_{\mu_k}(\x_k) &:= \mbE[f(\x_k,\xi)]+h_{\mu_k}(\G\x_k,\cX)=f(\x_k)+h_{\mu_k}(\G\x_k,\cX)\nonumber\\
&\leq f(\x_k)+h_{\mu_{k-1}}(\G\x_k,\cX)+\frac{\mu_{k-1}-\mu_k}{2}\|\lam^{\star}_{\mu_k}(\G\x_k)\|^2\nonumber\\
&=F_{\mu_{k-1}}(\x_{k})+\frac{\mu_{k-1}-\mu_k}{2}\|\lam^{\star}_{\mu_k}(\G\x_k)\|^2.
\end{align}

Substituting the bounds obtained in \eqref{half_grad} and \eqref{temp1} into \eqref{basic smoothness step} and subtracting $f(\x^{\star})$ from both sides, we obtain
\begin{align}\label{basic smoothness step3}
F_{\mu_k}(\x_{k+1})-f(\x^{\star})
&\leq (1-\eta_k)\left(F_{\mu_{k-1}}(\x_{k})-f(\x^{\star})\right)\nonumber\\&+\frac{\eta_k^2}{2}\left(L+\frac{L_G}{\mu_k}\right)D^2+\eta_k\|\n f(\x_k)-\y_k\|D\nonumber\\&\quad+ \frac{1}{2}\left((1-\eta_k)\mu_{k-1}-\mu_k\right)\|\lam^{\star}_{\mu_k}(\G\x_k)\|^2.
\end{align}
here we used the fact that $\one_{\cX}(\G\x^{\star})=0$. Now, since we have assumed that $\mu_k \geq \mu_{k-1}(1-\eta_k)$, the last term in \eqref{basic smoothness step3} is non-positive and can be dropped, yielding the desired result. 

\noindent (b) Taking expectation in \eqref{smgaprec} and using the inequality $\Ex{\norm{\sX}} \leq \sqrt{\Ex{\norm{\sX}^2}}$, we obtain
\begin{align}
\Ex{F_{\mu_k}(\x_{k+1})-f(\x^{\star})}
&\leq (1-\eta_k)\Ex{F_{\mu_{k-1}}(\x_{k})-f(\x^{\star})}\nonumber\\&+\frac{\eta_k^2}{2}\left(L+\frac{L_G}{\mu_k}\right)D^2+\eta_kD\sqrt{\Ex{\norm{\n f(\x_k)-\y_k}^2}}
\end{align} 
Substituting the result of Lemma \ref{track1} and that of the various parameters, we obtain
\begin{align}
	\Ex{F_{\mu_k}(\x_{k+1})-f(\x^{\star})}&\leq  (1-\frac{2}{k+1})\Ex{F_{\mu_{k-1}}(\x_{k})-f(\x^{\star})} \\
&+\frac{1}{(k+1)^2}\left(L+\frac{L_G\sqrt{k}}{\mu_c}\right)D^2+\frac{8D}{k+1}\sqrt{\frac{3\sigma^2+25L^2D^2}{k}}\nonumber
\end{align}
Note that for this choice of parameters, $\frac{\mu_k}{\mu_{k-1}} = \sqrt{\frac{k-1}{k}}\geq \frac{k-1}{k}=1-\eta_k$, so the condition required for \eqref{mainlem} is satisfied. The bound on the right can therefore be simplified as 
\begin{align}
\mbE[F_{\mu_k}(\x_{k+1})\!-f(\x^{\star})] \leq  (1-\frac{2}{k+1})\Ex{F_{\mu_{k-1}}(\x_{k})\!-f(\x^{\star})} \!+ \! \frac{8\sqrt{3}\sigma D\!+\!\!(41L\!+\!L_G\mu_c^{-1}\!)D^2}{k^{3/2}}
\end{align}
Finally, application of Lemma \ref{simplelemma} yields the required result.

\subsection{Proof of Theorem \ref{theorem1}:}\label{Proof:th1}
\noindent (a) The bound on the optimality gap follows directly from Lemma \ref{lemma:basic_F(x)_1st}:
\begin{align}\label{temp2}
\mbE[f(\x_{k+1})]-f(\x^{\star})&\leq  \mbE[f(\x_{k+1})]-f(\x^{\star})+\frac{1}{2\mu_k}\mbE [\cD^2_{\cX}(\G\x_{k+1})]\nonumber\\& = \mbE[F_{\mu_k}(\x_{k+1})]-f(\x^{\star})\leq \frac{8\sqrt{3}\sigma D+(41L+L_G\mu_c^{-1})D^2}{\sqrt{k}}.
\end{align}
\noindent (b) Since $\x_{k+1}\in \cc$, it follows from  corollary \ref{coro:borrow3}  that
\begin{align}\label{borrow3our}
\mbE[f(\x_{k+1})]-f(\x^{\star})\geq -\norm{\lam^{\star}}\mbE[\cD_{\cX}(\G\x_{k+1})].
\end{align}
Substituting \eqref{borrow3our} into \eqref{temp2} and using the Cauchy-Schwarz inequality, we obtain
\begin{align}\label{temp3}
-\norm{\lam^{\star}}\mbE[\cD_{\cX}(\G\x_{k+1})]+\frac{1}{2\mu_k}(\mbE [\cD_{\cX}(\G\x_{k+1})])^2 &\leq -\norm{\lam^{\star}}\mbE[\cD_{\cX}(\G\x_{k+1})]+\frac{1}{2\mu_k}\mbE [\cD^2_{\cX}(\G\x_{k+1})] \nonumber\\&\leq \frac{8\sqrt{3}\sigma D+(41L+L_G\mu_c^{-1})D^2}{\sqrt{k}} 
\end{align}
which is a quadratic inequality in $\Ex{\cD_\cX(\G\x_{k+1})}$, and can be solved to yield the required result:	
\begin{align}\label{cons_fes1}
\mbE[\cD_{\cX}(\G\x_{k+1})]&\leq \frac{1}{\sqrt{k}}\left(2\mu_c\norm{\lam^{\star}}+\sqrt{2(8\sqrt{3}\sigma D\mu_c+(41L\mu_c+L_G)D^2)}\right) \nonumber\\
&\leq \frac{1}{\sqrt{k}}\left(2\mu_c\norm{\lam^{\star}}+ 6\sqrt{\sigma D \mu_c} + 10D\sqrt{L\mu_c}+2D\sqrt{L_G}\right)\nonumber
\end{align}
\subsection{Proof of Lemma \ref{lemma:F-smooth}}\label{proof:lemma:F-smooth}
Since, $h_{\mu}(\G(\xi)\x,\cX_\xi)$ is $1/\mu$-smooth and $f(\x)$ is $L$-smooth (assumption \ref{smoothness}), it follows that $\hat{F}_{\mu}(\x)$ is $\left(L+\frac{L_G}{\mu}\right)$-smooth. Starting with the upper bound
\begin{align}
\hat{F}_{\mu_{k}}(\x_{k+1})\leq \hat{F}_{\mu_{k}}(\x_k)+\eta_k\ip{\nff,\z_k-\x_k}+\frac{\eta_k^2}{2}\left(L+\frac{L_G}{\mu_k}\right)D^2
\end{align}
and following the steps same as \eqref{basic smoothness step}-\eqref{basic smoothness step3}, we can obtain the following bound
\begin{align}\label{basic smoothness step33}
\hat{F}_{\mu_k}(\x_{k+1})-f(\x^*)&\leq (1-\eta_k)\hat{F}_{\mu_{k-1}}(\x_k)-f(\x^*)\\&+\frac{\eta^2}{2}D^2\left(L+\frac{L_G}{\mu_k}\right)+\eta_kD\norm{\nff-\y_k}\nonumber\\&+\left( \frac{(1-\eta_k)(\mu_{k-1}-\mu_k)-\eta_k\mu_k}{2}\right)||\lam^{\star}_{\mu_k}(\G(\xi_k)\x_k)||^2.
\end{align}
Now, setting $\eta_k$ and $\mu_k$ such that $
(1-\eta_k)(\mu_{k-1}-\mu_k)-\eta_k\mu_k\leq 0$, 
we get rid of last term of \eqref{basic smoothness step33}. Taking full expectation:
\begin{align*}
\mbE[\hat{F}_{\mu_k}(\x_{k+1})]-f(\x^*)&\leq (1-\eta_k)\mbE[\hat{F}_{\mu_{k-1}}(\x_k)]-f(\x^*)\nonumber\\&+\frac{\eta_k^2}{2}D^2\left(L+\frac{L_G}{\mu_k}\right)+\eta_kD\sqrt{\mbE||\nff-\y_k||^2}.
\end{align*}
Substituting  $\mu_k =\frac{\mu_c}{(k+1)^{1/4}}$, $\eta_k = \tfrac{2}{k+1}$, and $\rho_k \leq \frac{D}{\sqrt{m}(k+1)}$:
\begin{align}\label{smoothness5}
\mbE[\hat{F}_{\mu_k}(\x_{k+1})]-f(\x^\star)&\leq \left(1-\frac{2}{k+1}\right)\mbE[\hat{F}_{\mu_{k-1}}(\x_k)]-f(\x^\star)\nonumber\\&\quad+\frac{2}{(k+1)^2}\left(L+\frac{L_G(k+1)^\frac{1}{4}}{\mu_c}\right)D^2\nonumber\\&\quad+\frac{2D}{(k+1)}\frac{\sqrt{96(\sigma^2+12L^2D^2+9L_G^2\mu_c^{-2}D^2)}}{k^\frac{1}{4}}\nonumber\\&
\leq \left(1-\frac{2}{k+1}\right)\mbE[\hat{F}_{\mu_{k-1}}(\x_k)]-f(\x^\star)\nonumber\\&\quad+\frac{8\sqrt{6}\sigma D+2(35L+31L_G\mu_c^{-1})D^2}{k^{5/4}}
\end{align}
since $k^{5/4} \leq (k+1)^2$. Note that our selection of $\mu_k$ and $\eta_k$ satisfies the assumption $\mu_k \geq \mu_{k-1}(1-\eta_k)$. Finally, application of Lemma \ref{simplelemma} yields the required result. 
\subsection{Proof of Theorem \ref{theorem2}:}\label{proof:theorem2}
\noindent (a) The bound on the optimality gap follows from Lemma \ref{lemma:F-smooth}:
\begin{align}\label{temp22}
\mbE[f(\x_{k+1})]-f(\x^{\star})&\leq  \mbE[f(\x_{k+1})]-f(\x^{\star})+\frac{1}{2\mu_k}\mbE [\cD^2_{\cX(\xi)}(\G(\xi)\x_{k+1})]\nonumber\\& = \mbE[\hat{F}_{\mu_k}(\x_{k+1})]-f(\x^{\star})\nonumber\\&\leq  \frac{8(4\sqrt{6}\sigma D+35LD^2+31L_G\mu_c^{-1}D^2)}{k^{\frac{1}{4}}}
\end{align}

\noindent (b)	We use Jensen's inequality and Corollary \ref{coro:borrow4} to obtain bound on feasibility as	
\begin{align}\label{tem2}
&\mbE[\cD_{\cX(\xi)}(\G(\xi)\x_{k+1})]\leq\sqrt{\mbE[\cD^2_{\cX(\xi)}(\G(\xi)\x_{k+1})]}\\& \leq\sqrt{ 4\mu_k^2\norm{\lam^{\star}}^2+4\mu_k\mbE[\hat{F}_{\mu_k}(\x_{k+1})-f(\x^{\star})]}\nonumber\\&
\leq\frac{2\mu_c }{(k+1)^\frac{1}{4}}\norm{\lam^{\star}}+\frac{18\sqrt{\sigma  D \mu_c}+34D\sqrt{L\mu_c}+32D\sqrt{L_G}}{(k+1)^{\frac{1}{8}}(k)^{\frac{1}{8}}}\nonumber\\&
\leq \frac{2\left(\mu_c \norm{\lam^{\star}}+9\sqrt{\sigma  D \mu_c}+17D\sqrt{L\mu_c}+16D\sqrt{L_G}\right)}{k^{\frac{1}{4}}}.\nonumber
\end{align}	 

\subsection{Proof of Theorem \ref{theorem3}:}\label{proof:theorem3}
\noindent (a) Recall from  \eqref{grad_cons_funct} that
$\nf=\n f(\x_k)+\G^\T\n h_{\mu_k}(\G\x,\cX)$, which allows us to write the second term on the right of \eqref{basic smoothness step} as
\begin{align}\label{abc}
\ip{\nf,\z_k-\x_k}
&=\ip{\n f(\x_k),\z_k-\x_k}+\ip{\G^\T\n h_{\mu_k}(\G\x_k,\cX),\z_k-\x_k}\\
&=\ip{\n f(\x_k)-\v_k,\z_k-\x_k}+\ip{\v_k,\z_k-\x_k}+\ip{\G^\T\n h_{\mu_k}(\G\x_k,\cX),\z_k-\x_k}\label{xa2}\nonumber\\&\leq \ip{\n f(\x_k)-\v_k,\z_k-\x_k}+\ip{\v_k,\x^{\star}-\x_k}+\ip{\G^\T\n h_{\mu_k}(\G\x_k,\cX),\z_k-\x_k},\nonumber
\end{align}
Note that if trimming is not performed at iteration $k$, then we call the LMO to solve the problem 
$\z_k=\argmin_{\z\in \cc}\ip{\z,\v_k}$ while if trimming is done we set $\v_k=\v_{k-1}$ and $\z_k=\z_{k-1}$. Suppose that trimming was not  performed at iteration $k-1$, then we can say that $\z_k=\z_{k-1} = \argmin_{\z\in \cc}\ip{\z,\v_{k-1}} = \argmin_{\z\in \cc}\ip{\z,\v_{k}}$. The same argument can again be applied to conclude that $\z_k = \argmin_{\z\in \cc}\ip{\z,\v_{k}}$ if trimming was not performance at some iteration $k_0 \leq k$. We can assume that trimming is not performed at the first iteration. Hence, irrespective of whether trimming is done or not $\z_k$ is the minimizer of $\ip{\z_k,\v_k}$. Thus, at each iteration  we have $\ip{\z_k,\v_k}\leq \ip{\x^{\star},\v_k}$ using which we obtain \eqref{abc}.
Adding and subtracting $\ip{\n f(\x_k),\x^{\star}-\x_k}$, and using the convexity of $f$, we get
\begin{align}
\ip{\nf,\z_k-\x_k} &\leq \ip{\n f(\x_k)-\v_k,\z_k-\x^{\star}}+\ip{\n f(\x_k),\x^{\star}-\x_k}+\ip{\G^\T\n h_{\mu_k}(\G\x_k,\cX),\z_k-\x_k}\nonumber\\
&\leq \ip{\n f(\x_k)-\v_k,\z_k-\x^{\star}}+f(\x^{\star})- f(\x_k)+\ip{\G^\T\n h_{\mu_k}(\G\x_k,\cX),\z_k-\x_k}
\end{align}
Further, adding and subtracting $\ip{\w_k,\z_k-\x^{\star}}$ we get,
\begin{align}
\ip{\nf,\z_k-\x_k}
&\leq \ip{\n f(\x_k)-\v_k+\w_k-\w_k,\z_k-\x^{\star}}+f(\x^{\star})- f(\x_k)+\ip{\G^\T\n h_{\mu_k}(\G\x_k,\cX),\z_k-\x_k}\nonumber\\
&\leq \ip{\n f(\x_k)-\w_k,\z_k-\x^{\star}}+\ip{\w_k-\v_k,\z_k-\x^{\star}}\nonumber\\
&\quad+f(\x^{\star})- f(\x_k)+\ip{\G^\T\n h_{\mu_k}(\G\x_k,\cX),\z_k-\x_k}\nonumber\\
&\leq \ip{\n f(\x_k)-\y_k,\z_k-\x^{\star}}+\ip{\G^\T\n h_{\mu_k}(\G\x_k,\cX),\x^{\star}-\z_k}\nonumber\\&\quad+\ip{\w_k-\v_k,\z_k-\x^{\star}}+f(\x^{\star})- f(\x_k)\nonumber\\&\quad+\ip{\G^\T\n h_{\mu_k}(\G\x_k,\cX),\z_k-\x_k}\nonumber\\&= \ip{\n f(\x_k)-\y_k,\z_k-\x^{\star}}+\ip{\w_k-\v_k,\z_k-\x^{\star}}\nonumber\\&\quad+f(\x^{\star})- f(\x_k)+\ip{\G^\T\n h_{\mu_k}(\G\x_k,\cX),\x^{\star}-\x_k}\nonumber\\
&\leq \norm{\n f(\x_k)-\y_k}D+\norm{\w_k-\v_k}D+f(\x^{\star})- f(\x_k)\nonumber\\&\quad+\ip{\G^\T\n h_{\mu_k}(\G\x_k,\cX),\x^{\star}-\x_k}\label{ee}\\
&\leq \norm{\n f(\x_k)-\y_k}D+\tau_kD+f(\x^{\star})- f(\x_k)\nonumber\\&\quad+\ip{\G^\T\n h_{\mu_k}(\G\x_k,\cX),\x^{\star}-\x_k}\label{grad1t}
\end{align}
here  \eqref{ee} is obtained using the Cauchy-Schwartz inequality  the compactness assumption on $\cc$ while the last inequality comes from the use of Lemma \ref{cens-error}.

Now using \eqref{borrow1_our} in \eqref{grad1t} we  obtain
\begin{align}\label{half_gradt}
	\ip{\nf,\z_k-\x_k}
	&\leq \|\n f(\x_k)-\y_k\|D+\tau_kD+f(\x^{\star})- f(\x_k)+\one_{\cX}(\G\x^{\star})\nonumber\\	&-h_{\mu_k}(\G\x_k,\cX)-\frac{\mu_k}{2}\|\lam^{\star}_{\mu_k}(\G\x_k)\|^2.
\end{align}
Substituting the bounds obtained in \eqref{half_gradt} and \eqref{temp1} into \eqref{basic smoothness step} and subtracting $f(\x^{\star})$ from both sides, we obtain
\begin{align}\label{basic smoothness step3t}
F_{\mu_k}(\x_{k+1})-f(\x^{\star})
&\leq (1-\eta_k)\left(F_{\mu_{k-1}}(\x_{k})-f(\x^{\star})\right)+\frac{\eta_k^2}{2}\left(L+\frac{L_G}{\mu_k}\right)D^2\nonumber\\&\quad+\eta_k\|\n f(\x_k)-\y_k\|D+\eta_{k}\tau_kD+ \frac{1}{2}\left((1-\eta_k)\mu_{k-1}-\mu_k\right)\|\lam^{\star}_{\mu_k}(\G\x_k)\|^2,
\end{align}
here we used the fact that $\one_{\cX}(\G\x^{\star})=0$. Now, since we have assumed that $\mu_k \geq \mu_{k-1}(1-\eta_k)$, the last term in \eqref{basic smoothness step3t} is non-positive and can be dropped, yielding  
\begin{align}\label{basic smoothness step3t2}
F_{\mu_k}(\x_{k+1})-f(\x^{\star})
&\leq (1-\eta_k)\left(F_{\mu_{k-1}}(\x_{k})-f(\x^{\star})\right)+\frac{\eta_k^2}{2}\left(L+\frac{L_G}{\mu_k}\right)D^2\nonumber\\&\quad+\eta_k\|\n f(\x_k)-\y_k\|D+\eta_{k}\tau_kD.
\end{align}
Now, taking expectation on \eqref{basic smoothness step3t2}, we have
\begin{align}\label{flw2}
\Ex{F_{\mu_k}(\x_{k+1})-f(\x^{\star})}
&\leq (1-\eta_k)\Ex{F_{\mu_{k-1}}(\x_{k})-f(\x^{\star})}+\frac{\eta_k^2}{2}\left(L+\frac{L_G}{\mu_k}\right)D^2\nonumber\\&\quad+\eta_kD\sqrt{\Ex{\norm{\n f(\x_k)-\y_k}^2}}+\eta_{k}\tau_kD
\end{align} 
Substituting the result of Lemma \ref{track1} and that of the various parameters, we obtain
\begin{align}
\Ex{F_{\mu_k}(\x_{k+1})-f(\x^{\star})} &\leq  (1-\frac{2}{k+1})\Ex{F_{\mu_{k-1}}(\x_{k})-f(\x^{\star})} +\frac{2\tau_0D}{(k+1)^{3/2}}\nonumber\\
&+\frac{1}{(k+1)^2}\left(L+\frac{L_G\sqrt{k}}{\mu_c}\right)D^2+\frac{8D}{k+1}\sqrt{\frac{3\sigma^2+25L^2D^2}{k}}
\end{align}
Note that for this choice of parameters, $\frac{\mu_k}{\mu_{k-1}} = \sqrt{\frac{k-1}{k}}\geq \frac{k-1}{k}=1-\eta_k$, so the condition required for \eqref{mainlem} is satisfied. The bound on the right can therefore be simplified as 
\begin{align}\label{ss1}
&\Ex{F_{\mu_k}(\x_{k+1})-f(\x^{\star})}\\&\leq  (1-\frac{2}{k+1})\Ex{F_{\mu_{k-1}}(\x_{k})-f(\x^{\star})} +  \frac{8\sqrt{3}\sigma D+(41L+L_G\mu_c^{-1})D^2+2\tau_0D}{k^{3/2}}
\end{align}
Finally, application of Lemma \ref{simplelemma} yields the required result.

The bound on the optimality gap follows from \eqref{ss1}:
\begin{align}\label{temp22t}
&\mbE[f(\x_{k+1})]-f(\x^{\star})\leq \frac{4(8\sqrt{3}\sigma D+(41L+L_G\mu_c^{-1})D^2+2\tau_0D)}{\sqrt{k}}.
\end{align}
For constraint violation, following \eqref{temp3} and using \eqref{temp22t} we have 
\begin{align}\label{temp3t}
&-\norm{\lam^{\star}}\mbE[\cD_{\cX}(\G\x_{k+1})]+\frac{1}{2\mu_k}(\mbE [\cD_{\cX}(\G\x_{k+1})])^2 \leq \frac{4(8\sqrt{3}\sigma D+(41L+L_G\mu_c^{-1})D^2+2\tau_0D)}{\sqrt{k}}
\end{align}
which is a quadratic inequality in $\Ex{\cD_\cX(\G\x_{k+1})}$, and can be solved to yield the required result:
\begin{align}\label{cons_fes2}
	\mbE[\cD_{\cX}(\G\x_{k+1})]&\leq \frac{1}{\sqrt{k}}\bigg(2\mu_c\norm{\lam^{\star}}+\big(8(8\sqrt{3}\sigma D\mu_c+(41L\mu_c+L_G)D^2+2\tau_0D)\big)^{1/2}\bigg) \nonumber\\
&\leq \frac{1}{\sqrt{k}}\bigg(2\mu_c\norm{\lam^{\star}}+ 11\sqrt{\sigma D \mu_c} + 19D\sqrt{L\mu_c}+3D\sqrt{L_G}+4\sqrt{\tau_0D}\bigg)
\end{align}

\noindent (b) Following the steps same as \eqref{abc}-\eqref{flw2}, we can obtain the following bound for \ms
\begin{align}\label{basic smoothness step33t}
\hat{F}_{\mu_k}(\x_{k+1})-f(\x^*)&\leq (1-\eta_k)\hat{F}_{\mu_{k-1}}(\x_k)-f(\x^*)+\frac{\eta^2}{2}D^2\left(L+\frac{L_G}{\mu_k}\right)\nonumber\\&\quad+\eta_kD\norm{\nff-\y_k}+\eta_k\tau_kD\nonumber\\&\quad+\left( \frac{(1-\eta_k)(\mu_{k-1}-\mu_k)-\eta_k\mu_k}{2}\right)||\lam^{\star}_{\mu_k}(\G(\xi_k)\x_k)||^2.
\end{align}
Now, setting $\eta_k$ and $\mu_k$ such that $
(1-\eta_k)(\mu_{k-1}-\mu_k)-\eta_k\mu_k\leq 0$, 
we get rid of last term of \eqref{basic smoothness step33t}. Taking full expectation:
\begin{align}\label{basic smoothness4t}
\mbE[\hat{F}_{\mu_k}(\x_{k+1})]-f(\x^*)&\leq (1-\eta_k)\mbE[\hat{F}_{\mu_{k-1}}(\x_k)]-f(\x^*)+\frac{\eta_k^2}{2}D^2\left(L+\frac{L_G}{\mu_k}\right)\nonumber\\&\quad+\eta_kD\sqrt{\mbE||\nff-\y_k||^2}+\eta_k\tau_kD.
\end{align}
Substituting  $\mu_k =\frac{\mu_c}{(k+1)^{1/4}}$, $\tau_k=\frac{\tau_0}{(k+1)^{1/4}}$, $\eta_k = \tfrac{2}{k+1}$, and $\rho_k \leq \frac{D}{\sqrt{m}(k+1)}$, we get
\begin{align}\label{smoothness5t}
\mbE[\hat{F}_{\mu_k}(\x_{k+1})]-f(\x^\star)
&\leq \left(1-\frac{2}{k+1}\right)\mbE[\hat{F}_{\mu_{k-1}}(\x_k)]-f(\x^\star)\nonumber\\&\quad+\frac{2}{(k+1)^2}\left(L+\frac{L_G(k+1)^\frac{1}{4}}{\mu_c}\right)D^2+\frac{2\tau_0D}{(k+1)^{5/4}}\nonumber\\&\quad+\frac{2D}{(k+1)}\frac{\sqrt{96(\sigma^2+12L^2D^2+9L_G^2\mu_c^{-2}D^2)}}{k^\frac{1}{4}}\nonumber\\&
\leq \left(1-\frac{2}{k+1}\right)\mbE[\hat{F}_{\mu_{k-1}}(\x_k)]-f(\x^\star)\nonumber\\&\quad+\frac{8\sqrt{6}\sigma D+2(35L+31L_G\mu_c^{-1})D^2+2\tau_0D}{k^{5/4}}
\end{align}
Finally, application of Lemma \ref{simplelemma} yields the required results.

For constraint violation, 	following proof of Theorem \ref{theorem2}(b), we obtain the required bound as	
\begin{align}\label{tem2t}
&\mbE[\cD_{\cX(\xi)}(\G(\xi)\x_{k+1})]\\&
\leq\frac{2\mu_c \norm{\lam^{\star}}}{(k+1)^\frac{1}{4}}+\frac{18\sqrt{\sigma  D \mu_c}+34D\sqrt{L\mu_c}+32D\sqrt{L_G}+6\tau_0D}{(k+1)^{\frac{1}{8}}(k)^{\frac{1}{8}}}\nonumber\\&
\leq \frac{2\left(\mu_c \norm{\lam^{\star}}+9\sqrt{\sigma  D \mu_c}+17D\sqrt{L\mu_c}+16D\sqrt{L_G}+3\sqrt{\tau_0D}\right)}{k^{\frac{1}{4}}}.\nonumber
\end{align}

   \end{document}